\documentclass[leqno,12pt]{article}
\usepackage{amsmath}
\usepackage{amssymb}
\usepackage{amsthm}
\usepackage{amscd}
\usepackage{amsxtra}
\usepackage{verbatim}
\usepackage{xcolor}
\usepackage{color}
\usepackage{enumerate}
\usepackage{mathrsfs}
\usepackage[all]{xy} 

\usepackage{array,arydshln}
\usepackage{bm}

\allowdisplaybreaks

\numberwithin{equation}{section}

\definecolor{dblue}{rgb}{0,0,0.45}
\definecolor{red}{rgb}{0.7,0,0}

 \RequirePackage{geometry}
\newtheorem{theorem}{Theorem}[section]
\newtheorem{lemma}[theorem]{Lemma}
\newtheorem*{lemma*}{Lemma}

\newtheorem{proposition}[theorem]{Proposition}

\theoremstyle{definition}

\newtheorem{remark}[theorem]{Remark}
\newtheorem{definition}[theorem]{Definition}
\newtheorem{example}[theorem]{Example}

\theoremstyle{remark}

\newcommand{\E}{{\mathbb E}}
\newcommand{\N}{{\mathbb N}}
\newcommand{\R}{{\mathbb R}}

\newcommand{\cC}{{\mathcal C}}

\newcommand{\cF}{{\mathcal F}}

\newcommand{\cM}{{\mathcal M}}
\newcommand{\cN}{{\mathcal N}}

\newcommand{\cQ}{{\mathcal Q}}
\newcommand{\cS}{{\mathcal S}}

\newcommand{\cV}{{\mathcal V}}
\newcommand{\cW}{{\mathcal W}}

\newcommand{\la}{\langle}
\newcommand{\ra}{\rangle}
\newcommand{\nn}{\nonumber}
\newcommand{\ve}{\varepsilon}

\newcommand{\vertiii}[1]{{\left\vert\kern-0.25ex\left\vert\kern-0.25ex\left\vert #1 
    \right\vert\kern-0.25ex\right\vert\kern-0.25ex\right\vert}}

\date{}
\begin{document}

\pagestyle{myheadings}
\markright{Averaging for slow-fast RDEs} 

\title{Averaging principle for slow-fast systems of
rough differential equations 
via controlled paths
\footnote{ 
2020 \textit{Mathematics Subject Classification}.
Primary 60L90; Secondary 70K65, 70K70, 60F99.			
}
\footnote{ 
\textit{Key words and phrases}. 
Slow-fast system, averaging principle, 
rough path theory.
}
\footnote{ 
The author thanks Professor Pei Bin for helpful comments.
The author is supported by 
JSPS KAKENHI Grant No. 20H01807.
}
}
\author{   Yuzuru \textsc{Inahama} }
\maketitle

\begin{abstract}
\noindent
In this paper we prove the strong averaging principle
for a slow-fast system of rough differential equations.
The slow and the fast component of the system are
driven by a rather general random rough path 
and Brownian rough path, respectively.
These two driving noises are assumed to be independent.
A prominent example of the driver of the slow component is 
fractional Brownian rough path with Hurst parameter 
between $1/3$ and $1/2$.  
We work in the framework of  controlled path theory, 
which is one of the most widely-used frameworks
in rough path theory.
To prove our main theorem, we carry out Khas'minski\u{\i}'s 
time-discretizing method in this framework. 
%
%
\end{abstract}


\section{Introduction}

We study the averaging principle 
for slow-fast systems of stochastic differential equations (SDEs).
Although a few different limit theorems are called by the same name, 
the one we focus on in this paper is most typically formulated in the following way.

Let $(w_t)$ and $(b_t)$ be two independent standard
(finite-dimensional) Brownian motions (BMs).
A slow-fast system of (finite-dimensional) SDEs are given  by
\begin{equation} \nn%
\left\{
\begin{array}{lll}
X^{\ve}_t &=& x_0 + 
\int_0^{t}
 f(X^\varepsilon_s, Y^\varepsilon_s) ds
 +
\int_0^{t} \sigma (X^\varepsilon_s, Y^\varepsilon_s) db_s,
 \\
Y^{\ve}_t &=& y_0 + 
\varepsilon^{-1}\int_0^{t}
 g(X^\varepsilon_s, Y^\varepsilon_s) ds
 +
\varepsilon^{-1/2}\int_0^{t}
 h(X^\varepsilon_s, Y^\varepsilon_s) dw_s,
\end{array}
\right.
\end{equation}
where $0<\ve \ll 1$ is a small parameter.
The processes
$X^\varepsilon$ and $Y^\varepsilon$ are called 
the slow component and the fast component, respectively.
Suitable conditions are imposed on $g$ and $h$ so that 
the following frozen SDE satisfies certain ergodicity
for every $x$.
\[
Y^{x,y}_t =y+ \int_0^t g(x, Y^{x,y}_t)dt 
     +  \int_0^t  h(x, Y^{x,y}_t) dw_t,
\]
An associated unique invariant probability measure 
is denoted by $\mu^x$. Set 
$\bar{f} (x) =\int   f(x,y) \mu^x (dy)$ 
and $\bar{\sigma} (x)= \int   \sigma (x,y) \mu^x (dy)$
and consider the following 
averaged SDE:
\[
\bar{X}_t = x_0 + \int_0^{t}
 \bar{f} (\bar{X}_s) ds+
\int_0^{t} \bar\sigma (\bar{X}_s) db_s.
\]
The averaging principle of this type, 
which was initiated by Khas'minski\u{\i} \cite{khas},
 claims that 
$X^\varepsilon$ converges to $\bar{X}$ in an appropriate sense 
as $\ve \searrow 0$.
Even though the history is old and many papers have been written
(see \cite{fw, gkk, gol, gl, liu, lrsx, rsx1, ve, xlm} for example),
this research topic seems still quite active.  

It should also be recalled that
this averaging principle was generalized to  various kind of 
stochastic equations.
Examples include jump-type SDEs \cite{giv, liu2, sxx, xm, zfwl}, 
distribution-dependent SDEs  \cite{rsx2, xllm}, 
manifold-valued SDEs \cite{lixm},
functional-type SDEs
such as SDEs with delay \cite{bsyy, hy, wy1, wy2} among others.
(After a pioneering work \cite{cf},
the case of stochastic partial differential equations 
has also been studied extensively.
But, we do not discuss it in this paper.)


In all the preceding works mentioned above, the driving noises
are (semi)martingales. 
One naturally wonders what happens to the averaging principle
when the driving noise does not have (semi)martingale property.
A prominent example of such noises is fractional Brownian motion (fBM).
This research direction is fairly new and 
there are not many papers  at the moment of writing.

When $(b_t)$ is fBM with Hurst parameter $H \in (1/2, 1)$
and $(w_t)$ is the usual BM,
the averaging principle for slow-fasts systems like 
 \eqref{def.SFeq} was recently proved in \cite{hl, pix1, wxy, hxpw}
under various settings.
In these works, the integral $db_s$ is understood as a Young
(i.e.  a generalized Riemann-Stieltjes) integral.
Concerning this kind of slow-fast systems driven by 
fBM with $H \in (1/2, 1)$ and BM,
a few related problems have  
already been studied in \cite{bgs1, bgs2, bgs3}.
Though it looks quite difficult to study the case where
$(w_t)$ is also fBM, a recent preprint \cite{ls} made
a first attempt in that direction.

When $(b_t)$ is fBM with Hurst parameter $H <1/2$,
the problem becomes quite difficult, 
mainly because neither Young integration nor It\^o integration 
is available.
In \cite{pix2} the authors proved the averaging problem 
when $(b_t)$ is fBM with $1/3 <H \le 1/2$ by using rough path (RP) theory.
In the above mentioned paper, $db_s$ is actually understood
as a RP integral along a fractional Brownian RP.

The present paper is a continuation of \cite{pix2}
and generalizes its main theorem to a considerable extent
(see Remark \ref{rem.differ} below for details).
The framework of RP theory adopted in this paper and 
that in \cite{pix2} are different.
In this paper the controlled path (CP) theory is exclusively used,
while in \cite{pix2} a fractional calculus approach to 
RPs is mainly used.


The organization of this paper is as follows:
In Section 2, we introduce assumptions on the coefficients
and on the driving RP and then state our main theorem. 
A comparison with preceding works and examples
are also given. 
In Section 3, following the textbook \cite{fh},  we review CP theory.
We also slightly generalize
well-known results on the well-posedness of a
rough differential equation (RDE)
so that we can deal with the slow component of 
the slow-fast system \eqref{def.SFeq} of RDEs 
and the averaged RDE  \eqref{def.avRDE}. 
 All arguments in this section are deterministic.
The first half of Section 4 is devoted to an rigorous introduction of 
the slow-fast system of RDEs in a deterministic way.
An It\^o-Stratonovich correction formula 
for the fast component is also given.
The second half is devoted to showing non-explosion of 
the slow-fast system of RDEs when the driving RP is random. 
Some moment estimates of the solution,
which will be used in the proof of our main theorem,
 are also obtained.
In Section 5 we prove our main theorem by carrying out 
Khas'minski\u{\i}'s 
time-discretizing method in the framework of  CP theory.
In appendices some important known results are summarized.
The most important ones among them 
are basic facts on the frozen SDE
associated with (the fast component of) the slow-fast system.



Before closing Introduction, 
 we introduce the notation which will be used throughout the paper.
We denote the set of integers by $\N =\{0,1,2, \ldots\}$.
Let $T >0$ be arbitrary and we work on the time interval
$[0,T]$ unless otherwise specified.
For a subinterval $[a,b]\subset [0,T]$, 
we write $\triangle_{[a,b]} =\{ (s,t)\in \R^2 \mid  a\le s\le t \le b\}$.
When $[a,b]= [0,T]$, we simply write $\triangle_{T}$ for this set.
(The time horizon $T$ and the starting point $(x_0, y_0)$
are arbitrary but fixed. We will not keep track of the dependence 
on $T, x_0, y_0$.)

Below, $\cV$ and $\cW$ are (real) Banach spaces.
The set of bounded linear maps from $\cV$ to $\cW$ is denoted 
by $L(\cV, \cW)$.
When $\cV= \R^n$ and $\cW = \R^m$, $L(\cV, \cW)$
coincides with the set of all real $m \times n$ matrices
and is equipped with Hilbert-Schmidt norm instead of the operator norm.

\begin{itemize} 
\item
The set of all continuous paths $\varphi\colon [a,b] \to\cV$
is denoted by $\cC ([a,b], \cV)$. 
With the usual sup-norm $\|\varphi\|_{\infty, [a,b]}$ on the $[a,b]$-interval,
$\cC ([a,b], \cV)$ is a Banach space.
The difference of $\varphi$ is denoted by $\varphi^1$,
that is, $\varphi^1_{s,t} := \varphi_t - \varphi_s$ for $(s,t)\in \triangle_{[a,b]}$.

\item
Let $0< \gamma \le 1$.
For a path $\varphi \colon [a,b] \to \cV$,
the $\gamma$-H\"older seminorm is defined by 
\[
\|\varphi\|_{\gamma,[a, b]} :=\sup _{a \le s<t \le b} 
\frac{\left|\varphi_{t}-\varphi_{s}\right|_{\cV}}{(t-s)^{\gamma}}.
\]
If the right hand side is finite, we say $\varphi$ is 
$\gamma$-H\"older continuous on $[a,b]$.
The space of all $\gamma$-H\"older continuous paths on $[a,b]$ is denoted by 
$\cC^\gamma ([a,b], \cV)$.
The Banach norm on this space is
 $|\varphi_a|_{\cV}+\|\varphi\|_{\gamma,[a, b]}$.

\item
Let $0< \gamma \le 1$.
For a continuous map $\eta \colon \triangle_{[a,b]} \to \cV$, we set 
\[
\|\eta\|_{\gamma,[a, b]} :=\sup _{a\le s<t\le b} 
\frac{\left|\eta_{s,t}\right|_{\cV}}{(t-s)^{\gamma}}.
\]
If this is finite, then $\eta$ vanishes on the diagonal.
The set of all such $\eta$ with $\|\eta\|_{\gamma,[a, b]}<\infty$
is denoted by $\cC_2^\gamma ([a,b], \cV)$,
which is a Banach space with $\|\eta\|_{\gamma,[a, b]}$.

\item
When $[a,b] =[0,T]$, we write
$\cC (\cV)$, $\cC^\gamma (\cV)$, $\cC_2^\gamma (\cV)$
for these spaces 
and $\|\cdot\|_{\infty}$, $\|\cdot\|_{\gamma}$, $\|\cdot\|_{\gamma}$
for the corresponding (semi)norms for simplicity of notation.

\item
Let $U$ be an open set of $\R^m$.
For $k \in \N$,  $C^k (U, \R^n)$ stands for the set of 
$C^k$-functions from $U$ to $\R^n$.
(When $k=0$, we simply write $C (U, \R^n)$ 
instead of $C^0 (U, \R^n)$.)
The set of bounded $C^k$-functions $f \colon U\to \R^n$
whose derivatives up to order $k$ are all bounded 
is denoted by $C_b^k (U, \R^n)$, which is a Banach space with the norm
$\| f\|_{C_b^k } := \sum_{i=0}^k \|\nabla^i f\|_{\infty}$.
(Here, $ \|\cdot\|_{\infty}$ stands for the usual sup-norm on $U$.)
\end{itemize}

\section{Assumptions and main result}

In this section we first introduce natural assumptions 
on the coefficients and driving random RP
of the following slow-fast system.
Then, we state our main theorem.

Our slow-fast system of RDEs, defined on 
the time interval $[0,T]$, is given by 
\begin{equation} \label{def.SFeq}
\left\{
\begin{array}{lll}
X^{\ve}_t &=& x_0 + 
\int_0^{t}
 f(X^\varepsilon_s, Y^\varepsilon_s) ds
 +
\int_0^{t} \sigma (X^\varepsilon_s) dB_s,
 \\
Y^{\ve}_t &=& y_0 + 
\varepsilon^{-1}\int_0^{t}
 g(X^\varepsilon_s, Y^\varepsilon_s) ds
 +
\varepsilon^{-1/2}\int_0^{t}
 h(X^\varepsilon_s, Y^\varepsilon_s) dW_s.
\end{array}
\right.
\end{equation}
Here, $0 <\ve \le 1$ is a small constant
and (the first level path of)
$(X^{\ve}, Y^{\ve})$ takes values in $\R^m \times \R^n$. 
The starting point $(x_0, y_0)$ is always deterministic and arbitrary.
At the first stage, \eqref{def.SFeq} is a deterministic system of RDEs 
driven by an ($d+e$)-dimensional RP which is denoted by $(B, W)$. 
A precise definition of the system \eqref{def.SFeq}
will be given in Section 4.

When we consider this slow-fast system of RDEs,
we always impose the following conditions:
\begin{itemize}
\item
$\sigma\in C^3 (\R^m, L(\R^d, \R^m))$ and
$h \in C^3 (\R^m \times \R^n, L(\R^e, \R^n))$,
\item
$f \in C (\R^m \times \R^n, \R^m)$ 
is locally Lipschitz continuous
and so is $g \in C(\R^m \times \R^n, \R^n)$.
\end{itemize}
These guarantee that  
the slow-fast system of RDEs always has a unique 
local solution. (This fact is well-known. See also Remark \ref{rem.loc.sol} below.)
Since we show the strong version of the averaging principle in this work, 
we assume that $\sigma$  depends only on the slow variable.

To formulate our main result, we introduce several more 
assumptions on these coefficients.

\medskip
\noindent
${\bf (H1)}$~$\sigma$ is of $C^3_b$.

\medskip
\noindent
${\bf (H2)}$~$f$ is bounded and globally Lipschitz continuous.

\medskip
\noindent
${\bf (H3)}$~$h$ is globally Lipschitz continuous.

\medskip
\noindent
${\bf (H4)}$~There exist constants $\eta_1 \ge 0$ and $C>0$
such that, 
for all $x \in \R^m$ and $y \in \R^n$, 
\[
|g(x,y)| \le C( |x|^{\eta_1} + |y|^{\eta_1} +1).
\]

We set the following condition for $r\ge 0$:

\medskip
\noindent
${\bf (H5)}_r$~There exist constants $\eta_2 \ge 0$
and $C >0$  such that,
for all $x_1, x_2 \in \R^m$ and $y \in \R^n$, 
\[
|g(x_1,y) -g(x_2, y)| \le C  |x_1 -x_2|   
(1+  |x_1|^{\eta_2}+|x_2|^{\eta_2} +|y|^r).
\]

\medskip

Next, we set the following condition for $q \ge 2$.
If $q' >q$, then ${\bf (H6)}_{q'}$ obviously implies ${\bf (H6)}_q$
without changing the constants $\gamma_1, \eta_3, C$.

\medskip
\noindent
${\bf (H6)}_q$~There exist constants $\gamma_1 >0$,
$C >0$ and $\eta_3 \ge 0$ such that,
for all $x \in \R^m$ and $y \in \R^n$, 
\[
2\langle y, g (x,y)\rangle + (q-1) |h (x,y)|^2
\le -\gamma_1 |y|^2 + C(|x|^{\eta_3} +1).
\]

\medskip
\noindent
${\bf (H7)}$~There exists a constant $\gamma_2 >0$ such that,
for all $x \in \R^m$ and $y_1, y_2 \in \R^n$, 
\[
2\langle y_1- y_2, g(x, y_1)-g (x,y_2)\rangle
+  | h(x,y_1)- h (x,y_2) |^{2}
 \le  -\gamma_2 |y_1-y_2|^{2}.
 \]

\medskip

Let $T \in (0, \infty)$ and $\tfrac13 <\alpha_0 \le \tfrac12$. 
Let 
$(\Omega, \mathcal{F}, {\mathbb P};
 \{\cF_t\}_{0\le t\le T})$ be a filtered probability space
satisfying  the usual condition. 
On this probability space, 
the following two independent random variables 
$w$ and $B=(B^1, B^2)$  are defined.
The former, 
$w =(w_t)_{0\le t \le T}$, is a standard $e$-dimensional 
$\{\cF_t\}$-BM. 
The It\^o RP lift of $w$ is denoted by $W=(W^1, W^2)$.
The latter, $B=\{(B^1_{s,t}, B^2_{s,t})\}_{0\le s\le t\le T}$, is an 
$\Omega_{\alpha} (\R^d)$-valued 
random variable (i.e., random RP) for every $\alpha \in (1/3,\alpha_0)$.
Here, $\Omega_{\alpha} (\R^d)$ is the space of 
$\alpha$-H\"older RPs over $\R^d$.
We assume that $(B^1_{s,t}, B^2_{s,t})$ is $\cF_t$-measurable 
for every $(s,t)$ with $0\le s \le t \le T$.
Note that $B$ need not be (weakly) geometric.

We assume the following condition on  the integrability of $B$.
Below, $\vertiii{B}_{\alpha}:= 
\|B^1\|_{\alpha}+ \|B^2\|_{2\alpha}^{1/2}$ denotes the $\alpha$-H\"older homogeneous RP norm over the time interval $[0,T]$.

\medskip 
\noindent 
${\bf (A)}$~  For every $\alpha \in (1/3,\alpha_0)$ and 
$p \in [1,\infty)$,  we have $\E [\vertiii{B}_{\alpha}^p  ] <\infty$.

\medskip 

\noindent
Under this assumption, the mixed random RP $(B,W)$
and the slow-fast system \eqref{def.SFeq} of RDEs
driven by it 
can be defined in a natural way
(see Section \ref{sec.SF} for precise definitions).
We will show the averaging principle for \eqref{def.SFeq}
when it is driven by this random RP.

Before we provide our main theorem, we introduce 
the frozen SDE and the averaged RDE 
associated with the slow-fast system 
\eqref{def.SFeq} in the usual way.
The frozen SDE is given as follows:
\[
Y^{x,y}_t =y+ \int_0^t g(x, Y^{x,y}_t)dt 
     +  \int_0^t  h(x, Y^{x,y}_t) d^{{\rm I}}w_t,
\]
Here, $(x,y) \in \R^m \times \R^n$ are deterministic and arbitrary
and $d^{{\rm I}}w_t$ stands for the standard It\^o integral
with respect to a standard $e$-dimensional BM $(w_t)$.
(A more precise definition will be given in Eq. \eqref{def.frozenSDE}.)
Under suitable conditions, the Markov semigroup $(P^x_t)_{t\ge 0}$
defined by $P^x_t \varphi (y)= {\mathbb E} [\varphi (Y^{x,y}_t) ]$
has a unique invariant probability measure,
which is denoted by $\mu^{x}$ 
(see Lemma \ref{lem.lrsx3.8} below for details).
It should be noted that we are only interested in the law of $Y^{x,y}$
and hence any realization of BM will do.

Define the averaged drift by 
$\bar{f} (x) =\int_{\R^n}   f(x,y) \mu^x (dy)$ for $x\in \R^m$.
The averaged RDE is given as follows:
\begin{equation} \label{def.avRDE}
\bar{X}_t = x_0 + 
\int_0^{t}
 \bar{f} (\bar{X}_s) ds
 +
\int_0^{t} \sigma (\bar{X}_s) dB_s
\end{equation}
Here, $x_0 \in \R^m$ is the same as in \eqref{def.SFeq}.
It will be shown that under suitable conditions, 
this RDE has a unique global solution
(see Propositions \ref{prop.0702} and \ref{prop.av.drft} below for details).

Now we are in a position to state our main result.
It claims that (the first level path of) the slow component 
of the slow-fast system \eqref{def.SFeq} of RDEs
converges as $\ve\searrow 0$ to (the first level path of) the 
averaged RDE \eqref{def.avRDE} in $L^p$-sense.
This generalizes the main result of \cite{pix2}. 
Here, $\|\cdot\|_\beta$ stands for the $\beta$-H\"older (semi)norm of a usual path over the time interval $[0,T]$.
\begin{theorem} \label{thm.main}
Assume ${\bf (A)}$, ${\bf (H1)}$--${\bf (H4)}$, 
${\bf (H5)}_r$, ${\bf (H6)}_q$ and ${\bf (H7)}$
for some $q \ge 2$ and $r \ge 0$ such that $q> 2r$.
Then, for every $p\in [1,\infty)$ and $\beta \in (\tfrac13, \alpha_0)$, 
we have
\[
\lim_{\ve \searrow 0} \E [\| X^\ve - \bar{X} \|_\beta^p] =0.
\]
\end{theorem}

\begin{remark} \label{rem.differ}
(1)~In this paper Brownian RP which drives the fast component 
of the slow-fast RDEs is of It\^o type, while in \cite{pix2} it is  
of  Stratonovich type.
However, as will be explained in Lemma \ref{lem.det.IS} and
Remark \ref{rem.IS.corr}, one can easily switch 
between the two  formulations
by adding or subtracting a standard correction term 
(even at the deterministic level).
So, this is just a superficial difference.
\\
\noindent
(2)~
Our main theorem above is stronger than \cite[Theorem 1.2]{pix2}.
The main differences are as follows:
\begin{itemize}
\item
In \cite{pix2}, $B$ is fractional Brownian RP 
with $(\tfrac13, \tfrac12]$. In this paper, $B$ is a much more 
general random RP.
(See also Example \ref{exmpl.B}.)
\item
The conditions on the coefficients $\sigma, h, f, g$ are relaxed
in this paper. (See also Example \ref{exmpl.coeff}.)
\item
In \cite[Theorem 1.2]{pix2}, $L^1$-convergence of 
$\| X^\ve - \bar{X} \|_\infty$ was proved. 
In this paper, $L^p$-convergence of 
$\| X^\ve - \bar{X} \|_\beta$ is proved for every $1\le p< \infty$. 
\end{itemize}
Also, as was mentioned earlier, 
the framework of RDE theory used in this paper is different from 
that in \cite{pix2}.
\end{remark}

\begin{remark} \label{rem.probsp}
The law of $X^\ve - \bar{X}$ is uniquely determined by the law of 
$B=(B^1, B^2)$
and the $e$-dimensional Wiener measure.
So, the choice of a filtered probability space that carries 
$B$ and $w$ does not matter. 

Verifying the existence of such a filtered probability space
is not difficult.
First, a complete probability space $(\Omega, \mathcal{F}, {\mathbb P})$
that supports independent $B$ and $w$ clearly exists.
Then, we set $\cF_t = \sigma \{ w_u\mid 0\le u \le t\}\vee \sigma \{B\} \vee \cN$, $0\le t \le T$, for example.
Here, $\cN$ is the collection of ${\mathbb P}$-zero sets.
Then, $w$ is an $\{\cF_t\}$-BM.
The standard proof of  the right continuity of 
Brownian filtration still works for $\{\cF_t\}$ after trivial modifications.
Thus, $\{\cF_t\}$ satisfies the usual condition and 
$(\Omega, \mathcal{F}, {\mathbb P};
 \{\cF_t\}_{0\le t\le T})$ is a desired example.
 \end{remark}

\begin{example} \label{exmpl.B}
We provide some examples of $B =(B^1, B^2)$ 
satisfying Assumption ${\bf (A)}$ in Theorem \ref{thm.main}.
\begin{enumerate}
\item[(1)]
A deterministic RP $B \in \Omega_{\alpha_0} (\R^d)$ with
$\tfrac13 <\alpha_0 \le \tfrac12$. 
(In this case, we can actually take $\beta =\alpha_0$ 
in Theorem \ref{thm.main}.)
\item[(2)]
Fractional Brownian RP $B^H$ with Hurst parameter 
$\tfrac13 <H \le \tfrac12$. In this case $\alpha_0 =H$.
Note that when $H=1/2$, $B^H$ is Brownian RP of 
Stratonovich type.
\item[(3)]
Brownian RP $\tilde{B}$ of It\^o type. In this case $\alpha_0 =1/2$.
\item[(4)] Mixture of fractional Brownian RP and 
Brownian RP of It\^o type.
Let $B^H$ and $\tilde{B}$ be as in (2) and (3) above, respectively
and assume they are independent.
Then, the mixed RP $(B^H, \tilde{B})$ becomes an example with 
$\alpha_0 =H$.
(The precise definition of mixture will be given in Definition \ref{def.0520}.)
It should also be noted that mixture of 
$(B^H, \tilde{B})$ and $W$ equals mixture of $B^H$ and 
$(\tilde{B},W)$. Here, $W$ is another Brownian RP of 
It\^o type and $\{B^H, \tilde{B}, W\}$ are assumed to be independent.
\item[(5)] In \cite[Chapter 15]{fv}  and \cite[Chapters 10--11]{fh},
there are examples of 
Gaussian RP whose H\"older regularity is between $1/3$ and $1/2$.
Since they satisfy an integrability theorem of Fernique type, their RP norms 
have moments of all orders. So, these Gaussian RP are nice examples.
\end{enumerate}
\end{example}

A slow-fast system of usual SDEs basically corresponds to (3).
The driving process of a slow-fast system in \cite{pix2} is that of (2).
In \cite{wxy}, the slow component is driven by 
BM and fBM with Hurst parameter $H \in (\tfrac12,1)$.
So, the case (4) above can be viewed as a generalization 
of \cite{wxy} to a RP setting.

The conditions on the coefficients $\sigma, h, f, g$ in Theorem \ref{thm.main} are strictly weaker than those in the main theorem 
of the preceding work \cite[Theorem 1.2]{pix2}.
For example, the following ``superlinear" example  
satisfies the conditions of the former, 
but does not satisfy those of the latter.
\begin{example} \label{exmpl.coeff}
Let $\kappa >0$ and
$\lambda, \phi \colon  \R^m \to [\kappa,\infty)$ be 
$C^1$-functions such that $\nabla \lambda$ and 
$\nabla \phi$ are of at most polynomial growth in $x$.
Set 
\[
g (x, y)=  -\lambda (x) y|y|^2 - \phi (x) y,  \qquad 
x \in \R^m, y \in \R^n.
\]
Let $\sigma, h, f$ satisfy ${\bf (H1)}$--${\bf (H3)}$.
Then, ${\bf (H4)}$, ${\bf (H5)}_r$ for some $r \ge 0$,
and ${\bf (H6)}_q$ for every $q \ge 2$ are all satisfied.
If we assume in addition that $\|\nabla_y h \|_\infty <2\kappa$,
then ${\bf (H7)}$ is also satisfied.
Here, $\nabla_y h$ is the partial gradient of $h$ in the $y$-variable.
Therefore, these $\sigma, h, f, g$ satisfy the condition of 
 Theorem \ref{thm.main}.
  \end{example}

\section{Review of controlled path theory}

This section is devoted to recalling CP theory.
We basically follow the exposition in \cite{fh}.
However, we slightly generalize the setting and improve some of the results for our purpose.
Note that everything in this section is deterministic.
Throughout this section,
 $\cV$ and $\cW$ are Euclidean spaces and we assume
 $\tfrac13 <\alpha \le \tfrac12$.

\subsection{Rough paths and controlled paths}
First,  we recall the definition of 
$\alpha$-H\"older rough path ($\alpha$-RP or RP).
A continuous map 
$X=(X^{1}, X^{2})\colon \triangle_{T} \to \cV \oplus(\cV \otimes \cV)$
is called $\cV$-valued $\alpha$-RP if 
$\|X^i\|_{i \alpha} <\infty$ for $i=1,2$ and 
\begin{equation} \label{eq.0410-1}
X_{s, t}^{1}=X_{s, u}^{1}+X_{u, t}^{1}, 
\quad 
X_{s, t}^{2}
=X_{s, u}^{2}+X_{u, t}^{2}+X_{s,u}^{1} \otimes X_{u, t}^{1},
\qquad
s\le u \le t.
\end{equation}
The set of all $\cV$-valued $\alpha$-RPs is denoted by 
$\Omega_{\alpha} (\cV)$,
which is a complete metric space with the natural distance 
$d_{\alpha} (X, \hat{X}) :=\sum_{i=1}^2\|X^i - \hat{X}^i\|_{i \alpha}$.
(Note that $X^1$ and $X^2$ vanish on the diagonal.)
The homogeneous norm of $X$ is denoted by 
$\vertiii{X}_\alpha:= \|X^1\|_{\alpha} + \|X^2\|_{2 \alpha}^{1/2}$.
Obviously, $\Omega_{\alpha} (\cV)\subset \Omega_{\beta} (\cV)$
if $ \tfrac13 <\beta \le \alpha \le \tfrac12$.
The dilation by $\delta \in \R$ is defined by 
$\delta X = (\delta X^1, \delta^2 X^2)$. Clearly, 
$\vertiii{\delta X}_\alpha=|\delta|\cdot\vertiii{X}_\alpha$.

Now we recall the definition of a controlled path (CP)
with respect to a reference RP
$X =(X^{1}, X^{2})\in \Omega_{\alpha} (\cV)$.
Let $[a,b]\subset [0,T]$ be a subinterval.
We say that $(Y, Y^{\dagger}, Y^{\sharp})$ is a $\cW$-valued
CP with respect to $X$ on $[a,b]$ if
\[
(Y, Y^{\dagger}, Y^{\sharp}) \in 
\cC^\alpha ([a,b], \cW) \times 
\cC^\alpha ([a,b], L(\cV,\cW)) \times \cC_2^{2\alpha} ([a,b], \cW)
\]
and 
\begin{equation} \label{def_CP}
Y_{t}-Y_{s}
=Y_{s}^{\dagger} X_{s, t}^{1} +Y_{s, t}^{\sharp}, \qquad  (s, t) \in \triangle_{[a,b]}.
\end{equation}
The set of all such CPs with respect to $X$
is denoted by $\cQ^{\alpha}_X ([a,b], \cW)$.
For simplicity,
$(Y, Y^{\dagger}, Y^{\sharp})$ will often be written as $(Y, Y^{\dagger})$.
A natural seminorm of a CP is defined by
\[
\| (Y, Y^{\dagger}, Y^{\sharp}) \|_{\cQ^{\alpha}_X, [a,b]}
=\|Y^{\dagger} \|_{\alpha,[a, b]} + 
\|Y^{\sharp} \|_{2\alpha,[a, b]}
\]
$\cQ^{\alpha}_X ([a,b], \cW)$ is a Banach space with the norm 
$ |Y_a|_{\cW}+ |Y^\dagger_a|_{L(\cV,\cW)} +
\| (Y, Y^{\dagger}, Y^{\sharp}) \|_{\cQ^{\alpha}_X, [a,b]}
$.
When $[a,b]=[0,T]$, we simply write 
$\cQ^{\alpha}_X (\cW)$ and $\| \cdot \|_{\cQ^{\alpha}_X}$ instead.

\begin{example} \label{ex.0413}
Here are a few typical examples of CPs for 
a given RP $X \in \Omega_{\alpha} (\cV)$. 
(In the first three examples the time interval is $[0,T]$ just for simplicity.
It can be replaced by any subinterval $[a,b]$. )

\begin{enumerate} 
\item
For $\xi \in \cW$ and $\sigma\in L(\cV,\cW)$, 
$t \mapsto (\xi + \sigma X^1_{0,t}, \sigma)$ belongs to 
$\cQ^{\alpha}_X (\cW)$.
Note that $\cQ^{\alpha}_X$-seminorm of this element is zero.

\item
If $\varphi \in \cC^{2\alpha} (\cW)$, 
then obviously
$(\varphi, 0) \in \cQ^{\alpha}_X (\cW)$ with
$\|\varphi\|_{2\alpha} = \| (\varphi, 0) \|_{\cQ^{\alpha}_X}$.
In this way, we have a natural continuous embedding 
$\cC^{2\alpha} (\cW) \hookrightarrow \cQ^{\alpha}_X (\cW)$.

\item
Suppose that $(Y, Y^{\dagger})\in \cQ^{\alpha}_X (\cW)$ and 
$g\colon \cW \to \cW'$ is a $C^2$-function to 
another Euclidean space $\cW'$.
Then $(g(Y), g(Y)^\dagger) \in \cQ^{\alpha}_X (\cW')$
if we set $g(Y)_t:=g(Y_t)$ and
$g(Y)^\dagger_t := \nabla g(Y_t) Y^\dagger_t$,
where the right hand side is 
the composition of the two linear maps 
$\nabla g(Y_t)\in L(\cW,\cW')$ and $Y^\dagger_t \in L(\cV,\cW)$.
Using Taylor's  expansion, we can verify this fact as follows:
\begin{align} 
g(Y_{t})- g(Y_{s}) 
&=\nabla g(Y_{s})
\langle Y_{s, t}^{1}\rangle
   +\int_{0}^{1} d \theta(1-\theta) \nabla^{2} g(Y_{s}+\theta Y_{s, t}^{1})\langle Y_{s, t}^{1}, Y_{s, t}^{1}\rangle 
 \label{in.0411-1}  \\
&=\nabla g(Y_{s}) Y_{s}^{\dagger} X_{s, t}^{1} + g(Y)_{s, t}^{\sharp}
\nn
\end{align}
with
\[
g(Y)_{s, t}^{\sharp}:=\nabla g (Y_{s}) Y_{s, t}^{\sharp}
+\int_{0}^1 d \theta(1-\theta) \nabla^{2} g(Y_{s}+\theta Y_{s, t}^{1})\langle Y_{s, t}^{1}, Y_{s, t}^{1}\rangle.
\]
From the assumptions that $Y \in \cC^\alpha (\cW)$ and 
$Y^{\sharp}\in \cC_2^{2\alpha} (\cW)$, one can easily see 
that $g(Y)^{\sharp}\in \cC_2^{2\alpha} (\cW')$. 

\item
Concatenation of two CPs is also a CP.
Let $0\le a <b<c \le T$.
For $(Y, Y^{\dagger})\in \cQ^{\alpha}_X ([a,b], \cW)$ and
$(\hat{Y}, \hat{Y}^{\dagger})\in \cQ^{\alpha}_X ([b,c], \cW)$
with $(Y_b, Y^{\dagger}_b) = (\hat{Y}_b, \hat{Y}^{\dagger}_b)$,
their concatenation 
$(Z, Z^\dagger):=
(Y*\hat{Y}, Y^{\dagger}*\hat{Y}^{\dagger})$ can naturally be defined and 
belongs to $\cQ^{\alpha}_X ([a,c], \cW)$.
(Here, $*$ stands for the usual concatenation operation 
for two continuous paths.)
It is clear that $Z^{\dagger} \in \cC^\alpha ([a,c],\cW)$.
To prove that $Z^{\sharp} \in\cC_2^{2\alpha} (\cW)$,
it is sufficient to observe the following: 
For $a\le s \le b \le t \le c$, 
\begin{align} 
Z_{s, t}^{\sharp} &=Z_{t}-Z_{s}-Z_{s}^{\dagger} X_{s, t}^{1} 
\label{in.0411-2}\\
&=(Y_{b}-Y_{s}-Y_{s}^{\dagger} X_{s,b}^{1})
+(\hat{Y}_{t}-\hat{Y}_{b}-\hat{Y}_{b}^{\dagger} X_{b, t}^{1})
+(Y_{b}^{\dagger}-Y_{s}^{\dagger}) X_{b, t}^{1}
\nn\\
&=
Y^{\sharp}_{s, b}
 + \hat{Y}^{\sharp}_{b, t}
   +(Y_{b}^{\dagger}-Y_{s}^{\dagger}) X_{b, t}^{1}.
\nn
\end{align}
The right hand side is clearly dominated by a constant multiple of 
$(t-s)^{2\alpha}$.
\end{enumerate}
\end{example}

\subsection{Rough path integration of controlled paths}
Now we discuss integration of a CP 
$(Y, Y^{\dagger})\in \cQ^{\alpha}_X ([a,b], L(\cV, \cW))$
against a RP $X \in \Omega_{\alpha} (\cV)$.
It should be noted that $Y^{\dagger}$ takes values in 
\[
L ( \cV, L(\cV,\cW))\cong L^{(2)} (\cV \times\cV, \cW)
\cong
L ( \cV\otimes \cV,\cW),
\]
where $L^{(2)} (\cV \times\cV, \cW)$ stands for the vector space
of bounded bilinear maps from $\cV \times\cV$  to $\cW$.

First, we define 
\[
J_{s, t}=Y_{s} X_{s, t}^{1}+Y_{s}^{\dagger} X_{s, t}^{2},
\qquad
(s,t)\in \triangle_{[a,b]}.
\]
From \eqref{eq.0410-1} and \eqref{def_CP} we can easily see that
\begin{equation} \label{eq.0411-3}
J_{s, u}+J_{u, t}-J_{s, t}
=Y_{s, u}^{\sharp} X_{u, t}^{1}
 +Y_{s,u}^{\dagger, 1} X_{u, t}^{2},
\qquad 
a\le s\le u \le t \le b,
\end{equation}
where we set $Y_{s,u}^{\dagger, 1}= Y_{u}^{\dagger}-Y_{s}^{\dagger}$.
Let $\mathcal{P}=\{s=t_{0}<t_{1}<\cdots<t_{N}=t\}$
be a partition of $[s,t] \subset [a,b]$.
Its mesh size is denoted by $|\mathcal{P}|$.
We define 
$J_{s, t}(\mathcal{P})=\sum_{i=1}^{N} J_{t_{i-1}, t_{i}}$ and
\begin{equation} \label{eq.0411-4}
\int_{s}^{t} Y_{u} d X_{u}=\lim_{|\mathcal{P}| \searrow 0} J_{s, t}(\mathcal{P}),
\qquad (s,t)\in \triangle_{[a,b]}.
\end{equation}
The limit above is known to exist. This is called the RP integration.
It will turn out in the next proposition
that the RP integral against $X$ is again a CP with respect to $X$.
By the way it is defined,
this RP integral clearly has additivity with respect to the interval $[s,t]$.

\begin{proposition} \label{prop.gousei}
Let  $\tfrac13 < \alpha \le \tfrac12$ and $[a,b] \subset [0,T]$.
Suppose that $X \in \Omega_{\alpha} (\cV)$
and $(Y, Y^{\dagger})\in \cQ^{\alpha}_X ([a,b], L(\cV, \cW))$. Then,
the limit in \eqref{eq.0411-4} converges for all $(s,t)$.
Moreover, we have
\begin{equation}
\Bigl(\int_{a}^{\cdot} Y_{u} d X_{u}, Y\Bigr)
\in \cQ^{\alpha}_X ([a,b], \cW)
\label{eq.0413-3}
\end{equation}
with the following estimate:
\begin{align} 
\lefteqn{
\Bigl|
\int_{s}^{t} Y_{u} d X_{u} -(Y_{s} X_{s, t}^{1}+Y_{s}^{\dagger} X_{s, t}^{2})
\Bigr|_\cW
}
\label{in.0413-2} \\
&\le 
\kappa_{\alpha} (t-s)^{3\alpha} 
(\|Y^{\sharp} \|_{2\alpha,[a, b]} \|X^1\|_{\alpha,[a, b]}
+ \|Y^{\dagger} \|_{\alpha,[a, b]} \|X^2\|_{2\alpha,[a, b]}),
\,\, (s,t)\in \triangle_{[a,b]}.
\nn
\end{align}
Here, we set $\kappa_{\alpha}=2^{3\alpha}\zeta (3\alpha)$
with $\zeta$ being the usual Riemann zeta function.
\end{proposition}

\begin{proof} 
In this proof the norm of $\cW$ is denoted by $|\cdot|$ for brevity.
First we prove the convergence. 
For $\mathcal{P}$ given as above,
we can find $i ~(1 \le i \le N-1)$ such that $t_{i+1}-t_{i-1}\le 2(t-s)/(N-1)$
(see \cite[p. 52]{fh} for instance).
Then, we see from \eqref{eq.0411-3} that
\begin{align} 
|J_{s, t}(\mathcal{P}) - J_{s, t}(\mathcal{P}\setminus \{t_i\})|
&=
|J_{t_{i-1}, t_{i}} +  J_{t_{i}, t_{i+1}}  -  J_{t_{i-1}, t_{i+1}}|
\nn\\
&=
|Y_{t_{i-1}, t_{i}}^{\sharp} X_{t_{i}, t_{i+1}}^{1}
+Y_{t_{i-1}, t_{i}}^{\dagger, 1}X_{t_{i}, t_{i+1}}^{2}|
\nn\\
&\le
(\|Y^{\sharp} \|_{2\alpha,[a, b]} \|X^1\|_{\alpha,[a, b]}
+ \|Y^{\dagger} \|_{\alpha,[a, b]} \|X^2\|_{2\alpha,[a, b]}
 )
 \Bigl(\frac{2(t-s)}{N-1}\Bigr)^{3\alpha}.
 \nn
\end{align}
Extracting points one by one from $\mathcal{P}$ in this way
until it becomes the trivial partition $\{s,t\}$, we have
\begin{align} 
|J_{s, t}(\mathcal{P}) - J_{s, t}|
\le
\kappa_{\alpha} 
(\|Y^{\sharp} \|_{2\alpha,[a, b]} \|X^1\|_{\alpha,[a, b]}
+ \|Y^{\dagger} \|_{\alpha,[a, b]} \|X^2\|_{2\alpha,[a, b]})
(t-s)^{3\alpha}.
\label{in.0413-1}
\end{align}
Note that the condition $3\alpha >1$ is used here.
By a standard argument, \eqref{in.0413-1} implies that
$\{J_{s, t}(\mathcal{P}) \}_{\mathcal{P}}$ is Cauchy as 
$|\mathcal{P}| \searrow 0$.
Therefore, the limit in \eqref{eq.0411-4} exists and \eqref{in.0413-2} holds.
Since 
\[
|Y_{s}^{\dagger} X_{s, t}^{2}|
\le 
(|Y_{a}^{\dagger}|+ \|Y^{\dagger} \|_{\alpha,[a, b]}(b-a)^\alpha )
\|X^2\|_{2\alpha,[a, b]}(t-s)^{2\alpha},
\]
\eqref{in.0413-2} implies \eqref{eq.0413-3}.
\end{proof}

\subsection{Rough differential equations: The case of bounded and globally Lipschitz drift vector field}\label{sec.rde}

Now we discuss RDEs
in the framework of controlled path theory.
In this subsection we  assume that 
$\tfrac13 < \beta < \alpha \le \tfrac12$ and
$X \in \Omega_{\alpha} (\cV)\subset \Omega_{\beta} (\cV)$.

We set conditions on the coefficients of our RDE. 
Let $\sigma\colon \cW \to L(\cV,\cW)$ be of $C^3_b$
and let $f\colon \cW \times \cS \to \cW$ be a continuous map,
where $\cS$ is a metric space, satisfying the following condition:
\begin{equation}\label{cond.0506-1}
 \sup_{y \in \cW, z\in \cS} |f(y,z)|_{\cW}
+\sup_{y,y' \in \cW, y\neq y',z\in \cS} 
\frac{|f(y,z)-f(y',z) |_{\cW}}{|y-y'|_{\cW}} <\infty.
\end{equation}
The first and the second terms will be denoted by 
$\|f\|_{\infty}$ and $L_f$, respectively.

For an $\cS$-valued continuous path
$\psi \colon [0,T] \to \cS$, 
we consider the following RDE driven by $X$ with the initial point
 $\xi \in\cW$:
\begin{equation} \label{rde.0413}
Y_{t}
=\xi +\int_{0}^{t} f(Y_{s}, \psi_{s}) ds
+\int_{0}^{t} \sigma(Y_{s}) d X_{s},
\quad
Y^\dagger_t = \sigma(Y_t),  \qquad t \in [0,T].
\end{equation}
For every $(Y, Y^{\dagger})\in \cQ^{\beta}_X (\cW)$,
the right hand side of this system of equations also
belongs to $\cQ^{\beta}_X (\cW)$, due to Example \ref{ex.0413} and
Proposition \ref{prop.gousei}.
Therefore, \eqref{rde.0413} should be understood as
 an equality in $\cQ^{\beta}_X (\cW)$.
(Following \cite[Sections 8.5--8.6]{fh}, we slightly 
relax the H\"older topology 
of the space of CPs for quick proofs.
The estimate of $\|Y\|_{\beta}$ in the next proposition 
will be improved in Proposition \ref{prop.0702}.)

\begin{proposition} \label{prop.0429}
Let the assumptions be as above.
Then, for every $X\in \Omega_{\alpha} (\cV)$, $\xi \in \cW$ 
and $\psi$, there exists a unique global solution 
$(Y, Y^{\dagger})\in \cQ^{\beta}_X (\cW)$ of RDE \eqref{rde.0413}.
Moreover, it satisfies the following estimate:
there exist positive constants $c$ and $\nu$ 
independent of $X,  \xi, \psi, \sigma, f$ such that 
\[
\|Y\|_{\beta} \le c \{(K+1) (\vertiii{X}_{\alpha} +1)  \}^\nu,
\qquad 
X\in \Omega_{\alpha} (\cV). 
\]
Here, we set
$K := \| \sigma\|_{C_b^3}\vee \|f\|_{\infty}\vee L_f$. 
\end{proposition}

\begin{proof} 
In this proof the norm of $\cW$ is denoted by $|\cdot|$ for brevity.
Without loss of generality we may assume $T=1$.
Let $\tau \in (0,1]$ and $\xi \in \cW$. 

We define
$\cM^{\xi}_{[0, \tau]}, \cM^1_{[0, \tau]}, \cM^2_{[0, \tau]} \colon \cQ^{\beta}_X ([0,\tau], \cW) \to 
\cQ^{\beta}_X ([0,\tau], \cW)$ by 
\begin{align} 
\cM^{1}_{[0, \tau]} (Y, Y^\dagger) 
=
\Bigl( \int_{0}^{\cdot} \sigma(Y_{s}) d X_{s}, \,\, \sigma(Y)
\Bigr),
\quad
\cM^{2}_{[0, \tau]} (Y, Y^\dagger) = \Bigl(\int_{0}^{\cdot} f(Y_{s}, \psi_{s}) ds,
\,\, 0 \Bigr),
\nn
\end{align}
and $\cM^\xi_{[0, \tau]}=(\xi, 0) +\cM^{1}_{[0, \tau]}+ \cM^{2}_{[0, \tau]}$.
If $(Y, Y^\dagger)$ starts at $(\xi, \sigma(\xi))$, so does 
$\cM^\xi_{[0, \tau]}(Y, Y^\dagger)$.
A fixed point of $\cM^\xi_{[0, \tau]}$ is a solution of RDE \eqref{rde.0413}
on the interval $[0,\tau]$.

We also set 
\[
B_{[0, \tau]}^\xi
=\{(Y, Y^\dagger) \in \cQ^{\beta}_X ([0,\tau], \cW)\mid
\|(Y, Y^\dagger)\|_{ \cQ^{\beta}_X,  [0, \tau]} \le  1,  
  \,\, Y_0=\xi, Y^\dagger_0=\sigma(\xi) \}.
\]
This is something like a ball of radius $1$ centered at 
$t \mapsto (\xi + \sigma(\xi)X^1_{0,t}, \sigma(\xi))$.
Since the initial point $(Y_0, Y^\dagger_0)$ is fixed, 
$\|\cdot\|_{ \cQ^{\beta}_X,  [0, \tau]}$ works as a distance on this set.

For a while from now, we will work only on $[0,\tau]$ and therefore
omit $[0,\tau]$ from the subscript for notational simplicity.
We will often write  $Y_{s, t}^{1}:=Y_t -Y_s$.

For $(Y, Y^\dagger)\in B^\xi$ and $(\tilde{Y}, \tilde{Y}^\dagger) \in B^{\tilde\xi}$, the following estimates hold:
\begin{align} 
\|Y^{\dagger}\|_{\infty} &\leq
   |\sigma(\xi)|+\sup _{s \leq \tau} |Y_{s}^{\dagger}-Y_{0}^{\dagger}| 
\leq K+\|Y^{\dagger}\|_{\beta} \tau^{\beta} \leq K+1,
\label{sono1}\\
\|Y^{\dagger}- \tilde{Y}^\dagger\|_{\infty} 
 &\leq |Y^{\dagger}_0- \tilde{Y}^\dagger_0 |+\|Y^{\dagger}- \tilde{Y}^\dagger \|_{\beta}
 \le K |\xi -\tilde\xi| +\|Y^{\dagger}- \tilde{Y}^\dagger \|_{\beta},
\label{sono2}\\
|Y_{s, t}^{1}| &\le
   |Y_{s}^{\dagger} X_{s, t}^{1}|+|Y_{s, t}^{\sharp}|
\label{sono3}\\ 
  &\le (K+1)\|X^{1}\|_{\alpha} (t-s)^{\alpha}  
   +\|Y^{\sharp}\|_{2 \beta}(t-s)^{2 \beta}
   \nn\\
     &\le   (K+1)(\|X^{1}\|_{\alpha}+1)  (t-s)^{\alpha},
          \qquad (s,t) \in \triangle_\tau,
   \nn
   \\
|Y_{s, t}^{1}- \tilde{Y}_{s, t}^{1}| &\le
   |(Y_{s}^{\dagger} - \tilde{Y}_{s}^{\dagger})X_{s, t}^{1}|
     +|Y_{s, t}^{\sharp}-\tilde{Y}_{s, t}^{\sharp} |
     \label{sono4}\\
       &\le  (K |\xi -\tilde\xi| + \|Y^{\dagger}- \tilde{Y}^\dagger \|_{\beta})
        \|X^{1}\|_{\alpha} (t-s)^{\alpha}   
          +\|Y^{\sharp} -  \tilde{Y}^{\sharp} \|_{2 \beta}(t-s)^{2 \beta}
            \nn\\
       &\le  \bigl\{ K\|X^{1}\|_{\alpha} |\xi -\tilde\xi| 
       \nn\\
       &\quad +  (1+ \|X^{1}\|_{\alpha} ) \|  (Y, Y^\dagger)-(\tilde{Y},  \tilde{Y}^\dagger) \|_{\cQ^{\beta}_X}   \bigr\}     (t-s)^{\alpha},
               \quad (s,t) \in \triangle_\tau.
     \nn
     \end{align}
Note that $Y$ is in fact $\alpha$-H\"older continuous.
Hence, if $\tau$ is small, $\beta$-H\"older seminorms of 
$Y$ and $Y-\tilde{Y}$ can be made very small. 
From \eqref{sono4} we can easily see that
 \begin{align}
 \| Y- \tilde{Y} \|_\infty 
  &\le  |Y_0- \tilde{Y}_0| 
      + \| Y^1_{0,\cdot}- \tilde{Y}^1_{0,\cdot} \|_\infty
   \label{sono5}\\
     &\le (1 + K\|X^{1}\|_{\alpha}\tau^\alpha) |\xi -\tilde\xi|
     +(1+ \|X^{1}\|_{\alpha} ) \|  (Y, Y^\dagger)-(\tilde{Y},  \tilde{Y}^\dagger) \|_{\cQ^{\beta}_X}       
          \tau^{\alpha}.     \nn
     \end{align}

Let us consider $(\sigma(Y), \sigma(Y)^\dagger, \sigma(Y)^\sharp)$,
whose precise definition is given in the third item of Example \ref{ex.0413}.
First, $\sigma(Y)^\dagger_t = \nabla \sigma (Y_t) Y^\dagger_t
\in L(\cV \otimes \cV,\cW)$.
More precisely, it is defined by 
\begin{equation}\label{def.nabsigY}
\nabla \sigma (Y_t) Y^\dagger_t \la v\otimes v'\ra
:= \nabla \sigma (Y_t) \la Y^\dagger_t  v, v'\ra
\quad
\mbox{for $v, v' \in \cV$.}
\end{equation}
Then, we can see that
\begin{equation} \label{eq.0419-1}
\sigma(Y)^\dagger_t - \sigma(Y)^\dagger_s 
=
\int_{0}^{1} d \theta \nabla^{2} \sigma (Y_{s}+\theta Y_{s, t}^{1})\langle Y_{s, t}^{1}, Y_{t}^{\dagger}\rangle
+\nabla\sigma (Y_{s})(Y_{t}^{\dagger}-Y_{s}^{\dagger}).
\end{equation}
The second term on the right hand side is defined as in \eqref{def.nabsigY}.
Similarly, the precise meaning of the first term is given by 
$v\otimes v' \mapsto \int_0^1d\theta
\nabla^{2} \sigma (Y_{s}+\theta Y_{s, t}^{1})\langle Y_{s, t}^{1}, 
Y_{t}^{\dagger}v, v' \rangle$.
The remainder part reads:
\begin{align} 
\sigma(Y)^\sharp_{s,t}
&=\nabla \sigma (Y_{s}) \la Y_{s, t}^{\sharp}, \,\cdot\,\ra
\label{eq.0502-1}\\
&\qquad +\int_{0}^1 d \theta(1-\theta) \nabla^{2} \sigma (Y_{s}+\theta Y_{s, t}^{1})\langle Y_{s, t}^{1}, Y_{s, t}^{1}, \,\cdot\,\rangle 
   \in L(\cV,\cW).
   \nn
\end{align} 
From \eqref{sono1}, \eqref{sono3},
\eqref{eq.0419-1} and \eqref{eq.0502-1}, we can easily see that
\begin{align} 
\|\sigma(Y)^\dagger  \|_{\beta} 
&\le 
(K+1)^3(\|X^{1}\|_{\alpha}+1),
\nn\\
\|\sigma(Y)^\sharp  \|_{2\beta} 
&\le 
(K+1)^3(\|X^{1}\|_{\alpha}+1)^2.
\nn%
\end{align}

Next, we calculate the RP integral 
$(\int_0^\cdot  \sigma(Y)dX, \sigma(Y))$ 
using Proposition \ref{prop.gousei}.
It is easy to see that
\[
\|\sigma(Y)  \|_{\beta}  \le (K+1)^2(\|X^1\|_{\alpha}+1)
\tau^{\alpha-\beta}.
\]
We can see from \eqref{in.0413-2} and the above two estimates that 
\begin{align} 
\Bigl\|   \Bigl(\int_0^\cdot  \sigma(Y)dX
\Bigr)^\sharp \Bigr\|_{2\beta}
&\le   
\|\sigma(Y)^\dagger  \|_{\infty}  \|X^2\|_{2\beta}
\nn\\
&\qquad +
\kappa_\beta
(\|\sigma(Y)^{\sharp} \|_{2\beta} \|X^1\|_{\beta}
+ \|\sigma(Y)^{\dagger} \|_{\beta} \|X^2\|_{2\beta})
\nn\\
&\le
(1+2\kappa_\beta) (K+1)^3 (\vertiii{X}_{\alpha} +1)^3 \tau^{\alpha-\beta}.
\nn
\end{align}
Thus, we have an estimate of $\|\cM^{1} (Y, Y^\dagger) \|_{\cQ^{\beta}_X} $.
It is almost obvious from Example \ref{ex.0413} that
\begin{equation}
\|\cM^{2} (Y, Y^\dagger) \|_{\cQ^{\beta}_X}
\le 
\Bigl\|  \int_0^\cdot  f(Y_s, \psi_s)ds \Bigr\|_{2\beta}
\le 
K \tau^{1-2\beta} \le K \tau^{\alpha-\beta}.
\nn
\end{equation}
Combining these three estimates we obtain that
\[
\|\cM^{\xi} (Y, Y^\dagger) \|_{\cQ^{\beta}_X}
\le 
(2+2\kappa_\beta) (K+1)^3 (\vertiii{X}_{\alpha} +1)^3 \tau^{\alpha-\beta}.
\]
Hence, if 
\begin{equation}\label{in.0420-1}
\tau \le \lambda
\quad
 \mbox{ with } \quad
 \lambda :=\{8 \kappa_\beta (K+1)^3 (\vertiii{X}_{\alpha} +1)^3  \}^{-1/(\alpha-\beta)},
\end{equation}
then $\cM^\xi$ leaves $B^\xi$ invariant.
 Note also that $\kappa_\beta \ge 2$.
(For the rest of the proof
 we will assume $\tau$ satisfies this inequality.
 The constant ``$8$" in \eqref{in.0420-1} has no particular meaning.)

Now we will prove that $\cM^\xi$ is a contraction on
$B^\xi$ for $\tau$ as in \eqref{in.0420-1}.
We will calculate 
$(\sigma(Y)- \sigma(\tilde{Y}), 
  \sigma(Y)^\dagger -\sigma(\tilde{Y})^\dagger)$
  for $(Y, Y^\dagger)\in B^\xi$ and $(\tilde{Y}, \tilde{Y}^\dagger) \in B^{\tilde\xi}$.  
 
By straightforward (and slightly cumbersome) computations 
we obtain from \eqref{sono1}--\eqref{sono5} and \eqref{eq.0419-1} that
\begin{align*}
\lefteqn{
|\{\sigma(Y)^\dagger_t - \sigma(\tilde{Y})^\dagger_t\}
-
\{\sigma(Y)^\dagger_s - \sigma(\tilde{Y})^\dagger_s\}|
}
\nn\\
&\le
\int_{0}^{1} d \theta \bigl|
\nabla^{2} \sigma(Y_{s}+\theta Y_{s, t}^{1})\langle Y_{s, t}^{1}, Y_{t}^{\dagger}\rangle
- 
\nabla^{2} \sigma(\tilde{Y}_{s}+\theta \tilde{Y}_{s, t}^{1})
\langle \tilde{Y}_{s, t}^{1}, \tilde{Y}_{t}^{\dagger}\rangle
\bigr|
\nn\\
&\qquad
+
|\sigma(Y_{s}) (Y_{t}^{\dagger}-Y_{s}^{\dagger})
  -  \sigma(\tilde{Y}_{s}) 
  (\tilde{Y}_{t}^{\dagger}-\tilde{Y}_{s}^{\dagger}) |
  \nn\\
&\le
\{2(K+1)^3(\|X^{1}\|_{\alpha}+1)^2  \tau^{\alpha -\beta} 
( |\xi -\tilde\xi| +\|  (Y, Y^\dagger)-(\tilde{Y},  \tilde{Y}^\dagger) \|_{\cQ^{\beta}_X})
  \nn\\
  & \qquad\qquad 
    +K( |\xi -\tilde\xi| +\| Y^{\dagger}-\tilde{Y}^{\dagger} \|_{\beta} ) \} 
      (t-s)^\beta.
\end{align*}
This estimate and  \eqref{in.0420-1} imply that
\begin{equation} \label{in.0419-1}
\|\sigma(Y)^\dagger - \sigma(\tilde{Y})^\dagger\|_{\beta} 
\le 
(K+1) ( |\xi -\tilde\xi|+
  \|  (Y, Y^\dagger)-(\tilde{Y},  \tilde{Y}^\dagger) \|_{\cQ^{\beta}_X} ).
\end{equation}

In a very similar way, using \eqref{sono1}--\eqref{sono5} and \eqref{eq.0502-1}, we can estimate the difference of the remainder parts as follows:
\begin{align} 
 | \sigma(Y)^\sharp_{s,t} - \sigma(\tilde{Y})^\sharp_{s,t}|
 &\le 
 |\nabla \sigma (Y_{s}) \la Y_{s, t}^{\sharp}, \,\cdot\,\ra
 - \nabla \sigma (\tilde{Y}_{s}) \la \tilde{Y}_{s, t}^{\sharp}, \,\cdot\,\ra |
 \nn\\
 & +
 \int_{0}^1 d \theta |\nabla^{2} \sigma (Y_{s}+\theta Y_{s, t}^{1})\langle Y_{s, t}^{1}, Y_{s, t}^{1}, \,\cdot\,\rangle 
 -  
\nabla^{2} \sigma (\tilde{Y}_{s}+\theta \tilde{Y}_{s, t}^{1})\langle \tilde{Y}_{s, t}^{1}, \tilde{Y}_{s, t}^{1}, \,\cdot\,\rangle |
 \nn\\
 &\le 
 (K+1) ( |\xi -\tilde\xi|+
  \|  (Y, Y^\dagger)-(\tilde{Y},  \tilde{Y}^\dagger) \|_{\cQ^{\beta}_X})
 (t-s)^{2\beta},
  \nn
\end{align}
which immediately implies 
\begin{equation} \label{in.0421-1}
\|\sigma(Y)^\sharp - \sigma(\tilde{Y})^\sharp\|_{2\beta} 
\le 
(K+1) ( |\xi -\tilde\xi|+
  \|  (Y, Y^\dagger)-(\tilde{Y},  \tilde{Y}^\dagger) \|_{\cQ^{\beta}_X}).
\end{equation}

Using Proposition \ref{prop.gousei} again, we estimate
$(\int_0^\cdot  \{\sigma(Y)- \sigma(\tilde{Y})\}dX, \sigma(Y)- \sigma(\tilde{Y}))$.
From \eqref{in.0419-1}, \eqref{in.0421-1},  \eqref{sono4} and the fact 
$\sigma(Y_t)-\sigma(Y_s)=\int_0^1d\theta
 \nabla\sigma(Y_s+\theta Y^1_{s,t}) \langle Y^1_{s,t}\rangle$,
we have 
\begin{align} 
\|\sigma(Y)- \sigma(\tilde{Y})  \|_{\beta}  
&\le 
(K+1)^2(\|X\|_{\alpha}+1)^2 
\tau^{\alpha-\beta}
(  |\xi -\tilde\xi|+
\|  (Y, Y^\dagger)-(\tilde{Y},  \tilde{Y}^\dagger) \|_{\cQ^{\beta}_X})
\nn
\end{align}
and
\begin{align} 
\lefteqn{
\Bigl\|   \Bigl(\int_0^\cdot  \{ \sigma(Y)- \sigma(\tilde{Y})\}dX
\Bigr)^\sharp \Bigr\|_{2\beta}
}\nn\\
&\le   
\|\sigma(Y)^\dagger - \sigma(\tilde{Y})^\dagger \|_{\infty}  \|X^2\|_{2\beta}
\nn\\
&\quad +\kappa_\beta
(\|\sigma(Y)^{\sharp} - \sigma(\tilde{Y})^\sharp\|_{2\beta} \|X^1\|_{\beta}
+ \|\sigma(Y)^{\dagger}- \sigma(\tilde{Y})^\dagger\|_{\beta} \|X^2\|_{2\beta})
\nn\\
&\le
2(1+ \kappa_\beta) (K+1) (\vertiii{X}_{\alpha} +1)^2 \tau^{\alpha-\beta}
( |\xi -\tilde\xi|+\|  (Y, Y^\dagger)-(\tilde{Y},  \tilde{Y}^\dagger) \|_{\cQ^{\beta}_X}).
\nn
\end{align}
Using \eqref{in.0420-1}  we  obtain from the above estimates that 
 \[
 \|\cM^{1} (Y, Y^\dagger)-  \cM^{1} (\tilde{Y}, \tilde{Y}^\dagger)\|_{\cQ^{\beta}_X} \le 
\frac12 ( |\xi -\tilde\xi|+
   \|  (Y, Y^\dagger)-(\tilde{Y},  \tilde{Y}^\dagger) \|_{\cQ^{\beta}_X}).
\]

Since $f$ is  Lipschitz in the first argument, we can easily see 
from \eqref{sono5} that
\begin{align} 
\|\cM^{2} (Y, Y^\dagger) -  \cM^{2} (\tilde{Y}, \tilde{Y}^\dagger)\|_{\cQ^{\beta}_X}
&\le 
\Bigl\|  \int_0^\cdot \{ f(Y_s, \psi_s)-   f(\tilde{Y}_s, \psi_s)\}ds \Bigr\|_{2\beta}
\nn\\
&\le
K \|Y -\tilde{Y} \|_\infty  \tau^{1-2\beta} 
\nn\\
&\le
K (2 |\xi -\tilde\xi|+
\|  (Y, Y^\dagger)-(\tilde{Y},  \tilde{Y}^\dagger) \|_{\cQ^{\beta}_X})\tau^{\alpha-\beta}
\nn\\
&\le \frac14
( |\xi -\tilde\xi|+ \| (Y, Y^\dagger)-(\tilde{Y},  \tilde{Y}^\dagger) \|_{\cQ^{\beta}_X}).
 \nn
\end{align}
Summing up, we have 
\begin{equation}\label{key1}
\|\cM^\xi (Y, Y^\dagger) -  \cM^{\tilde\xi} (\tilde{Y}, \tilde{Y}^\dagger)\|_{\cQ^{\beta}_X}
\le 
\frac34 |\xi - \tilde\xi|+
\frac34
 \| (Y, Y^\dagger)-(\tilde{Y},  \tilde{Y}^\dagger) \|_{\cQ^{\beta}_X}.
\end{equation}
If $\xi =\tilde\xi$ in particular, this estimate implies that
$\cM^\xi$ is a contraction on $B^\xi=B^\xi_{[0,\tau]}$
and has a unique fixed point in this set, which is a local solution 
of the RDE.
Thus, we have obtained a solution on $[0,\lambda]$.
Note that $\tau$ (and $\lambda$) is determined by $\vertiii{X}_{\alpha}$ and $K$, 
but is chosen independent of $\xi$ and  $\psi$.

Next,  we do the same thing on the second interval 
$[\lambda, 2\lambda]$ with the initial condition $\xi$ 
at time $0$ being replaced by 
$Y_{\lambda}$ at time $\lambda$.
Since all the estimates above is independent of $\xi$ and $\psi$,
$(Y_s, Y^\dagger_s)_{s \in [\lambda, 2\lambda]}$
satisfies the same estimates as those for $(Y_s, Y^\dagger_s)_{s \in [0,\lambda]}$.
Concatenating them as in Example \ref{ex.0413}, 
we obtain  a solution on $[0, 2\lambda]$

We can continue this procedure to obtain a global 
$(Y_s, Y^\dagger_s)_{s \in [0, 1]}$.
There are $\lfloor \lambda^{-1}\rfloor +1$ subintervals,
where $\lfloor \cdot \rfloor$ stands for the integer part.
Except the last one, the length of each interval equals $\lambda$.
On each subinterval, $(Y, Y^\dagger)$ satisfies the same estimates.
In particular, Inequality \eqref{sono3} implies that
$\beta$-H\"older norm of $Y$ on each subinterval
 is dominated by 
$\{8 \kappa_\beta  (K+1)^2 (\vertiii{X}_{\alpha} +1)^2  \}^{-1}$.
By H\"older's inequality for finite sums, we can easily see that 
\begin{align*}
\|Y\|_{\beta, [0,1]} &\le \{8 \kappa_\beta  (K+1)^2 (\vertiii{X}_{\alpha} +1)^2  \}^{-1}(\lfloor\lambda^{-1}\rfloor +1)^{1-\beta}
\nn\\
&\le
c_{\alpha,\beta}  \{(K+1) (\vertiii{X}_{\alpha} +1)  \}^{\nu},
\end{align*}
where  $\nu :=3(1-\beta)/(\alpha-\beta) -2 >0$ and $c_{\alpha,\beta}>0$
is a constant which depends only on $\alpha, \beta$.

The uniqueness for this type of RDEs is well-known.
So we just give a quick explanation.
The uniqueness is a time-local issue, so it suffices to prove 
that any two solutions,
 $(Y, \sigma (Y))$ and $(\tilde{Y}, \sigma (\tilde{Y}))$,
of RDE \eqref{rde.0413}  must coincide near $t=0$. 
Take any $\beta' \in (1/3, \beta)$ and we work in $\beta'$-H\"older
topology instead of $\beta$-H\"older topology. 
If $\tau$ is small enough, both 
$(Y, \sigma (Y))$ and $(\tilde{Y}, \sigma (\tilde{Y}))$ 
restricted to $[0,\tau]$ belong to $B^{\xi}_{[0,\tau]}$.
But, we have already proved that there is only one fixed point 
of $\cM^{\xi}_{[0,\tau]}$ in this ball. 
Hence, $(Y, \sigma (Y))$ and $(\tilde{Y}, \sigma (\tilde{Y}))$ 
must be identically equal on $[0,\tau]$. 
 \end{proof}

 \begin{remark} \label{rem.loc.sol}
 {\rm (i)}~By examining the proof of  Proposition \ref{prop.0429},
 one naturally realize the following:
  Just to prove the existence of a unique global solution
of RDE \eqref{rde.0413} for any given $\psi$, $X$ and $\xi$, 
 it suffices to assume that $\sigma$ is of $C^3_b$ and $f$
 satisfies that
 \[
 \sup_{y \in \cW, t \in [0,T]} |f(y,\psi_t)|_{\cW} + 
\sup_{y,y' \in \cW, y\neq y', t \in [0,T]} 
\frac{|f(y,\psi_t)-f(y',\psi_t) |_{\cW}}{|y-y'|_{\cW}} <\infty.
 \]
 \noindent
 {\rm (ii)}~By a standard cut-off argument, 
 it immediately follows from {\rm (i)} above that if
 $\sigma$ is of $C^3$ and $f$ is locally Lipschitz 
 continuous in the following sense
 \[
\sup_{|y| \vee |y'| \le N, y\neq y', t \in [0,T]} 
\frac{|f(y,\psi_t)-f(y',\psi_t) |_{\cW}}{|y-y'|_{\cW}} <\infty,
\quad
N \in \N,
\]
then RDE \eqref{rde.0413} has a unique
local solution for any given $\psi$, $X$ and $\xi$.
Hence, a unique solution exists up to either the explosion time 
or the time horizon $T$.
  \end{remark}


Together with RDE \eqref{rde.0413}, we also consider the following
RDE with the same $X$, $\sigma$ and $\xi$:
\begin{equation} \label{rde.0506}
\tilde{Y}_{t}
=\xi +\int_{0}^{t} \tilde{f} (\tilde{Y}_{s}, \tilde{\psi}_{s}) ds
+\int_{0}^{t} \sigma (\tilde{Y}_{s}) d X_{s},
\quad
\tilde{Y}^\dagger_t = \sigma(\tilde{Y}_t),  \qquad t \in [0,T].
\end{equation}
We assume that 
$ \tilde{f}\colon \cW \times \cS \to \cW$ is also 
continous and satisfies Condition \eqref{cond.0506-1}.
Let $\tilde{\psi}\colon [0,T] \to \cS$ be another continuous path in $\cS$.

\begin{proposition} \label{prop.0506}
Let $\sigma, f, \tilde{f}, \xi$ be as above.
For $X\in \Omega_{\alpha} (\cV)$, $\xi \in \cW$ 
and $\psi, \tilde{\psi}$, 
denote by $(Y, \sigma(Y))$ and $(\tilde{Y}, \sigma(\tilde{Y}))$ 
the unique solutions of RDEs \eqref{rde.0413} and \eqref{rde.0506}
on $[0,T]$, respectively.
For a bounded, globally Lipschitz map
$g \colon \cW \to \cW$,  set
\begin{equation} \label{eq.0506-1}
M_t :=
(Y_t - \tilde{Y}_{t})
-
\int_{0}^{t} \{g(Y_{s})  - g(\tilde{Y}_{s}) \}ds
-
\int_{0}^{t} \{\sigma (Y_{s}) - \sigma (\tilde{Y}_{s}) \}d X_{s},
 \quad t \in [0,T].
 \end{equation}
 Then, $M \in \cC^{1} (\cW)$ and we have the following estimate
  for every $\beta \in (\tfrac13, \alpha)$: there exist positive constants $c$ and $\nu$ such that   %
 \begin{equation} \label{in.0506-2}
 \|Y - \tilde{Y} \|_{\beta} \le c
  \exp \bigl[c (K'+1)^\nu (\vertiii{X}_{\alpha} +1)^\nu \bigr]
     \| M\|_{2\beta}.
 \end{equation}
 Here, we set $K'= \max\{\| \sigma\|_{C_b^3}, \|f\|_{\infty}, L_f, 
 \|\tilde{f}\|_{\infty}, L_{\tilde{f}}, \|g\|_{\infty}, L_g\}$ and
the constants $c$ and $\nu$ are independent of 
 $X,  \xi, \psi, \tilde{\psi},  \sigma, f, \tilde{f}, g, M$. 
\end{proposition} 

\begin{proof} 
Without loss of generality we assume $T=1$.
For simplicity we write $(Y, Y^\dagger)$ and $(\tilde{Y}, \tilde{Y}^\dagger)$
for $(Y, \sigma(Y))$ and $(\tilde{Y}, \sigma(\tilde{Y}))$, respectively.
It is easy to see that $M \in \cC^{1} (\cW)$.
Hence, \eqref{eq.0506-1} is in fact an equality in $\cQ^{\beta}_X (\cW)$
with the $\dagger$-parts being clearly equal.

Mimicking \eqref{in.0420-1}, we set 
$\lambda' :=\{8 \kappa_\beta (K'+1)^3 (\vertiii{X}_{\alpha} +1)^3  \}^{-1/(\alpha-\beta)}$.
Set $s_j : = j \lambda'$ for $0 \le j \le \lfloor 1/\lambda'\rfloor$
and $s_N :=1$ with $N:= \lfloor 1/\lambda'\rfloor +1$.
Then, on each subinterval $[s_{j-1}, s_j]$,
$(Y, Y^\dagger) \in B_{[s_{j-1}, s_j]}^{\xi_j}$,  
 $(\tilde{Y}, \tilde{Y}^\dagger) \in B_{[s_{j-1}, s_j]}^{\tilde{\xi}_j}$
 and the 
estimates in the proof of Proposition \ref{prop.0429} are available
(with $K$ being replaced by $K'$). 
From \eqref{key1} we have for all $j$ that
\begin{align}
\Bigl\|
\int_{s_{j-1}}^{\cdot} \{g(Y_{s})  - g(\tilde{Y}_{s}) \}ds
-
\int_{s_{j-1}}^{\cdot} \{\sigma (Y_{s}) - \sigma (\tilde{Y}_{s}) \}d X_{s}
\Bigr\|_{\cQ^{\beta}_X,  [s_{j-1}, s_j]}
\nn\\
\le 
\frac34 |\xi_{j-1} - \tilde\xi_{j-1}|+
\frac34
 \| (Y, Y^\dagger)-(\tilde{Y},  \tilde{Y}^\dagger) \|_{\cQ^{\beta}_X, [s_{j-1}, s_j]},
 \nn
\end{align}
where we write $\xi_j = Y_{s_j}$ and $\tilde{\xi}_j = \tilde{Y}_{s_j}$.
Taking the seminorms of both sides of \eqref{eq.0506-1} on 
each subinterval, 
we can easily see that
\begin{equation}\label{key3}
\| (Y, Y^\dagger)-(\tilde{Y},  \tilde{Y}^\dagger) \|_{\cQ^{\beta}_X, [s_{j-1}, s_j]}
\le
4\|M\|_{2\beta} + 3 |\xi_{j-1} - \tilde\xi_{j-1}|.
\end{equation}
Plugging this into \eqref{sono4} , we obtain for all $j$ that
\begin{align}
|Y_{s, t}^{1}- \tilde{Y}_{s, t}^{1}| 
&\le
\{K' \|X^{1}\|_{\alpha} |\xi_{j-1} -\tilde\xi_{j-1}| 
       \label{in.0509-2}\\
      &\qquad \quad +  (1+ \|X^{1}\|_{\alpha} ) 
       (4\|M\|_{2\beta} + 3 |\xi_{j-1} - \tilde\xi_{j-1}|)   \}
         (t-s)^{\alpha}
          \nn\\
           &\le   |\xi_{j-1} -\tilde\xi_{j-1}| 
              + \|M\|_{2\beta},
                                        \qquad (s,t) \in \triangle_{[s_{j-1}, s_j]} 
                                         \nn
                                         \end{align}
               and in particular
               \begin{equation}\nn
               | \xi_{j} -\tilde\xi_{j} | 
                \le 2 |\xi_{j-1} -\tilde\xi_{j-1}| 
                    + \|M\|_{2\beta}.               
                   \end{equation}
               By mathematical induction we have 
               \[
               | \xi_{j} -\tilde\xi_{j} | \le (1+2^1 +\cdots + 2^{j-1}) \|M\|_{2\beta}
                 = (2^{j}-1)\|M\|_{2\beta}, 
                  \quad 1\le j \le N.                          
                   \]
               Then, we see from \eqref{in.0509-2}  that 
               \begin{align} 
               \|Y^{1}- \tilde{Y}^{1}\|_{\beta, [s_{j-1}, s_j]} 
               &\le
               \{K' \|X^{1}\|_{\alpha} |\xi_{j-1} -\tilde\xi_{j-1}| 
               \label{in.0510-1}
\\
&\qquad        +  (1+ \|X^{1}\|_{\alpha} ) 
       (4\|M\|_{2\beta} + 3 |\xi_{j-1} - \tilde\xi_{j-1}|)   \}
         (\lambda')^{\alpha -\beta}
         \nn\\
          & \le  4 (1+K')(1+ \|X^{1}\|_{\alpha})
          2^N (\lambda')^{\alpha -\beta} \|M\|_{2\beta}, 
                  \quad 1\le j \le N.                                 
               \nn
\end{align}
By H\"older's inequality for finite sums
and the trivial inequality $N^{1-\beta} 2^N \le 2^{2N}$, 
we see that 
\begin{align} \|Y^{1}- \tilde{Y}^{1}\|_{\beta, [0,1]}
&\le 
4(1+K') (1+ \|X^{1}\|_{\alpha})
          2^N (\lambda')^{\alpha -\beta} \|M\|_{2\beta}N^{1-\beta} 
                \label{in.0510-2}
\\
           &\le 
             c_{\alpha,\beta}  
             (K'+1) (\vertiii{X}_{\alpha} +1)       
             \nn\\
              & \qquad \times 
             \exp \bigl[  c_{\alpha,\beta}\{(K'+1) (\vertiii{X}_{\alpha} +1)  \}^{3/(\alpha-\beta)}  \bigr]       \|M\|_{2\beta},
                       \nn
\end{align}
where $c_{\alpha,\beta}>0$
is a constant which depends only on $\alpha, \beta$.
By adjusting the constant $c_{\alpha,\beta}$, 
we can easily obtain \eqref{in.0506-2} from \eqref{in.0510-2}.
               \end{proof}

\subsection{Rough differential equations:  The case of bounded and locally Lipschitz drift vector field}
\label{sec.rde.local}
 
In this subsection  we continue to study RDE  \eqref{rde.0413},
where $\psi \colon [0,T] \to \cS$ is a continuous path in $\cS$.
We still assume that $\tfrac13 < \beta <\alpha \le \tfrac12$
and  $\sigma$ is of $C^3_b$,
but relax the global Lipschitz condition \eqref{cond.0506-1} on $f$.

\begin{proposition} \label{prop.0702}
Suppose that $\sigma$ is of $C^3_b$,
$\|f\|_{\infty} := \sup_{y \in \cW, z \in \cS} |f(y,z)|_{\cW}<\infty$
and 
\begin{equation}\nn
\sup_{|y| \vee |y'| \le N, y\neq y', t \in [0,T]} 
\frac{|f(y,\psi_t)-f(y',\psi_t) |_{\cW}}{|y-y'|_{\cW}} <\infty,
\quad
N \in \N,
\end{equation}
then RDE  \eqref{rde.0413} has a unique global solution 
$(Y, Y^{\dagger})\in \cQ^{\beta}_X (\cW)$ for every $X \in  \Omega_{\alpha} (\cV)$, $\xi \in \cW$ and $\psi$. 
Moreover, there 
exists a constant $C>0$ such that 
\begin{equation} \label{in.0616-1}
\| Y\|_{\beta} \le 
C \{(\|\sigma\|_{C^2_b} \vertiii{X}_{\alpha})^{1/\beta} +\|\sigma\|_{C^2_b} \vertiii{X}_{\alpha}+\|f\|_\infty\},
\quad X \in  \Omega_{\alpha} (\cV).
\end{equation}
Here, $C$ is independent of $X, \xi, \psi, f, \sigma$.
\end{proposition}

\begin{proof} 
We mimic the proof of a priori estimates in \cite[Section 8.4]{fh}. 
In this proof positive constants $c_i$'s are independent of 
$X, \xi, \psi, f, \sigma$.
Without loss of generality we may assume $T=1$.
We will write  $L:= \|f\|_\infty$
and $L':=\|\sigma\|_{C^2_b}$ for simplicity.
The norms of finite dimensional vector spaces are
simply denoted by $|\cdot|$.

By Remark \ref{rem.loc.sol},
 a unique local solution exists, which is denoted by 
 $\{(Y_t, Y^{\dagger}_t)\}_{0\le t \le S}$, $S \in (0,1]$.
 If $\| Y\|_{\beta, [0,S]}$ is dominated by the right hand side of 
  \eqref{in.0616-1}, which is independent of $S$, 
  then $Y$  does not explode and therefore extends to a 
  global solution.

First, we consider the case $L' \le 1/2$. For $(s,t)\in \triangle_S$,
\begin{align} 
|Y^\sharp_{s,t}|
 &= |Y_{s,t} - \sigma (Y_s) X^1_{s,t} |
 \nn\\
    &\le
     \Bigl|  
       \int_s^t \sigma (Y_u) dX_u -   \sigma (Y_s) X^1_{s,t}
            -  \nabla\sigma (Y_s)  \sigma (Y_s) X^2_{s,t}       
              \Bigr| 
             \nn\\
               &\qquad 
               + |\nabla\sigma (Y_s)  \sigma (Y_s) X^2_{s,t}|
                  +\Bigl|  
                       \int_s^t   f (Y_u, \psi_u) du \Bigr| 
                         \nn\\
                          &\le 
                           \kappa_{\beta} 
(\|\sigma(Y)^{\sharp} \|_{2\beta,[s,t]} \|X^1\|_{\beta,[s,t]}
+ \|\sigma(Y)^\dagger \|_{\beta,[s,t]} \|X^2\|_{2\beta,[s,t]})(t-s)^{3\beta}
  \nn\\
     &\qquad + \|X^2\|_{2\beta,[s,t]} (t-s)^{2\beta}  +L(t-s),
       \nn
\end{align}  
  where Proposition \ref{prop.gousei} is used  for the last inequality.
 Note that $\|\sigma(Y)^\dagger \|_{\beta,[s,t]}\le \|Y \|_{\beta,[s,t]}$.
For $h \in (0,S]$, we set 
$\| \cdot\|_{\beta; h}:=\sup_{0<t-s\le h} \| \cdot\|_{\beta, [s,t]}$.
Then,  it immediately follows that
\begin{align} 
\| Y^\sharp\|_{2\beta; h} 
&\le 
\kappa_{\beta} 
(\|\sigma(Y)^{\sharp} \|_{2\beta;h} \|X^1\|_{\beta;h}
+ \|Y \|_{\beta;h} \|X^2\|_{2\beta;h})  h^\beta
\label{in.0617-1}
\\
 &\qquad 
     + \|X^2\|_{2\beta;h}   +Lh^{1-2\beta}.
\nn
\end{align}

Next we calculate $\sigma(Y)^{\sharp}$. By definition, we have
\begin{align} 
\sigma(Y)^{\sharp}_{s,t} 
&= 
\sigma(Y_t)- \sigma(Y_s) - \nabla\sigma(Y_s)  Y^\dagger_s X^1_{s,t}
\nn\\
&=
\sigma(Y_t)- \sigma(Y_s) - \nabla\sigma(Y_s)\langle Y^1_{s,t}\rangle
+
\nabla\sigma(Y_s)\langle Y^{\sharp}_{s,t}\rangle.
\nn
\end{align}
From this and Taylor's formula, 
\begin{align} 
\| \sigma(Y)^\sharp\|_{2\beta; h} 
&\le
\frac12 \|Y\|_{\beta;h}^2 + \|Y^\sharp\|_{2\beta; h}. 
\nn
\end{align}

We put the above inequality back into \eqref{in.0617-1}. 
Then we can easily see that
\begin{align} 
\| Y^\sharp\|_{2\beta; h} 
&\le 
\kappa_{\beta} 
  ( \frac12 \|Y\|_{\beta;h}^2 + \|Y^\sharp\|_{2\beta; h}) 
   \|X^1\|_{\beta;h} h^\beta
 + \kappa_{\beta}   \|Y \|_{\beta;h} \|X^2\|_{2\beta;h}  h^\beta
\label{in.0617-2}\\
 &\qquad 
     + \|X^2\|_{2\beta;h}   +Lh^{1-2\beta}.
\nn
\end{align}
Hence, if $h$ is so small that 
$\kappa_{\beta}  \vertiii{X}_{\beta} h^\beta 
= \kappa_{\beta} (  \|X^1\|_{\beta} + \|X^2\|_{2\beta}^{1/2} ) h^\beta \le 1/2$ 
holds,  we see from \eqref{in.0617-2} that
\begin{align} 
\| Y^\sharp\|_{2\beta; h} 
&\le 
\kappa_{\beta} \|Y\|_{\beta;h}^2  \|X^1\|_{\beta;h} h^\beta
 + 2\kappa_{\beta}\|Y \|_{\beta;h} \|X^2\|_{2\beta;h}  h^\beta
     + 2\|X^2\|_{2\beta;h}   +2Lh^{1-2\beta}
     \label{in.0617-3}\\
      &\le \frac12 \|Y\|_{\beta;h}^2 + \|Y \|_{\beta;h} \|X^2\|_{2\beta;h}^{1/2}
        + 2\|X^2\|_{2\beta;h}   +2Lh^{1-2\beta}
          \nn\\
      &\le 
        \|Y\|_{\beta;h}^2         + 3\|X^2\|_{2\beta;h}   +2Lh^{1-2\beta}.
                            \nn
\end{align}
Since $Y^1_{s,t}= Y^\dagger_s X^1_{s,t} + Y^{\sharp}_{s,t}= \sigma (Y_s) X^1_{s,t} + Y^{\sharp}_{s,t}$ by definition, 
we see from \eqref{in.0617-3} and $\kappa_\beta \ge 2$ that
\begin{align} 
\|Y\|_{\beta;h} h^\beta
\le \|X^1\|_{\beta;h} h^\beta + \|Y^\sharp\|_{2\beta; h}h^{2\beta}
\le
\vertiii{X}_{\beta} h^\beta+2Lh+
(\|Y\|_{\beta;h} h^\beta)^2.
\nn
\end{align}
If we set $\lambda_h = \vertiii{X}_{\beta} h^\beta+2Lh$
and 
$\varphi_h =\|Y\|_{\beta;h} h^\beta$ for $h\in (0,S]$, then we have
\begin{equation} \label{in.0617-4}
0 \le \varphi_h \le \lambda_h +  \varphi_h^2
\qquad\quad
\mbox{if $\vertiii{X}_{\beta}
h^\beta \le 1/(2\kappa_{\beta})$.}
\end{equation}
Note that $\varphi_h$ is left-continuous and
non-decreasing in $h$. Moreover, $\lim_{h\searrow 0}\varphi_h =0$.
So, $\varphi_h \le 1/4$ if $h \in (0,\eta)$ for some very small $\eta >0$.
However, this $\eta$ seems to depend on $Y$.
We will show that there exists $\eta >0$ 
which depends only on $\vertiii{X}_{\beta}$.

It is easy to see that
$(1/2) -\sqrt{ (1/4) - u} \le 2u$ for $u \in [0,1/8]$.
Hence, it immediately follows from \eqref{in.0617-4} that
one and only one of the following two conditions,  \eqref{in.0811a} and 
\eqref{in.0811b}, must hold if $\lambda_h \le 1/8$:
\begin{align} 
\varphi_h &\ge \frac12 + \sqrt{\frac14 - \lambda_h} \ge \frac34,
\label{in.0811a}
\\
\varphi_h &\le \frac12 - \sqrt{\frac14 - \lambda_h} \le 2\lambda_h
\le \frac14.
\label{in.0811b}
\end{align}
If $0<\delta <h$, then $\|Y\|_{\beta;h+\delta} \le 2^{1-\beta} \|Y\|_{\beta;h}$.
This implies that the relative jump size of $\varphi$ is
at most $2^{1-\beta}$, i.e. 
$\lim_{\delta\searrow 0}\varphi_{h+\delta} \le 2^{1-\beta}\varphi_h$.
Hence, as $h$ increases from $0$ under the conditions that
$\lambda_h \le 1/8$ and $\vertiii{X}_{\beta} 
h^\beta \le 1/(2\kappa_{\beta})$, 
$\varphi_h$
cannot jump from $[0, 1/4]$ to $[3/4, \infty)$.
So $\varphi_h$ must always satisfy \eqref{in.0811b} for all such $h$.
By the above argument, we see that if 
\[
h \le \Lambda, 
\qquad
\mbox{where we set} \quad 
\Lambda :=
\min\{ [16 \vee (2\kappa_{\beta}) ]^{-1/\beta} \vertiii{X}_{\beta}^{-1/\beta}, 
(32 L)^{-1}, S\}
\]
(this includes the case $\vertiii{X}_{\beta}=0$ or $L=0$), then \eqref{in.0811b} holds and hence
\begin{align} 
\|Y\|_{\beta;h}
&\le 
\varphi_h  h^{-\beta}
\le 2 \lambda_h h^{-\beta}
\le
 2\vertiii{X}_{\beta} +4Lh^{1-\beta}
\label{in.0617-5} 
\\
    &\le
\left\{
\begin{array}{ll}
2 \vertiii{X}_{\beta} + c_1 L \min\{ \vertiii{X}_{\beta}^{-(1-\beta)/\beta}, 1\} 
 & (\mbox{if $L\le 1$}),
 \\
2 \vertiii{X}_{\beta} + c_1 L \min\{ \vertiii{X}_{\beta}^{-(1-\beta)/\beta}, 1/L^{1-\beta}\} 
 & (\mbox{if $L\ge 1$}),
\nn
\end{array}
\right.
 \end{align}
also holds.

Set $s_j : = j \Lambda$ for $0 \le j \le \lfloor S/\Lambda\rfloor$
and $s_N :=S$ with $N:= \lfloor S/\Lambda\rfloor +1$. Then, 
\begin{align} 
N\le
c_2\max\{\vertiii{X}_{\beta}^{1/\beta}, L, 1\}+1
\le
\left\{
\begin{array}{ll}
c_3 ( \vertiii{X}_{\beta}^{1/\beta} + 1)
 & (\mbox{if $L\le 1$}),
 \\
 c_3 ( \vertiii{X}_{\beta}^{1/\beta} + L)
 & (\mbox{if $L\ge 1$}).
 \label{in.0617-6} 
\end{array}
\right.
 \end{align}
On  each subinterval $[s_{j-1}, s_j]$, 
Inequality \eqref{in.0617-5} is available.

First we calculate the case $L\le 1$.
Using \eqref{in.0617-5}, \eqref{in.0617-6}  and
H\"older's inequality for finite sums again, we have
\begin{align*}
\|Y\|_{\beta, [0,S]}
&\le  \|Y\|_{\beta;\Lambda} N^{1-\beta}
\nn\\
&\le c_4 (\vertiii{X}_{\beta} + L \min\{ \vertiii{X}_{\beta}^{-(1-\beta)/\beta}, 1\} 
 )  (\vertiii{X}_{\beta}^{1/\beta} + 1)^{1-\beta}
 \nn\\
&\le c_5 (\vertiii{X}_{\beta}^{1/\beta} + \vertiii{X}_{\beta}  +L).
\end{align*}
The case $L\ge 1$ is quite similar;
\begin{align*}
\|Y\|_{\beta, [0,S]}
&\le  \|Y\|_{\beta;\Lambda} N^{1-\beta}
\nn\\
&\le  c_6 (\vertiii{X}_{\beta} + L \min\{ \vertiii{X}_{\beta}^{-(1-\beta)/\beta}, 
   1/L^{1-\beta}\} 
 )    (\vertiii{X}_{\beta}^{1/\beta} +L)^{1-\beta}
 \nn\\
&\le  c_7 (\vertiii{X}_{\beta}^{1/\beta} + \vertiii{X}_{\beta}L^{1-\beta}  +L)
 \nn\\
&\le  c_8 (\vertiii{X}_{\beta}^{1/\beta}  +L)
\le  c_8 (\vertiii{X}_{\beta}^{1/\beta} + \vertiii{X}_{\beta}+L).
\end{align*}
In the second to the last inequality, we used Young's inequality.
Either way, we have proved that $Y$ does not explode and 
\begin{equation}\label{in.0617-7}
\|Y\|_{\beta, [0,1]}
\le  c_9 (\vertiii{X}_{\beta}^{1/\beta} + \vertiii{X}_{\beta}+L).
\end{equation}
Thus, the case $L' \le 1/2$ is done.

Finally, we consider the general case.  
Since the $L' =0$ is trivial, we may assume $L' >0$.
Note that $(Y, \sigma(Y))$ solves RDE \eqref{rde.0413}
if and only if $(Y, \sigma(Y)/2L')$ 
solves RDE \eqref{rde.0413} with $\sigma$ and $X$ being 
replaced by $\sigma/2L'$ and $2L' X$, respectively.
Here, $2L' X$ is the dilation of $X$ by $2L'$. 
Therefore, we have 
\begin{equation}\nn
\|Y\|_{\beta, [0,1]}
\le  c_9 \{(2L'\vertiii{X}_{\beta})^{1/\beta} +2L' \vertiii{X}_{\beta}+L\}.
\end{equation}
Noting that $\vertiii{X}_{\beta}\le \vertiii{X}_{\alpha}$ since $\beta<\alpha$,
we complete the proof.
  \end{proof}

Next we generalize Proposition \ref{prop.0506} to the case of 
bounded and locally Lipschitz drift vector fields.
We consider RDEs  \eqref{rde.0413}   and  \eqref{rde.0506} again 
and then define $M$ as in \eqref{eq.0506-1}
for a continuous map $g\colon \cW \to \cW$.
This time we do not impose the global Lipschitz  property 
on $f, \tilde{f}, g$, however. 
For $N \in \N$, we set
\begin{equation}\nn
L_f^N :=
\sup_{|y| \vee |y'| \le N, y\neq y', z \in \cS} 
\frac{|f(y,z)-f(y',z) |_{\cW}}{|y-y'|_{\cW}}.
\end{equation}
In the essentially same way we set  $L_{\tilde{f}}^N$  and $L_g^N$, too.

\begin{proposition}\label{prop.0703}
For $X\in \Omega_{\alpha} (\cV)$, $\xi \in \cW$ 
and continuous paths $\psi$ and $\tilde{\psi}$, we
denote by $(Y, \sigma(Y))$ and $(\tilde{Y}, \sigma(\tilde{Y}))$ 
the unique solutions of RDEs \eqref{rde.0413} and \eqref{rde.0506}
on $[0,T]$, respectively.
We assume that $\sigma$ is of $C^3_b$ and that $f, \tilde{f}, g$ are all bounded and satisfy 
\[
L_{f}^N +L_{\tilde{f}}^N+ L_{g}^N = O (N^r) \qquad \mbox{as $N\to\infty$} 
\]
for some $r>0$. Here, $O$ stands for Landau's large $O$.

Then, $M \in \cC^{1} (\cW)$ and we have the following estimate
  for every $\beta \in (\tfrac13, \alpha)$: there exist positive constants 
  $C$ and $\nu$ such that   %
 \begin{equation} \nn
 \|Y - \tilde{Y} \|_{\beta} \le 
 C\exp (C  \vertiii{X}_{\alpha}^{\nu} ) \| M\|_{2\beta}.
 \end{equation}
Here, $C$ and $\nu$ are independent of 
 $X,  \xi, \psi, \tilde{\psi},  M$. 
 Recall that $M$ was defined in \eqref{eq.0506-1}.
\end{proposition}

\begin{proof}
In this proof, $c_i$'s are positive constants independent of 
$X,  \xi, \psi, \tilde{\psi},  M$.

We use Proposition \ref{prop.0702}. From \eqref{in.0616-1}  we see that 
\[
\| Y\|_{\beta} \vee \|\tilde{Y}\|_{\beta} \le 
c_1 ( \vertiii{X}_{\alpha}^{1/\beta} +\vertiii{X}_{\alpha}+1)
\]
for some constant $c_1 \ge 1$.
We denote the right hand side by $R^X \,(\ge 1)$.

Take a Lipschitz continuous function $\chi \colon \cW \to [0,1]$ 
such that {\rm (i)} $\chi (y) =1$ if $|y| \le R^X$, 
{\rm (ii)} $\chi (y) =0$ if $|y| \ge 2R^X$ and
{\rm (iii)} the Lipschitz constant of $\chi$ is at most $1$.
Then, $\chi f$ is globally Lipschitz (in $y$) since 
\begin{align*}
|\chi (y) f (y,z) - \chi (y') f (y',z)|  
&\le 
|\chi (y)| | f (y,z) - f (y',z)|  
+
|\chi (y) - \chi (y')|| f (y',z)|  
\nn\\
&\le 
\{c_2 (R^X)^r + \|f\|_{\infty}\} |y-y'|_{\cW}
\nn\\
&\le c_3 \{(R^X)^r +1\}|y-y'|_{\cW}.
\end{align*}
It is obvious that 
$\chi \tilde{f}$ and $\chi g$ also satisfy the same property.

Noting that $Y$ (resp. $\tilde{Y}$) solves RDEs 
 \eqref{rde.0413} (resp. \eqref{rde.0506}) whose drift 
 vector field is replaced by $\chi f$ (resp. $\chi \tilde{f}$).
Hence, we can use Inequality \eqref{in.0506-2} in Proposition \ref{prop.0506}  
with 
\[
K' \le c_3 \{(R^X)^r +1\} \le c_4 ( \vertiii{X}_{\alpha}^{1/\beta} +1).
 \]
This proves our assertion.
 \end{proof}

\begin{remark}
In a recent preprint \cite{duc},  well-posedness of an RDE
with globally Lipschitz drift is proved (in the framework of CP theory). 
If one uses this result, one may be able to drop the 
boundedness assumption on $f$ in this section 
and in our main theorem (Theorem \ref{thm.main}).
Though this problem looks interesting, 
we do not pursue it in this paper.
\end{remark}

%

\section{Slow-fast system of rough differential equations}
\label{sec.SF}

In this section we define the slow-fast system \eqref{def.SFeq} 
of RDEs precisely and study it in details. 
In this and the next sections, 
we always assume that $\sigma$ and $h$ are of $C^3$ and 
$f$ and $g$ are locally Lipschitz continuous. 
Under these assumptions,  the slow-fast system always has a
unique time-local solution.

\subsection{Deterministic aspects}
As before, $1/3< \beta< \alpha \le 1/2$ is assumed.
Respecting the direct sum decomposition 
$\R^{d+e}= \R^{d}\oplus \R^{e}$,
a generic path taking values in this space is denoted by
 $(b, w)=(b_t, w_t)_{0 \le t \le T}$.
Similarly, a generic element $\Xi = (\Xi^1, \Xi^2)$
of $\Omega_\alpha (\R^{d+e})$  is denoted by
\begin{equation} \label{def.0518-1}
\Xi^{1}=\left(B^{1}, W^{1}\right), 
\qquad
\Xi^{2}=\left(\begin{array}{cc}
B^{2} & I[B, W] \\
I[W, B] & W^{2}
\end{array}\right).
\end{equation}
It is clear that $B=(B^{1}, B^{2}) \in \Omega_\alpha (\R^{d})$
and  $W=(W^{1}, W^{2}) \in \Omega_\alpha (\R^{e})$.
$I[B, W]$ takes values in $\R^{d}\otimes \R^{e}$ and,
loosely speaking, plays the role of the iterated integral
$(s,t) \mapsto \int_s^t B^{1}_{s,u} \otimes d_u W^{1}_{s,u}$.
Also, $I[W,B]$ can be explained in an analogous way.

In the same manner, respecting the direct sum decomposition 
$\R^{m+n}= \R^{m}\oplus \R^{n}$, 
a generic element of $\R^{m+n}$ is denoted by $z =(x,y)$.
We set 
\[
F_{\varepsilon}(x, y)=\left(\begin{array}{c}f(x, y) 
\\ 
\varepsilon^{-1} g(x, y)\end{array}\right), 
\qquad
 \Sigma_{\varepsilon} (x, y)=\left(\begin{array}{cc}\sigma(x) & O 
 \\ 
 O & \varepsilon^{-1/2} h(x, y)\end{array}\right).
\]
Then, $F_{\varepsilon} \colon \R^{m+n}\to \R^{m+n}$ and 
$\Sigma_{\varepsilon} \colon  \R^{m+n}\to L( \R^{d+e},  \R^{m+n})$.

Let $(Z, Z^\dagger)\in \cQ^{\alpha}_\Xi (\R^{m+n})$ be a 
CP with respect to $\Xi \in \Omega_\alpha (\R^{d+e})$.
We write $Z =(X,Y)$.
Since $Z^\dagger$ takes values in  $L(\R^{d+e}, \R^{m+n})$, 
it can be written as a block matrix;
\[
Z^\dagger =  
\left(\begin{array}{cc}
Z^{\dagger, 11} &  Z^{\dagger, 12}
 \\ 
Z^{\dagger, 21} &  Z^{\dagger, 22}
 \end{array}\right).
 \]
Using this notation, we can write 
the remainder part $Z^\sharp$ as follows;
\[
Z^\sharp_{s,t} 
=
 \left(\begin{array}{c} X^1_{s,t} 
\\ 
Y^1_{s,t} 
\end{array}\right)
-
\left(\begin{array}{cc}
Z^{\dagger, 11}_s &  Z^{\dagger, 12}_s
 \\ 
Z^{\dagger, 21}_s &  Z^{\dagger, 22}_s
 \end{array}\right)
  \left(\begin{array}{c} B^1_{s,t} 
\\ 
W^1_{s,t} 
\end{array}\right).
 \]

The precise definition of the slow-fast system \eqref{def.SFeq}
of RDEs (in the deterministic sense) is given by  
\begin{equation}\label{rde.SF}
Z_{t}^{\varepsilon} 
=z_0 +\int_{0}^{t} F_{\varepsilon} (Z^{\varepsilon}_{s}) ds
+\int_{0}^{t} \Sigma_{\varepsilon} (Z^{\varepsilon}_{s}) d \Xi_{s},
\quad
(Z^{\varepsilon})^\dagger_t = \Sigma_{\varepsilon} (Z^{\varepsilon}_t),  \qquad t \in [0,T].
\end{equation}
We consider this RDE in the $\beta$-H\"older topology.

\begin{remark}
An element of a direct sum space  is denoted by
both a ``column vector" and a``row vector." 
These are not precisely distinguished.
\end{remark}

If $Z^{\varepsilon}$ solves the above slow-fast system \eqref{rde.SF}, 
then its slow component $X^\ve$ solves an RDE 
driven by $B$ alone on $\R^m$ 
with $Y^\ve$ viewed as a parameter.
(See \eqref{rde.Fast} below. This type was introduced in \eqref{rde.0413}
and  discussed in Subsection \ref{sec.rde}.)
Since a solution of an RDE is not a usual path, this fact is not obvious.
\begin{lemma}\label{lem.rde.Xve}
Let $(Z^{\varepsilon}, \Sigma_{\varepsilon} (Z^{\varepsilon}) )
\in \cQ^{\beta}_\Xi ([0, \tau], \R^{m+n})$, $0<\tau \le T$,
be  a unique solution of RDE \eqref{rde.SF} on $[0, \tau]$.
Then, $(X^{\varepsilon}, \sigma (X^{\varepsilon}) )$ belongs to
$\cQ^{\beta}_B ([0, \tau], \R^{m})$ and is 
a unique local solution of the following RDE driven by $B$ on $[0, \tau]$:
\begin{equation}\label{rde.Fast}
X_{t}^{\varepsilon} 
=x_0 +\int_{0}^{t} f (X^{\varepsilon}_{s}, Y^{\varepsilon}_{s}) ds
+\int_{0}^{t} \sigma (X^{\varepsilon}_{s}) dB_{s},
\quad
(X^{\varepsilon})^\dagger_t = \sigma(X^{\varepsilon}_t),  \qquad t \in [0,T].
\end{equation}
\end{lemma}

\begin{proof}
Applying Proposition \ref{prop.gousei} to the right hand side of
RDE \eqref{rde.SF}, we see that
\[
Z^{\varepsilon, 1}_{s,t} - \Sigma_{\varepsilon} (Z^{\varepsilon}_s)\Xi^1_{s,t}
=
 \left(\begin{array}{c} X^{\varepsilon,1}_{s,t} 
\\ 
Y^{\varepsilon, 1}_{s,t} 
\end{array}\right)
-
\left(\begin{array}{cc}\sigma( X^\varepsilon_s ) & O 
 \\ 
 O & \varepsilon^{-1/2} h(X^\varepsilon_s, Y^\varepsilon_s)
 \end{array}\right)
  \left(\begin{array}{c} B^1_{s,t} 
\\ 
W^1_{s,t} 
\end{array}\right)
 \]
 is of $2\beta$-H\"older as a function of $(s,t) \in \triangle_\tau$
 and so is its first component $X^{\varepsilon,1}_{s,t}  -\sigma( X^\varepsilon_s ) B^1_{s,t}$. This means 
 $(X^{\varepsilon}, \sigma (X^{\varepsilon}) ) \in
\cQ^{\beta}_B ([0, \tau], \R^{m})$.

It suffices to show that the first component of
$\int_{0}^{t} \Sigma_{\varepsilon} (Z^{\varepsilon}_{s}) d \Xi_{s}$
equals 
$\int_{0}^{t} \sigma (X^{\varepsilon}_{s}) d B_{s}$.
A straightforward but slightly cumbersome computation 
of block matrices yields that 
\begin{align} 
\Sigma_{\varepsilon} (Z^{\varepsilon}_{s}) \Xi^1_{s,t}
+
\{\Sigma_{\varepsilon} (Z^{\varepsilon})\}^{\dagger}_{s} \Xi^2_{s,t}
=
\left(\begin{array}{c} 
\sigma( X^\varepsilon_s )B^1_{s,t} 
+\nabla_x \sigma( X^\varepsilon_s ) \sigma( X^\varepsilon_s ) 
\langle B^2_{s,t}  \rangle
\\ 
\varepsilon^{-1/2} h(X^\varepsilon_s, Y^\varepsilon_s) W^1_{s,t} 
\,\,+ \,\, \star
\end{array}\right),
\label{eq.smmnd}
\end{align}
where we set
\begin{equation}
\star = \varepsilon^{-1/2} \nabla_x h(X^\varepsilon_s, Y^\varepsilon_s)
\sigma( X^\varepsilon_s ) 
\langle I[B,W]_{s,t}  \rangle
+
\varepsilon^{-1} \nabla_y h(X^\varepsilon_s, Y^\varepsilon_s)
h(X^\varepsilon_s, Y^\varepsilon_s)
\langle W^2_{s,t}  \rangle.
\label{eq.smmnd2}
\end{equation}
In above, $\nabla_x$ and $\nabla_y$ stands for the partial gradient operator 
in $x$ and $y$, respectively.
The precise meaning of $\nabla_x \sigma(x) \sigma(x)$ was already essentially
explained in \eqref{def.nabsigY}. Also, 
$\nabla_x h(x, y)\sigma(x)$ and $\nabla_y h(x, y)h(x, y)$
should be understood in a similar way.

The left hand side of  \eqref{eq.smmnd} is a summand
in the modified Riemann sum that approximates 
$\int_{a}^{b} \Sigma_{\varepsilon} (Z^{\varepsilon}_{s}) d \Xi_{s}$,
$(a,b) \in \triangle_\tau$.
The first component of the right hand side of \eqref{eq.smmnd} equals
a summand
in the modified Riemann sum that approximates 
$\int_{a}^{b} \sigma (X^{\varepsilon}_{s}) d B_{s}$.
This verifies our assertion.
\end{proof}

Now we calculate the It\^o-Stratonovich correction
at a deterministic level.
Set 
${\rm Id}_e := \sum_{i=1}^e  \mathbf{a}_i \otimes  \mathbf{a}_i$
for an orthonormal basis $\{ \mathbf{a}_i\}_{i=1}^e$ of $\R^e$.
Note that this definition is independent of the choice of 
$\{ \mathbf{a}_i\}$.
For  $\lambda\in\R$, we set 
\begin{align*}
\tilde{g} (x,y) &
=g (x,y)  - \lambda 
\nabla_y h(x,y)
h(x,y)
\langle  {\rm Id}_e \rangle
\\
&=
g (x,y) 
- \lambda {\rm Trace}\bigl[
\nabla_y h(x,y)
h(x,y)
\langle \bullet, \star   \rangle \bigr]
\end{align*}
and $\tilde{F}_{\varepsilon}(x, y)= (f (x,y), \ve^{-1} \tilde{g} (x,y))^T$.
We also define $\tilde\Xi = (\tilde{\Xi}^1, \tilde{\Xi}^2)$ by 
$\tilde{\Xi}^1= \Xi^1$ and 
\[
\tilde{\Xi}^{2}_{s,t}=\left(\begin{array}{cc}
B^{2}_{s,t} & I[B, W]_{s,t} \\
I[W, B]_{s,t} & \tilde{W}^{2}_{s,t} 
\end{array}\right),
\quad \mbox{where}\quad
\tilde{W}^{2}_{s,t} :=
W^{2}_{s,t} +\lambda\, {\rm Id}_e(t-s).
\]
Since the definition of CP
depends only on the first level of the reference RP, 
$(Y, Y^\dagger)$ belongs to $\cQ^{\beta}_\Xi ([a,b], \cW)$
if and only if it belongs to $\cQ^{\beta}_{\tilde\Xi} ([a,b], \cW)$.

\begin{lemma} \label{lem.det.IS}
Let the notation be as above and let $\tau \in (0, T]$. 
Then, the following are equivalent:
\begin{enumerate} 
\item[{\rm (1)}]
$(Z^{\varepsilon}, \Sigma_{\varepsilon} (Z^{\varepsilon}) )
\in \cQ^{\beta}_\Xi ([0, \tau], \R^{m+n})$ and it
solves RDE \eqref{rde.SF} on $[0, \tau]$.
\item[{\rm (2)}]
$(Z^{\varepsilon}, \Sigma_{\varepsilon} (Z^{\varepsilon}) )
\in \cQ^{\beta}_{\tilde\Xi} ([0, \tau], \R^{m+n})$ and it
solves the following RDE on $[0, \tau]$.
\[
Z_{t}^{\varepsilon} 
=z_0 +\int_{0}^{t} \tilde{F}_{\varepsilon} (Z^{\varepsilon}_{s}) ds
+\int_{0}^{t} \Sigma_{\varepsilon} (Z^{\varepsilon}_{s}) d \tilde\Xi_{s},
\quad
(Z^{\varepsilon})^\dagger_t = \Sigma(Z^{\varepsilon}_t).
\]
\end{enumerate}
\end{lemma}

\begin{proof} 
We only prove (1) implies (2). 
(The proof of the converse is essentially the same.)
The second term on the right hand side of \eqref{eq.smmnd2}
equals
\[
\varepsilon^{-1} \nabla_y h(X^\varepsilon_s, Y^\varepsilon_s)
h(X^\varepsilon_s, Y^\varepsilon_s)
\langle \tilde{W}^2_{s,t}  \rangle
+
\varepsilon^{-1} \nabla_y h(X^\varepsilon_s, Y^\varepsilon_s)
h(X^\varepsilon_s, Y^\varepsilon_s)
\langle {\rm Id}_e\rangle (t-s).
\]
The second term above is a summand in the Riemann sum 
which approximates 
\[
\ve^{-1}\int_0^{\cdot}   \nabla_y h(X^\varepsilon_s, Y^\varepsilon_s)
h(X^\varepsilon_s, Y^\varepsilon_s)
\langle {\rm Id}_e\rangle ds.
\]
Thus, we have proved that  (1) implies (2).
\end{proof}

\begin{remark}\label{rem.IS.corr}
If $\tilde{W}^2$ and $W^2$ are the second level paths of 
Stratonovich and It\^o Brownian RP, respectively, then
it is well-known that $\tilde{W}^2_{s,t}=W^2_{s,t} + (1/2) {\rm Id}_e(t-s)$
and Lemma \ref{lem.det.IS} can be used with $\lambda =1/2$.
Therefore, we can also use Stratonovich Brownian RP 
to formulate the slow-fast system \eqref{def.SFeq} of random RDEs.
(Indeed, in the previous paper \cite{pix2}, the Stratonovich-type
formulation is used.)
  \end{remark}  

%

\subsection{Probabilistic aspects}

Let $\tfrac13 <\alpha_0 \le \tfrac12$ and 
$(\Omega, \mathcal{F}, {\mathbb P})$ be a probability space. 
Let $w =(w_t)_{0\le t \le T}$ be a standard $e$-dimensional Brownian
motion 
and let $B=\{(B^1_{s,t}, B^2_{s,t})\}_{0\le s\le t\le T}$ be an 
$\Omega_{\alpha} (\R^d)$-valued 
random variable (i.e., random RP)
defined on $(\Omega, \mathcal{F}, {\mathbb P})$
for every $\alpha \in (1/3,\alpha_0)$.
We assume that $w$ and $B$ are independent.
As for the integrability of $B$, Condition ${\bf (A)}$ is assumed.
Let $\{\cF_t\}_{0\le t\le T}$ be a filtration satisfying the usual condition
as well as the following two conditions:
 {\rm (i)} $w$ is an $\{\cF_t\}$-BM and 
{\rm (ii)} $t \mapsto (B^1_{0,t}, B^2_{0,t})$ is $\{\cF_t\}$-adapted.

Our random RP $\Xi$ is defined in the following way.
The symbols for 
the components of $\Xi$ were set in \eqref{def.0518-1}.
\begin{definition}\label{def.0520}
$B=(B^1, B^2)$ is given a priori in our assumptions.
First, define $W =(W^1, W^2)$ as the It\^o-type Brownian RP, that is,
\[
W^1_{s,t} = w_t-w_s, \qquad  W^2_{s,t} = \int_s^t (w_u -w_s)\otimes d^{{\rm I}}w_u,
\]
where $d^{{\rm I}}w$ stands for the standard It\^o integral.
It is well-known that $W \in \Omega_{(1/2) -\delta} (\R^e)$ 
for every (suffiiently small) $\delta >0$.
Next, we  set
\begin{equation} \label{in.0519-1}
 I[B, W]_{s,t} = \int_s^t B^1_{s,u}\otimes d^{{\rm I}}w_u,
 \qquad
 I[W, B]_{s,t} =  W^1_{s,u}\otimes B^1_{s,u} 
   - \int_s^t ( d^{{\rm I}}w_u)\otimes B^1_{s,u}.
 \end{equation}
Thus, we have set all components of $\Xi =(\Xi^1, \Xi^2)$.
\end{definition}

\begin{lemma} 
Suppose $\alpha \in (1/3,\alpha_0)$
and let $\Xi$ be as in Definition \ref{def.0520}.
Under ${\bf (A)}$ 
it holds that $\Xi \in \Omega_\alpha (\R^{d+e})$, a.s.
Moreover, $\E [\vertiii{\Xi}_{\alpha}^q  ] <\infty$
for every $q \in [1,\infty)$.
\end{lemma}

\begin{proof} 
We use a RP-version of Kolmogorov-\v{C}entsov's continuity criterion 
(see \cite[Theorem 3.1]{fh}).
Take any $\alpha' \in (\alpha,\alpha_0)$.
By the criterion it suffices to show that 
for every sufficiently large $q\ge 1/\alpha'$, there exists a constant $C=C_q>0$
such that
\begin{equation} \label{in.0520-1}
\E[ |\Xi^1_{s,t}|^{q}]^{1/q} \le C (t-s)^{\alpha'},
\qquad
\E[ |\Xi^2_{s,t}|^{q/2}]^{2/q} \le C (t-s)^{2\alpha'},
\quad
(s,t) \in \triangle_T
\end{equation}
holds. We will verify \eqref{in.0520-1} componentwise.

Due to ${\bf (A)}$, $B=(B^1, B^2)$ clearly satisfies \eqref{in.0520-1}.
It is well-known that Brownian RP $W=(W^1, W^2)$
satisfies \eqref{in.0520-1}, too.
Now, we estimate the $(i,j)$-component of $I[B, W]$
($1\le i \le d$, $1\le j \le e$).
By Burkholder's inequality, we see that
\begin{align*} 
\E \Bigl[  \Bigl| \int_s^t B^{1,i}_{s,u} d^{{\rm I}}w_u^j \Bigr|^q \Bigr]
&\le 
C_1 
\E \Bigl[  \Bigl| \int_s^t (B^{1,i}_{s,u})^2  du \Bigr|^{q/2} \Bigr]
\nn\\
&\le 
C_1 \E [ \| B^1\|_{\alpha'}^q]  
 \Bigl| \int_s^t (u-s)^{2\alpha'}du  \Bigr|^{q/2}
 \le 
 C_2 (u-s)^{(1+2\alpha' )q/2}
 \end{align*}
 for certain positive constants $C_1, C_2$ independent of $(s,t)$.
 Hence, $I[B, W]$ satisfies \eqref{in.0520-1}.
Due to \eqref{in.0519-1}, the proof for $I[W, B]$
 is essentially the same.
    \end{proof}


In the sequel we work under ${\bf (A)}$ and
 assume that $1/3< \beta< \alpha <\alpha_0 (\le 1/2)$.
For the rest of this section,
$\Xi$ is as in Definition \ref{def.0520}.  
The precise meaning of the random RDE in our main theorem is
 RDE \eqref{rde.SF} driven by this $\Xi$.

We extend the time interval of the filtration 
$\{\mathcal{F}_t\}$
by setting $\mathcal{F}_{t}=\mathcal{F}_{t\wedge T}$ for $t \ge 0$.
Denote by $\hat{\R}^{m+n} := \R^{m+n}\cup \{\infty\}$
the one-point compactification of $\R^{m+n}$.
If a global solution $(Z^{\ve}_t)_{t \in [0,T]}$ exists, 
then we set $Z^{\ve}_t =Z^{\ve}_{t\wedge T}$ for $t \ge 0$.
Otherwise, denote by
$(Z^{\ve}_t)_{t \in [0,u^\ve)}$, $0< u^\ve \le T$, 
be a maximal local solution
and set $Z^{\ve}_t =\infty$ for $t \in [u^\ve, \infty)$.
Either way, $(Z^{\ve}_t)$ is constant in $t$ on $[T,\infty)$, a.s.

Define
$\tau^\ve_N =\inf\{ t \ge 0 \mid  |Z^{\ve}_t|\ge N\}$
for each $N \in \N$
and $\tau^\ve_\infty = \lim_{N\to \infty}\tau^\ve_N$.
(As usual $\inf \emptyset :=\infty$.)
These are $\{\mathcal{F}_t\}$-stopping times. 
Then, the following are equivalent:

\begin{itemize} 
\item
A global solution $(Z^{\ve}_t)_{t \in [0,T]}$ of RDE \eqref{rde.SF} exists.
\item
$(Z^{\ve}_t)_{t \in [0,T)}$ defined as above is bounded in $\R^{m+n}$.
\item
$\tau^\ve_N =\infty$ for some $N$.
\item
$\tau^\ve_\infty > T$.
\end{itemize}
It should be noted that while a solution of RDE 
\eqref{rde.SF} moves in a bounded set, its trajectory is uniformly 
continuous in $t$ (because its H\"older norm is bounded). 
Hence, if $(Z^{\ve}_t)_{t \in [0,s)}$, $0<s\le T$, is bounded, 
then $(Z^{\ve}_t)_{t \in [0,s]}$ solves  RDE \eqref{rde.SF}.

On the other hand, if no global solution exists, 
then we have $u^\ve = \tau^\ve_\infty \in (0,T]$ and 
$\limsup_{t \nearrow \tau^\ve_\infty} |Z_t| =\infty$.
Moreover, $\lim_{t \nearrow \tau^\ve_\infty} |Z_t| =\infty$
because of the uniform continuity mentioned above.
Therefore, $(Z^{\ve}_t)_{t \ge 0}$ is a continuous process 
that takes values in $\hat{\R}^{m+n}$.


\begin{proposition} \label{pr.0520}
Let the notation be as above
and assume ${\bf (A)}$. Then, for every $\ve \in (0,1]$,
$Y^{\ve}$ satisfies
the following It\^o SDE up to the explosion time of 
$Z^{\ve}=(X^{\ve}, Y^{\ve})$:
\[
Y^{\ve}_t = y_0 + 
\frac{1}{\varepsilon}\int_0^{t \wedge T}
 g(X^\varepsilon_s, Y^\varepsilon_s) ds
 +
\frac{1}{\sqrt{\varepsilon}}\int_0^{t \wedge T}
 h(X^\varepsilon_s, Y^\varepsilon_s) d^{{\rm I}}w_s,
 \qquad
 0 \le t < \tau^{\ve}_{\infty}.
\]
\end{proposition}

\begin{proof} 
Let $\mathcal{P}=\{ 0=t_0 <t_1 <\cdots < t_K =t\}$ be a partition of 
$[0,t]$ for $0<t \le T$.
The summand of the Riemann sum for the RP integral in
RDE \eqref{rde.SF} was given in 
\eqref{eq.smmnd} and \eqref{eq.smmnd2}.

First, we prove the lemma 
when $h, \sigma$ are of $C^3_b$
and $f, g$ are bounded and globally Lipschitz continuous.
In this case the solution never explodes, i.e. 
$\tau_{\infty}^{\varepsilon} =\infty$, a.s.
It is easy to see that 
\begin{equation}\label{eq.0531-1}
\lim_{|\mathcal{P} | \searrow 0}
\sum_{i=1}^K
h(X^\varepsilon_{t_{i-1}}, Y^\varepsilon_{t_{i-1}}) W^1_{t_{i-1},t_i}
=
\int_0^t  h(X^\varepsilon_s, Y^\varepsilon_s)  d^{{\rm I}}w_s  
\qquad
\mbox{in $L^2 ({\mathbb P})$.}
\end{equation}

Denote by $\{  \mathbf{e}_k\}_{k=1}^e$ and 
$\{  \mathbf{f}_p\}_{p=1}^n$ the canonical 
orthonormal basis of $\R^e$ and $\R^n$, respectively.
We write  
$A^{k,l}_s := 
\nabla_y h(X^\varepsilon_s, Y^\varepsilon_s)
h(X^\varepsilon_s, Y^\varepsilon_s)
\langle \mathbf{e}_k \otimes  \mathbf{e}_l \rangle$,
$A^{p;k,l}_s := \langle A^{k,l}_s,  \mathbf{f}_p\rangle_{\R^n}$
and $W^2_{s,t} = \sum_{k, l=1}^e W^{2;k,l}_{s,t} \mathbf{e}_k \otimes  \mathbf{e}_l $.
Then, we obviously see that the $p$th component of
$\nabla_y h(X^\varepsilon_s, Y^\varepsilon_s)
h(X^\varepsilon_s, Y^\varepsilon_s)
\langle W^2_{s,t}  \rangle$ equals $\sum_{k, l=1}^e A^{p;k,l}_s W^{2;k,l}_{s,t}$.
Noting that $\E [W^{2;k,l}_{s,t} | \cF_s ]=0$ for 
all $k,l$ and $s \le t$, we have for all $i<j$ and $p$ that
\begin{align} 
\lefteqn{
\E \Bigl[ 
\sum_{k, l=1}^e A^{p;k,l}_{t_{i-1}} W^{2;k,l}_{t_{i-1},t_i}
\cdot
\sum_{k', l' =1}^e A^{p;k',l'}_{t_{j-1}} W^{2;k',l'}_{t_{j-1},t_j}
\Bigr]
}
\label{eq.0531-2}\\
&= \sum_{k, l}\sum_{k', l'}
\E \bigl[
A^{p;k,l}_{t_{i-1}} W^{2;k,l}_{t_{i-1},t_i}
 A^{p;k',l'}_{t_{j-1}} 
 \E[
  W^{2;k',l'}_{t_{j-1},t_j}  | \cF_{t_{j-1}} ] 
\bigr]
=0.
\nn
\end{align}
Using $\E[(W^{2;k,l}_{s,t})^2] = (t-s)^2/2$, we can easily see 
from \eqref{eq.0531-2} that
\begin{align} 
\lefteqn{
\E \Bigl[ 
\Bigl|  
\sum_{i=1}^K 
\nabla_y h(X^\varepsilon_{t_{i-1}}, Y^\varepsilon_{t_{i-1}})
h(X^\varepsilon_{t_{i-1}}, Y^\varepsilon_{t_{i-1}})
\langle W^2_{t_{i-1},t_i}  \rangle
\Bigr|^2
\Bigr]
}
\label{eq.0531-9}
\\
&=\sum_{p,k,l} \sum_{i=1}^K 
\E [ 
|  A^{p;k,l}_{t_{i-1}} W^{2;k,l}_{t_{i-1},t_i}|^2]
\nn\\
&\le 
\frac12 \sum_{p,k,l} \sum_{i=1}^K \|h\|_{\infty}  \|\nabla h\|_{\infty} 
(t_i -t_{i-1})^2
\le 
\frac12 ne^2 \|h\|_{\infty}  \|\nabla h\|_{\infty} |\mathcal{P}| \to 0
\nn
\end{align}
as the mesh size $|\mathcal{P}|$ tends to $0$.

In a similar way, 
we write  
$\tilde{A}^{k,l}_s := 
\nabla_y h(X^\varepsilon_s, Y^\varepsilon_s)
\sigma (X^\varepsilon_s, Y^\varepsilon_s)
\langle \tilde{\mathbf{e}}_k \otimes  \mathbf{e}_l \rangle$,
$\tilde{A}^{p;k,l}_s := \langle \tilde{A}^{k,l}_s,  \mathbf{f}_p\rangle_{\R^n}$
and $I[B,W]_{s,t} =\sum_{k=1}^d
 \sum_{l=1}^e I[B,W]^{k,l}_{s,t} 
\tilde{\mathbf{e}}_k \otimes  \mathbf{e}_l $,
where $\{  \tilde{\mathbf{e}}_k\}_{k=1}^d$ is
the canonical  orthonormal basis of $\R^d$.
Then,  the $p$th component of
$\nabla_y h(X^\varepsilon_s, Y^\varepsilon_s)
\sigma (X^\varepsilon_s, Y^\varepsilon_s)
\langle I[B,W]_{s,t}  \rangle$ equals 
$\sum_{k=1}^d\sum_{l=1}^e\tilde{A}^{p;k,l}_s  I[B,W]^{k,l}_{s,t}$.
Since $I[B,W]_{s,t}$ is an It\^o integral, 
$\E [I[B,W]^{k,l}_{s,t} | \cF_s ]=0$ for all $k,l$ and $s \le t$.
Hence, in the same way as in \eqref{eq.0531-2} we have
\[
\E \Bigl[ 
\sum_{k, l} \tilde{A}^{p;k,l}_{t_{i-1}} I[B,W]^{k,l}_{t_{i-1},t_i}
\cdot
\sum_{k', l'} \tilde{A}^{p;k',l'}_{t_{j-1}} I[B,W]^{k',l'}_{t_{j-1},t_j}
\Bigr]
=0
\]
for all $i<j$ and $p$.
Since $\E[(I[B,W]^{k,l}_{s,t})^2] =(t-s)^{2\alpha' +1}/ (2\alpha' +1)$,
we can calculate in the same way as in \eqref{eq.0531-9} to see that
\begin{align} 
\lefteqn{
\E \Bigl[ 
\Bigl|  
\sum_{i=1}^K 
\nabla_y h(X^\varepsilon_{t_{i-1}}, Y^\varepsilon_{t_{i-1}})
\sigma (X^\varepsilon_{t_{i-1}}, Y^\varepsilon_{t_{i-1}})
\langle I[B,W]_{t_{i-1},t_i}  \rangle
\Bigr|^2
\Bigr]
}
\label{eq.0531-4}\\
&=\sum_{p,k,l} \sum_{i=1}^K 
\E [ 
|  \tilde{A}^{p;k,l}_{t_{i-1}} I[B,W]^{k,l}_{t_{i-1},t_i}|^2]
\nn\\
&\le 
\sum_{p,k,l} \sum_{i=1}^K \|\sigma\|_{\infty}  \|\nabla h\|_{\infty} 
\frac{(t_i -t_{i-1})^{2\alpha' +1}}{2\alpha'  +1}
\le 
\frac{ne^2}{2\alpha' +1} \|\sigma\|_{\infty}  \|\nabla h\|_{\infty} |\mathcal{P}|^{2\alpha'} \to 0
\nn
\end{align}
as the mesh size $|\mathcal{P}|$ tends to $0$.

Combining \eqref{eq.0531-1}, \eqref{eq.0531-9} and \eqref{eq.0531-4},
we have shown that 
\begin{equation}\label{eq.0531-5}
Y^{\ve}_t = y_0 + 
\frac{1}{\varepsilon}\int_0^{t \wedge T}
 g(X^\varepsilon_s, Y^\varepsilon_s) ds
 +
\frac{1}{\sqrt{\varepsilon}}\int_0^{t \wedge T}
 h(X^\varepsilon_s, Y^\varepsilon_s) d^{{\rm I}}w_s,
 \qquad
t \ge 0
\end{equation}
holds almost surely in this case.

From here only the standing assumption is assumed 
on the coefficients $h, \sigma, f, g$.
Take any sufficiently large $N$.  
Let $\phi_N \colon \R^{m+n}\to [0,1]$ be a smooth function 
with compact support such that $\phi_N \equiv 1$ on 
the ball $\{ z\in \R^{m+n}\mid |z| \le N \}$
and set $\hat{h} := h \phi_N$. Also,
$\hat{\sigma}, \hat{f}, \hat{g}$ are defined in the same way.
We replace the coefficients of 
RDE \eqref{rde.SF} by these corresponding data with ``hat" and 
denote a unique solution by $\hat{Z}^\ve =(\hat{X}^\ve, \hat{Y}^\ve)$.
Then, \eqref{eq.0531-5} holds with 
$\hat{X}^\ve, \hat{Y}^\ve,  \hat{h}, \hat{g}$ 
in place of $X^\ve, Y^\ve, h, g$.
By the uniqueness of the RDE, it holds that 
$\hat{Z}^{\ve}_{t \wedge \tau^\ve_N} = Z^\ve_{t \wedge \tau^\ve_N}$
for all $0 \le t \le T$.
Therefore, we almost surely have 
\begin{align} 
Y^\ve_{t \wedge \tau^\ve_N  \wedge T}
&=
\hat{Y}^\ve_{t \wedge \tau^\ve_N  \wedge T}
\label{eq.0601-1}\\
&=
 y_0 + 
\frac{1}{\varepsilon}\int_0^{t \wedge \tau^\ve_N \wedge T}
 \hat{g}(\hat{X}^\varepsilon_s, \hat{Y}^\varepsilon_s) ds
 +
\frac{1}{\sqrt{\varepsilon}}\int_0^{t \wedge \tau^\ve_N \wedge T}
 \hat{h} (\hat{X}^\varepsilon_s, \hat{Y}^\varepsilon_s) d^{{\rm I}}w_s,
\nn\\
&=
 y_0 + 
\frac{1}{\varepsilon}
\int_0^{t \wedge \tau^\ve_N \wedge T}
g (X^\varepsilon_s, Y^\varepsilon_s) ds
 +
\frac{1}{\sqrt{\varepsilon}}
\int_0^{t \wedge \tau^\ve_N\wedge T}
 h(X^\varepsilon_s, Y^\varepsilon_s) d^{{\rm I}}w_s,
  \quad
t \ge 0.
\nn
\end{align}
Since $\tau^\ve_N \nearrow \tau^\ve_\infty$ as $N \to \infty$ a.s. 
on  the set $\{  \tau^\ve_\infty \le T\}$, 
we  finish the proof by letting $N \to \infty$.
\end{proof}

\begin{proposition} \label{prop.0601}
Assume ${\bf (A)}$, ${\bf (H1)}$--${\bf (H4)}$
and  ${\bf (H6)}_q$ for some $q\ge 2$.
Then, the probability that 
$Z^{\ve}=(X^{\ve}, Y^{\ve})$ explodes on $[0,T]$ is zero. 
Moreover, we have 
\begin{align}
\sup_{0<\ve \le 1}
\E [\|X^{\ve} \|_{\beta, [0,T]}^p] &<\infty, \qquad 1\le p<\infty,
\label{in.0602-1}
\\
\sup_{0<\ve \le 1}\sup_{0\le t \le T}\E [ |Y^{\ve}_t|^q]  &<\infty.
\label{in.0602-2}
\end{align}
\end{proposition}

\begin{proof} 
In this proof, $c_i~(i \ge 1)$  denotes a positive constant  
independent of $\ve, N$, $s, t$ and 
the sample $\omega \in \Omega$.

Set $\lambda^\ve_N =\inf\{ t \ge 0 \mid  |Y^{\ve}_t|\ge N\}$
for each $N \in \N$
and $\lambda^\ve_\infty = \lim_{N\to \infty}\lambda^\ve_N$.
Clearly, $\lambda^\ve_N \ge \tau^\ve_N$, a.s.
We first show $\lambda^\ve_\infty =\tau^\ve_\infty$, a.s.
(In other words, $X^{\ve}$ does not explode before $Y^{\ve}$ does.)
Suppose that $\lambda^\ve_\infty   > \tau^\ve_\infty$ 
holds for some RP $\Xi$,
which automatically implies $\tau^\ve_\infty \le T$
and $Y^{\ve}$ stay bounded on $[0,\tau^\ve_\infty)$.
Due to Lemma \ref{lem.rde.Xve} and Proposition \ref{prop.0702},
$X^{\ve}$ also stay bounded on $[0, \tau^\ve_\infty)$
and hence so does $Z^{\ve}$.
However, this is a contradiction.

Lemma \ref{lem.rde.Xve} and Proposition \ref{prop.0702}
also imply the following.
If $\tau^\ve_\infty \le T$ (i.e. $Z^{\ve}$ explodes),
then $X^{\ve}$ stay bounded on $[0,\tau^\ve_\infty)$ and hence
$\lim_{t\nearrow \tau^\ve_\infty} |Y^{\ve}_t| =\infty$.
In particular, there exists $c_1, c_2 >0$ such that
$\|X^{\ve}\|_{\beta, [0, \tau^\ve_\infty)} 
\le c_1 (\vertiii{B}_{\alpha}^{c_2} +1)$.

For a while we use ${\bf (H6)}_2$ instead of ${\bf (H6)}_q$.
First, we prove non-explosion for any fixed $\ve$
by using Proposition \ref{pr.0520}, in particular \eqref{eq.0601-1}.
Lemma \ref{lem.apen.ito} (with $q=2$) and ${\bf (H6)}_2$ yield that
\begin{align} 
\E [|Y^\ve_{t \wedge \tau^\ve_N}|^2]
&= |y_0|^2 +
\ve^{-1}\E \Bigl[
\int_0^{t \wedge \tau^\ve_N} 
\{ 2\langle Y^\varepsilon_s, g(X^\varepsilon_s, Y^\varepsilon_s) \rangle
 + |h(X^\varepsilon_s, Y^\varepsilon_s) |^2 \}ds
 \Bigr]
\label{in.0602-3}\\
&\le 
|y_0|^2  
+\ve^{-1} \E \Bigl[ \int_0^t
\{  -\gamma_1 |Y^\varepsilon_s|^2 + C |X^\varepsilon_s|^{\eta_3} +C\}
 {\bf 1}_{\{s \le \tau^\ve_N\}}ds \Bigr]
\nn\\
&\le  
|y_0|^2   +\ve^{-1} CT 
\{ c_1^{\eta_3}   \E [(\vertiii{B}_{\alpha}^{c_2} +1)^{\eta_3}] +1\}
\nn\\
&\le  
|y_0|^2   +\ve^{-1}c_3,
\qquad  t \in [0,T].
\nn
\end{align}
Since we stopped the time by $\tau^\ve_N$,
the martingale part in the It\^o formula is a true martingale
and its expectation vanishes.
Recall that, as $N\to\infty$,
$|Y^\varepsilon_{T \wedge \tau^\ve_N}| \to \infty$ 
if 
$\tau^\ve_\infty \le T$ 
and 
$|Y^\varepsilon_{T \wedge \tau^\ve_N}| \to |Y^\varepsilon_T|$ if 
$\tau^\ve_\infty =\infty$.
Applying Fatou's lemma to \eqref{in.0602-3} with $t=T$,
we see that $\{ \tau^\ve_\infty \le T\}$ is a zero set,
i.e. explosion never occurs.
From this \eqref{in.0602-1} immediately follows.
Again by Fatou's lemma, $\sup_{t\le T} 
\E [|Y^\ve_{t}|^2] \le |y_0|^2   +\ve^{-1}c_3$ holds.

In a similar way as above, we use Lemma \ref{lem.apen.ito} 
again to obtain that, for every $p \ge 2$,
\begin{align} 
\E [|Y^\ve_{t}|^{p}  {\bf 1}_{\{t \le \tau^\ve_N \}}]
&\le
\E [|Y^\ve_{t \wedge \tau^\ve_N}|^{p}]
\label{in.0602-4}
\\
&=
|y_0|^{p} +\ve^{-1}\E \Bigl[ \int_0^{t \wedge \tau^\ve_N} 
p \bigl\{  | Y^\varepsilon_s |^{p-2}
 \langle Y^\ve_s,  g(X^\varepsilon_s, Y^\varepsilon_s)\rangle 
 \nn\\
 &\qquad
 +\frac12 |Y^\ve_s|^{p-2} |h(X^\varepsilon_s, Y^\varepsilon_s)|^2
  \nn\\
 & \qquad 
+
\frac{p-2}{2} |Y^\ve_s|^{p-4} \cdot (Y^\ve_s)^\top h(X^\varepsilon_s, Y^\varepsilon_s) h(X^\varepsilon_s, Y^\varepsilon_s)^\top Y^\ve_s 
\bigr\}ds 
\Bigr] %
\nn\\
&\le 
|y_0|^{p} +\ve^{-1}\E \Bigl[ \int_0^{t \wedge \tau^\ve_N} 
\frac{p}{2}  | Y^\varepsilon_s |^{p-2}
\nn\\
&\qquad\qquad \times\{
2\langle Y^\varepsilon_s, g(X^\varepsilon_s, Y^\varepsilon_s) \rangle
+ (p-1) |h(X^{\ve}_s,Y^{\ve}_s)|^2
\}
ds
\Bigr]
\nn\\
&\le 
|y_0|^{p} +
\ve^{-1}\E \Bigl[ \int_0^{t \wedge \tau^\ve_N} 
 c_4 | Y^\varepsilon_s |^{p-2} 
 \{ | Y^\varepsilon_s |^2 + |  X^\varepsilon_s|^{\eta_3 \vee 2} +1
  \}ds  \Bigr]
\nn\\
&\le 
|y_0|^{p} +
\ve^{-1}\E \Bigl[ \int_0^{t \wedge \tau^\ve_N} 
 c_5  \{ | Y^\varepsilon_s |^{p} + |  X^\varepsilon_s|^{(\eta_3 \vee 2)p/2} +1
  \}ds  \Bigr]
\nn\\
&\le 
|y_0|^{p} +
\ve^{-1}c_5  \int_0^{t } 
\E [|Y^\ve_{s}|^{p}  {\bf 1}_{\{s \le \tau^\ve_N \}}]ds  +\ve^{-1} c_6,
\quad
t \in [0,T].
\nn
\end{align}
Here, we have also used ${\bf (H3)}$, ${\bf (H6)}_2$, \eqref{in.0602-1}
and Young's inequality.

Using Lemma \ref{lem.Gron} (1) (Gronwall's inequality),
we have
\[
\E [|Y^\ve_{t}|^{p}  {\bf 1}_{\{t \le \tau^\ve_N \}}]
\le
(|y_0|^{p} + \ve^{-1} c_6) e^{\ve^{-1}c_5 T},
\qquad  t \in [0,T].
\]
Putting this inequality on the right hand side of \eqref{in.0602-4},
we have  
\[
\E [|Y^\ve_{t \wedge \tau^\ve_N}|^{p} ]
\le
C_{\ve},
\qquad  t \in [0,T],
\]
where $C_{\ve}>0$ is a constant which depends on $\ve$,
but not on $N, t$.
Since $p$ can be arbitrarily large, 
this implies that $\{ |Y^\ve_{t \wedge \tau^\ve_N}|^p \}_{N=1}^\infty$
is uniformly integrable.
Since $ |Y^\ve_{t \wedge \tau^\ve_N}| \to |Y^\ve_{t}|$
as $N\to\infty$, a.s., it follows that 
$\E [|Y^\ve_{t \wedge \tau^\ve_N}|^{p}] 
   \to \E [|Y^\ve_{t}|^{p}]$ as $N\to\infty$ and 
\begin{equation}\label{in.0603-2}
\sup_{0\le t\le T} \E [|Y^\ve_{t}|^{p} ]\ <\infty
\qquad
\mbox{for every $\ve \in (0,1]$ and $2\le p<\infty$.}
\end{equation}
%
 %
 
 On the other hand, it  immediately follows from ${\bf (H3)}$, 
${\bf (H4)}$ and ${\bf (H6)}_2$ that
\begin{align*}
\lefteqn{
\Bigl|   | Y^\varepsilon_s |^{p-2}
 \langle Y^\ve_s,  g(X^\varepsilon_s, Y^\varepsilon_s)\rangle 
 +\frac12 |Y^\ve_s|^{p-2} |h(X^\varepsilon_s, Y^\varepsilon_s)|^2
 }
\nn\\
&\qquad +
\frac{p-2}{2} |Y^\ve_s|^{p-4} \cdot (Y^\ve_s)^\top h(X^\varepsilon_s, Y^\varepsilon_s) h(X^\varepsilon_s, Y^\varepsilon_s)^\top Y^\ve_s 
\Bigr|
  {\bf 1}_{\{s \le \tau^\ve_N \}}
\\
&\le   |Y^\ve_s|^{p-1}  C ( |X^\ve_s|^{\eta_1} +|Y^\ve_s|^{\eta_1} +1)
+ 
\frac{p-1}{2} |Y^\ve_s|^{p-2} ( |X^\ve_s| +|Y^\ve_s| +|h(0,0)|)^2
  \nn\\
  &\le 
  c_7 ( |Y^\varepsilon_s|^{c_7} +  |X^\varepsilon_s|^{c_7} +1).
   \end{align*}
The right hand side is independent of $N$ and
integrable on $\Omega \times [0,t]$
with respect to the measure ${\mathbb P} \otimes ds$,
 thanks to \eqref{in.0603-2}.
Hence, we can use the dominated convergence theorem 
to \eqref{in.0602-4} to obtain that
\begin{align*}
\E [|Y^\ve_{t}|^{p}] 
&=
\lim_{N\to\infty}
\E [|Y^\ve_{t \wedge \tau^\ve_N}|^{p}] 
\\
&=
|y_0|^p+ \frac{1}{\ve} \lim_{N\to\infty} 
\E \Bigl[
\int_0^{t \wedge \tau^\ve_N}  
\{ p  | Y^\varepsilon_s |^{p-2}
 \langle Y^\ve_s,  g(X^\varepsilon_s, Y^\varepsilon_s)\rangle 
 +\frac12 |Y^\ve_s|^{p-2} |h(X^\varepsilon_s, Y^\varepsilon_s)|^2
\nn\\
&\qquad \qquad +
\frac{p-2}{2} |Y^\ve_s|^{p-4} \cdot (Y^\ve_s)^\top h(X^\varepsilon_s, Y^\varepsilon_s) h(X^\varepsilon_s, Y^\varepsilon_s)^\top Y^\ve_s 
\} ds
 \Bigr]
 \\
 &=
|y_0|^p  +  \frac{1}{\ve}
\E \Bigl[
\int_0^{t }  
\{ p  | Y^\varepsilon_s |^{p-2}
 \langle Y^\ve_s,  g(X^\varepsilon_s, Y^\varepsilon_s)\rangle 
 +\frac12 |Y^\ve_s|^{p-2} |h(X^\varepsilon_s, Y^\varepsilon_s)|^2
\nn\\
&\qquad \qquad +
\frac{p-2}{2} |Y^\ve_s|^{p-4} \cdot (Y^\ve_s)^\top h(X^\varepsilon_s, Y^\varepsilon_s) h(X^\varepsilon_s, Y^\varepsilon_s)^\top Y^\ve_s 
\}ds \Bigr].
 \end{align*}


After differentiating both sides with  respect to $t$, 
we set $p=q$ and use  ${\bf (H6)}_q$ as follows:
\begin{align*}
\frac{d}{dt}\E [|Y^\ve_{t}|^{q}] 
 &=
 \frac{1}{\ve}\E \Bigl[q  | Y^\varepsilon_t |^{q-2}
 \langle Y^\ve_t,  g(X^\varepsilon_t, Y^\varepsilon_t)\rangle 
 +\frac12 |Y^\ve_t|^{q-2} |h(X^\varepsilon_t, Y^\varepsilon_t)|^2
\nn\\
&\qquad \qquad +
\frac{q-2}{2} |Y^\ve_t|^{q-4} \cdot (Y^\ve_t)^\top h(X^\varepsilon_t, Y^\varepsilon_t) h(X^\varepsilon_t, Y^\varepsilon_t)^\top Y^\ve_t
\Bigr]
 \nn\\
 &\le
 \frac{q}{2\ve}\E \bigl[   | Y^\varepsilon_t |^{q-2}
 \{ 2\langle Y^\varepsilon_t, g(X^\varepsilon_t, Y^\varepsilon_t) \rangle
 + (q-1) |h(X^\varepsilon_t, Y^\varepsilon_t) |^2 \}
 \bigr]
 \nn\\
 &\le
 \frac{q}{2\ve}\E \bigl[   | Y^\varepsilon_t |^{q-2}
 \{ -\gamma_1 | Y^\varepsilon_t |^2 + C|  X^\varepsilon_t|^{\eta_3}+C  \}
 \bigr]
 \nn\\
 &\le 
 - \frac{q\gamma_1}{4\ve}  
     \E[|Y^\varepsilon_t|^q] + 
          \frac{c_8}{\ve}  \E[|  X^\varepsilon_t|^{\eta_3q/2}+1]  
  \nn\\
 &\le
 -  \frac{q\gamma_1}{4\ve}\E[|Y^\varepsilon_t|^q] + \frac{c_9}{\ve},
 \qquad t \in [0,T].
  \end{align*}
  Here, we have also used Young's inequality and \eqref{in.0602-1}.
We can now use Lemma \ref{lem.Gron} (2)
(the differential version of Gronwall's inequality) 
with $b = -q \gamma_1 /(4\ve)$ and $c=c_ 9/\ve$
to obtain that
\[
\E [|Y^\ve_{t}|^{q}] \le |y_0|^q e^{-t q \gamma_1 /(4\ve)} 
+ \frac{4c_9}{q\gamma_1} (1- e^{-t q \gamma_1 /(4\ve)} )
\le 
|y_0|^q + \frac{4c_9}{q\gamma_1}, \qquad  t \in [0,T].
\]
Note that the right hand side is independent of $\ve$ as desired.
Thus, we have obtained \eqref{in.0602-2}. 
This completes the proof of Proposition \ref{prop.0601}.
\end{proof}

%
\section{Proof of main result}

This section is devoted to proving our main result
(Theorem \ref{thm.main}).

First, we introduce a new parameter $\delta$ with $0<\ve <\delta \le 1$.
(In spirit, $0<\ve \ll \delta \ll 1$.
Later, we will set $\delta := \ve^{1/(4\beta)} \log \ve^{-1}$.)
We devide $[0,T]$ into subintervals of equal length $\delta$.
For $s\ge 0$,  set $s(\delta) := \lfloor s/\delta\rfloor \delta$,
which is the nearest breaking point preceding (or equal to) $s$.

We set
\begin{equation}\label{def.0603}
\hat{Y}^{\ve}_t = y_0 + 
\frac{1}{\varepsilon}\int_0^{t }
 g(X^\varepsilon_{s(\delta)}, \hat{Y}^\varepsilon_s) ds
 +
\frac{1}{\sqrt{\varepsilon}}\int_0^{t}
 h(X^\varepsilon_{s(\delta)}, \hat{Y}^\varepsilon_s) d^{{\rm I}}w_s,
 \qquad
t \in [0,T].
\end{equation}
Note that $\hat{Y}^{\ve}$'s dependence on $\delta$ 
is suppressed in the notation.
This approximation process satisfies the following two estimates.

\begin{lemma} \label{lem.0606-1}
Under the same assumptions as in Proposition \ref{prop.0601},
we have the following: 
For every $\delta$ and $\ve$ with $0<\ve <\delta \le 1$, 
the above process $\hat{Y}^{\ve}$ does not explode 
and satisfies
\begin{equation}\label{in.0707-1}
\sup_{0<\ve <\delta \le 1}\sup_{0\le t \le T}
\E [ |\hat{Y}^{\ve}_t|^{q}]  <\infty.
\end{equation}
\end{lemma}

\begin{proof}
The proof is essentially the same as that of Proposition \ref{prop.0601}.
(In fact, this one is easier because we already know 
$X^\varepsilon_{s(\delta)}$ exists and satisfies the estimate 
\eqref{in.0602-1}.)
\end{proof}

\begin{lemma} \label{lem.0603-1}
Assume ${\bf (A)}$, ${\bf (H1)}$--${\bf (H4)}$, 
${\bf (H5)}_r$, ${\bf (H6)}_q$ and ${\bf (H7)}$
for some $q \ge 2$ and $r \ge 0$ such that $q> 2r$.
Then, there exists a positive constant $C$ independent of $\delta$
such that
\[
\sup_{\ve\in (0,\delta)}\sup_{0\le t \le T} 
\E [|Y^\varepsilon_t -  \hat{Y}^\varepsilon_t|^2 ] 
\le C \delta^{2\beta}.
\]
\end{lemma}

\begin{proof} 
By It\^o's  formula, we can easily see that
\begin{align*}
\mathbb{E} [|Y_{t}^{\varepsilon}-\hat{Y}_{t}^{\varepsilon}|^{2}]
&= 
\frac{2}{\varepsilon} \mathbb{E} \Bigl[
\int_{0}^{t}  
 \langle g\left(X_{s}^{\varepsilon}, Y_{s}^{\varepsilon}\right) 
   -g (X_{s(\delta)}^{\varepsilon}, \hat{Y}^{\varepsilon}_s), Y_{s}^{\varepsilon}-\hat{Y}_{s}^{\varepsilon}
   \rangle ds
\Bigr] 
\\
&\quad +\frac{1}{\varepsilon} \mathbb{E} \Bigl[
 \int_{0}^{t} 
   |h (X_{s}^{\varepsilon}, Y_{s}^{\varepsilon})
      - h(X_{s(\delta)}^{\varepsilon}, \hat{Y}_{s}^{\varepsilon})|^{2}
       d s\Bigr].
\end{align*}
Note that, due to \eqref{in.0602-1}, \eqref{in.0603-2}, 
\eqref{in.0707-1}
and ${\bf (H3)}$
the martingale part is actually a true martingale.
(Therefore, stopping times are not needed here.)

By differentiating with  respect to $t$, we get
\begin{align*}
\frac{d}{dt}
\mathbb{E} [|Y_{t}^{\varepsilon}-\hat{Y}_{t}^{\varepsilon}|^{2}]
&= 
\frac{2}{\varepsilon} \mathbb{E} \Bigl[
 \langle g (X_{t}^{\varepsilon}, Y_{t}^{\varepsilon}) 
   -g (X_{t(\delta)}^{\varepsilon}, \hat{Y}^{\varepsilon}_t), Y_{t}^{\varepsilon}-\hat{Y}_{t}^{\varepsilon}
   \rangle 
\Bigr] 
\\
&\quad +\frac{1}{\varepsilon} \mathbb{E} \Bigl[
   |h (X_{t}^{\varepsilon}, Y_{t}^{\varepsilon})
      - h(X_{t(\delta)}^{\varepsilon}, \hat{Y}_{t}^{\varepsilon})|^{2}
       \Bigr]
\\
&=\frac{1}{\varepsilon} \mathbb{E}\Bigl[
\{2\langle g(X_{t}^{\varepsilon}, Y_{t}^{\varepsilon})
-g (X_{t}^{\varepsilon}, \hat{Y}_{t}^{\varepsilon}), 
Y_{t}^{\varepsilon}-\hat{Y}_{t}^{\varepsilon} \rangle
\\
& \qquad \qquad\qquad +
|h(X_{t}^{\varepsilon}, Y_{t}^{\varepsilon})-h(X_{t}^{\varepsilon}, \hat{Y}_{t}^{\varepsilon})|^{2}\} \Bigr]
\\
&\quad + 
\frac{2}{\varepsilon} \mathbb{E} \Bigl[
 \langle g(X_{t}^{\varepsilon}, \hat{Y}_{t}^{\varepsilon}) 
   -g (X_{t(\delta)}^{\varepsilon}, \hat{Y}^{\varepsilon}_t), Y_{t}^{\varepsilon}-\hat{Y}_{t}^{\varepsilon}
   \rangle 
\Bigr] 
\\
&\quad + 
\frac{2}{\varepsilon} \mathbb{E} \Bigl[
 \langle h (X_{t}^{\varepsilon}, Y_{t}^{\varepsilon}) 
   -h (X_{t}^{\varepsilon}, \hat{Y}^{\varepsilon}_t),  \,
       h (X_{t}^{\varepsilon}, \hat{Y}_{t}^{\varepsilon}) 
   -h (X_{t(\delta)}^{\varepsilon}, \hat{Y}^{\varepsilon}_t)   \rangle 
\Bigr] 
\\
&\quad +\frac{1}{\varepsilon} \mathbb{E} \Bigl[
   |h (X_{t}^{\varepsilon}, \hat{Y}_{t}^{\varepsilon})
      - h(X_{t(\delta)}^{\varepsilon}, \hat{Y}_{t}^{\varepsilon})|^{2}
       \Bigr]
       \\
       &=: I_1 +I_2+ I_3+I_4. 
          \end{align*}

We will estimate $I_i~(1\le i \le 4)$.
Below $c_i~(i\ge 1)$ 
are positive constants independent of $t, \ve, \delta$.
We will often use Young's inequality.
It is clear from ${\bf (H7)}$ that
\begin{equation} \label{eq.0606-1}
I_1
 \le 
- \frac{\gamma_2}{\varepsilon} \mathbb{E} [
|Y_{t}^{\varepsilon}-\hat{Y}_{t}^{\varepsilon} |^2].
\end{equation}
Using ${\bf (H3)}$ and \eqref{in.0602-1}, we have
\begin{align} 
I_3+I_4 &\le \frac{c_1}{\ve}  {\rm Lip}(h)^2 \,
  \mathbb{E} [
 |Y_{t}^{\varepsilon}-\hat{Y}_{t}^{\varepsilon} | 
 \cdot 
 |X_{t}^{\varepsilon}- X_{t(\delta)}^{\varepsilon} |
  +
|X_{t}^{\varepsilon}- X_{t(\delta)}^{\varepsilon} |^2]
\label{eq.0606-3}
\\
&\le
 \frac{\gamma_2}{4\varepsilon} \mathbb{E} [
|Y_{t}^{\varepsilon}-\hat{Y}_{t}^{\varepsilon} |^2]
+\frac{c_2}{\ve} \mathbb{E} [
|X_{t}^{\varepsilon}- X_{t(\delta)}^{\varepsilon} |^2]
\nn\\
&\le 
\frac{\gamma_2}{4\varepsilon} \mathbb{E} [
|Y_{t}^{\varepsilon}-\hat{Y}_{t}^{\varepsilon} |^2]
+\frac{c_3}{\ve} \delta^{2\beta},
\nn
\end{align}
where ${\rm Lip}(h)$ stands for the Lipschitz constant of $h$.

Next, 
using ${\bf (H5)}_r$ with $2r <q$, Proposition \ref{prop.0601}
and Lemma \ref{lem.0606-1}, 
we estimate $I_2$ as follows:
\begin{align}
I_2 &\le 
\frac{2C}{\varepsilon}
 \mathbb{E} [
(1+ \| X^\ve\|_\infty^{\eta_2} + |\hat{Y}_{t}^{\varepsilon} |^r) |X_{t}^{\varepsilon}- X_{t(\delta)}^{\varepsilon} |
\cdot
|Y_{t}^{\varepsilon}-\hat{Y}_{t}^{\varepsilon} |
]
\label{eq.0606-2}\\
&\le
 \frac{\gamma_2}{4\varepsilon} \mathbb{E} [
|Y_{t}^{\varepsilon}-\hat{Y}_{t}^{\varepsilon} |^2]
+\frac{c_4}{\ve} \mathbb{E} [
\{ (1+ \| X^\ve\|_\infty)^{2\eta_2}
+ |\hat{Y}_{t}^{\varepsilon} |^{2r}
\} |X_{t}^{\varepsilon}- X_{t(\delta)}^{\varepsilon} |^2]
\nn\\
&\le
 \frac{\gamma_2}{4\varepsilon} \mathbb{E} [
|Y_{t}^{\varepsilon}-\hat{Y}_{t}^{\varepsilon} |^2]
+\frac{c_4}{\ve} \mathbb{E} [
 (1+ \| X^\ve\|_\infty)^{4\eta_2} ]^{1/2}
\mathbb{E} [|X_{t}^{\varepsilon}- X_{t(\delta)}^{\varepsilon} |^4]^{1/2}
\nn\\
&\qquad \qquad + 
\frac{c_4}{\ve} \mathbb{E} [|\hat{Y}_{t}^{\varepsilon} |^{q}]^{2r/q}
\mathbb{E} [|X_{t}^{\varepsilon}- X_{t(\delta)}^{\varepsilon} |^{2q/(q-2r)}
]^{(q-2r)/q}
\nn\\
&\le 
\frac{\gamma_2}{4\varepsilon} \mathbb{E} [
|Y_{t}^{\varepsilon}-\hat{Y}_{t}^{\varepsilon} |^2]
+\frac{c_5}{\ve} \delta^{2\beta}.
\nn
\end{align}
Note that we also used H\"older's inequality above.

Combining \eqref{eq.0606-1}--\eqref{eq.0606-2}, we see that
\[
\frac{d}{dt}
\mathbb{E} [|Y_{t}^{\varepsilon}-\hat{Y}_{t}^{\varepsilon}|^{2}]
\le 
- \frac{\gamma_2}{2\varepsilon} \mathbb{E} [
|Y_{t}^{\varepsilon}-\hat{Y}_{t}^{\varepsilon} |^2]
+\frac{c_6}{\ve} \delta^{2\beta},
\qquad
t\in [0,T].
\]
Applying the differential version of 
Gronwall's inequality (Lemma \ref{lem.Gron} (2)) 
with $f_0=0$,
$b = -\gamma_2/(2\ve)$ and $c=c_6\delta^{2\beta}/\ve$,
we have 
\[
\mathbb{E} [|Y_{t}^{\varepsilon}-\hat{Y}_{t}^{\varepsilon}|^{2}]
\le 
\frac{c_6\delta^{2\beta}}{\gamma_2} 
 \{1 -e^{-t\gamma_2/(2\ve)}\} 
     \le c_7\delta^{2\beta}, \qquad  t \in [0,T].
\]
This completes the proof of the lemma.
\end{proof}

%

It is easy to see that, if we define 
\begin{align}
M_t &=
           \int_{0}^{t} \{
             f (X_{s}^{\varepsilon}, Y_{s}^{\varepsilon} )
                -f (X_{s(\delta)}^{\varepsilon}, Y_{s}^{\varepsilon}) 
           \} d s
             +  \int_{0}^{t} \{f(X_{s(\delta)}^{\varepsilon}, 
                      Y_{s}^{\varepsilon})
              -f(X_{s(\delta)}^{\varepsilon}, \hat{Y}_{s}^{\varepsilon})
                \} d s
                       &\quad 
                       \label{def.230428}\\
                       &\qquad +
         \int_{0}^{t}  \{ 
            f(X_{s(\delta)}^{\varepsilon}, \hat{Y}_{s}^{\varepsilon})
              -\bar{f} (X_{s(\delta)}^{\varepsilon} )  \} d s
               +  
                 \int_{0}^{t} \{\bar{f} (X_{s(\delta)}^{\ve} )-\bar{f} (X_{s}^{\ve})) 
                 \} ds,
                 \nn
\end{align}
then 
\begin{align} 
X_{t}^{\varepsilon}-\bar{X}_{t} 
&=
\int_{0}^{t} f (X_{s}^{\varepsilon}, Y_{s}^{\varepsilon}) d s
   +\int_{0}^{t} \sigma (X_{s}^{\varepsilon} ) d B_{s} 
    - 
       \int_{0}^{t} \bar{f} (\bar{X}_{s}) d s
          -\int_{0}^{t} \sigma (\bar{X}_{s}) d B_{s}
          \label{eq.0603-1}\\ 
                   \nn\\
                   &= M_t + \Bigl( 
                       \int_{0}^{t} \{ \bar{f} (X_{s}^{\ve})
                           - \bar{f}(\bar{X}_{s})\}ds 
                                                + \int_{0}^{t} \{ \sigma (X_{s}^{\ve}) 
                               - \sigma (\bar{X}_{s}) \}d B_{s} \Bigr)
                            \nn
                                \end{align}
holds as an equality of CPs with respect to $B$.
We will apply Proposition \ref{prop.0703}  to \eqref{eq.0603-1}.

In the next lemma we estimate their H\"older norms.
\begin{lemma} \label{lem.0606-2}
Assume the same condition as in Lemma \ref{lem.0603-1}.
Then, there exists a positive constant $C$ independent of $\ve, \delta$
such that
\begin{align*}
\E \Bigl[
\Bigl\| \int_{0}^{\cdot} \{
             f (X_{s}^{\varepsilon}, Y_{s}^{\varepsilon} )
                -f (X_{s(\delta)}^{\varepsilon}, Y_{s}^{\varepsilon}) 
           \} d s \Bigr\|_{1}^2
           +
           \Bigl\| 
                \int_{0}^{\cdot} \{\bar{f} (X_{s(\delta)}^{\ve} )-\bar{f} (X_{s}^{\ve})) 
                 \} ds
                                \Bigr\|_{1}^2 
\\  +
\Bigl\| 
 \int_{0}^{\cdot} \{f(X_{s(\delta)}^{\varepsilon}, 
                      Y_{s}^{\varepsilon})
              -f(X_{s(\delta)}^{\varepsilon}, \hat{Y}_{s}^{\varepsilon})
                \} d s
 \Bigr\|_{1}^2
 \Bigr]   \le        
                          C\delta^{2\beta}.
\end{align*}
Here, $\|\cdot \|_{1}$ stands  for the $1$-H\"older (i.e. Lipschitz) norm.
\end{lemma}

\begin{proof} 
We use the globally Lipschitz property of $f$. 
The first term on the left hand side is dominated by 
 \begin{align*}
\E \Bigl[ \Bigl(
\int_{0}^{T} 
            | f (X_{s}^{\varepsilon}, Y_{s}^{\varepsilon} )
                -f (X_{s(\delta)}^{\varepsilon}, Y_{s}^{\varepsilon})| d s
                  \Bigr)^2 \Bigr]
                &\le 
                 T\int_{0}^{T} \E [
            | f (X_{s}^{\varepsilon}, Y_{s}^{\varepsilon} )
                -f (X_{s(\delta)}^{\varepsilon}, Y_{s}^{\varepsilon})|^2 ] ds
                                  \nn\\
                &\le
                 T^2 {\rm Lip} (f)^2 \sup_{0 \le s \le T} \E[
                   |X_{s}^{\varepsilon} - X_{s(\delta)}^{\varepsilon} |^2]
                   \nn\\
                   &\le 
T^2 {\rm Lip} (f)^2  \E[ \|X^{\varepsilon}\|_{\beta}^2] \delta^{2\beta}
 \le C \delta^{2\beta}.
                                                 \end{align*}
Here, we used \eqref{in.0602-1}. 
Similarly, the third term is dominated by 
\[
\E \Bigl[ \Bigl(
   \int_{0}^{T} | f(X_{s(\delta)}^{\varepsilon}, 
                      Y_{s}^{\varepsilon})
              -f(X_{s(\delta)}^{\varepsilon}, \hat{Y}_{s}^{\varepsilon})| d s
              \Bigr)^2 \Bigr] 
                \le 
               T^2 {\rm Lip} (f)^2\sup_{0 \le s \le T}
                  \E[ |Y_{s}^{\varepsilon} - \hat{Y}_{s}^{\varepsilon} |^2]           
                    \le C \delta^{2\beta}.
                     \]
                   Here, we used Lemma \ref{lem.0603-1}.

The estimate of the second term is more difficult 
since $\bar{f}$ may not be globally Lipschitz continuous. 
Since  $q \ge 2\eta_2$, we see from Proposition \ref{prop.av.drft}
that 
\[
| \bar{f} (X_{s(\delta)}^{\ve} ) -  \bar{f} (X_{s}^{\ve})|
 \le C' | X_{s(\delta)}^{\ve}- X_{s}^{\ve}| 
 (1+ |X_{s(\delta)}^{\ve} |^{\xi}+|X_{s}^{\ve}|^{\xi})
 \]
 for certain positive constants $C'$ and $\xi$ which do not depend on 
 $\ve, \delta, s$. (Below, the value of $C'$ may change from line to line.)
Then, the second term is dominated by
\begin{align*}
\E \Bigl[ \Bigl(
\int_{0}^{T} | \bar{f} (X_{s(\delta)}^{\ve} ) -  \bar{f} (X_{s}^{\ve})|  ds
                  \Bigr)^2 \Bigr]
                &\le  
                  T^2  C'       \sup_{0 \le s \le T} 
                    \E[ | X_{s(\delta)}^{\ve}- X_{s}^{\ve}|^2 
 (1+ |X_{s(\delta)}^{\ve} |^{\xi}+|X_{s}^{\ve}|^{\xi})^2
  ]      
   \nn\\ 
     &\le  T^2  C'    \E [\|X^{\varepsilon}\|_{\beta}^4]^{1/2}  \delta^{2\beta} 
       \sup_{0 \le s \le T} \E [ 1 + |X_{s}^{\ve}|^{4\xi}]  \le C' \delta^{2\beta}.\end{align*}
       Here, we used \eqref{in.0602-1} again. This completes the proof.
\end{proof}

%
\begin{lemma} \label{lem.0606-3}
Assume the same condition as in Lemma \ref{lem.0603-1}
and let $0< \gamma <1$.
Then, there exists a positive constant $C$ independent of $\ve, \delta$
such that
\begin{align*}
\E \Bigl[
           \Bigl\| 
                \int_{0}^{\cdot} \{f (X_{s(\delta)}^{\ve}, \hat{Y}_{s}^{\varepsilon})-\bar{f} (X_{s(\delta)}^{\ve}) 
                 \} ds
                                \Bigr\|_{\gamma}^2 
 \Bigr]   \le        
                          C(\delta^{2(1-\gamma)} 
+ \delta^{-2\gamma} \ve ).
\end{align*}
\end{lemma}

\begin{proof} 
We write $G_t =  \int_{0}^{t} \{f (X_{s(\delta)}^{\ve}, \hat{Y}_{s}^{\varepsilon})-\bar{f} (X_{s(\delta)}^{\ve}) 
                 \} ds$ 
                   and $G^1_{s,t} = G_t -G_s$ for simplicity.
It suffices to prove the case $1/2 <\gamma <1$.
Below, $c_i~(i \ge 1)$ are positive constants independent of 
$s, t, \delta, \ve, k$.

First, consider the case $0< t-s \le 2\delta$.
Since both $f$ and $\bar{f}$ are bounded, we can easily see that 
\begin{equation}\label{in.0607-1}
|G^1_{s,t}| \le (\|f\|_\infty + \|\bar{f}\|_\infty)(2\delta)^{1-\gamma} 
(t-s)^\gamma.
\end{equation}

Next, consider the case $t-s > 2\delta$.
In this case, there exists $k \in \N$ such that 
$[k\delta, (k+1)\delta] \subset [s,t]$. 
In the following estimates, 
we suppose that $s\delta^{-1}, t\delta^{-1} \notin \N$ 
(because the other cases are actually easier to be dealt with).
Using Schwarz's inequality and noting that the number of 
subintervals that are contained in $[s,t]$ does not exceed $(t-s)/\delta$,
we have
\begin{align} 
|G^1_{s,t}|^2 &\le
 \Bigl| G^1_{s, \, (\lfloor s/\delta\rfloor +1)\delta}
 +\sum_{k=\lfloor s/\delta\rfloor +1}^{\lfloor t/\delta\rfloor -1} 
 G^1_{k\delta,\, (k+1)\delta} +G^1_{\lfloor t/\delta\rfloor \delta,\, t}
 \Bigr|^2
\nn\\
&\le
c_1\delta^2 + 2 
\sum_{k=\lfloor s/\delta\rfloor +1}^{\lfloor t/\delta\rfloor -1} 1^2
\times
\sum_{k=\lfloor s/\delta\rfloor +1}^{\lfloor t/\delta\rfloor -1} 
 |G^1_{k\delta,\, (k+1)\delta} |^2
 \nn\\
&\le
c_1\delta^2 +\frac{2(t-s)}{\delta}
     \sum_{k=0}^{\lfloor T/\delta\rfloor -1} 
      |G^1_{k\delta,\, (k+1)\delta} |^2
                 \nn
 \end{align}
and therefore 
\begin{equation} \label{in.0608-1}
\frac{|G^1_{s,t}|^2}{(t-s)^{2\gamma}}
\le 
c_2 \delta^{2(1-\gamma)} + c_2 \delta^{-2\gamma}
  \sum_{k=0}^{\lfloor T/\delta\rfloor -1}  |G^1_{k\delta,\, (k+1)\delta} |^2.
                    \end{equation}

Combining \eqref{in.0607-1} and \eqref{in.0608-1}, we obtain that
\begin{align} 
\E [ \|G\|_{\gamma}^2 ] 
&\le 
c_2 \delta^{2(1-\gamma)} 
\label{in.0608-2}\\
&\quad
+ c_2T \delta^{-(1+2\gamma)} 
 \max_{0\le k \le \lfloor T/\delta\rfloor -1} 
 \E \Bigl[ \Bigl| 
            \int_{k\delta}^{(k+1)\delta} \{ f (X_{k\delta}^{\ve}, \hat{Y}_{s}^{\varepsilon})-\bar{f} (X_{k\delta}^{\ve})\} ds  
                 \Bigr|^2\Bigr].
                   \nn
 \end{align}
As will be seen in \eqref{in.0609-2} below, 
the expectation above is dominated by $c_3 \ve\delta$.
Hence, we have 
\[
\E [ \|G\|_{\gamma}^2 ] 
\le 
c_2 \delta^{2(1-\gamma)} 
+ c_4 \delta^{-2\gamma} \ve,
\]
which is the desired estimate.

Now, using a result of the frozen SDE in Appendix \ref{appen.FR},
we estimate the expectation on the right hand side 
of \eqref{in.0608-2}.
Let $\bar{w} =(\bar{w}_t)_{t \ge 0}$ be another standard $e$-dimensional 
BM which is independent of $(B, w)$
and let $(Y^{x,y}_t)_{t \ge 0}$ be a unique solution of 
the frozen SDE \eqref{def.frozenSDE} driven by $(\bar{w}_t)_{t \ge 0}$.
As will be shown in Proposition \ref{prop.estJ}, 
it holds that 
\begin{equation}\label{in.0609-1}
 \E \Bigl[ \Bigl| 
            \int_{0}^{t} \{ f (x, Y^{x,y}_{s})-\bar{f} (x)\} ds  
                 \Bigr|^2\Bigr]
                   \le C'  (1 + |x|^{\eta_3/2} + |y|) t,
                   \quad
                    (x, y)\in \R^{m+n}, \, t \ge 0                                     \end{equation}
                 for some constant $C' >0$ independent of $x, y, t$.
                   Recall that $\eta_3\ge 0$ is the constant in ${\bf (H6)}_q$.

We define  $\bar{Y}_t = Y^{x,y}_t \vert_{(x,y)= (X_{k\delta}^{\ve}, \hat{Y}_{k\delta}^{\varepsilon})}$.
(This means that $(X_{k\delta}^{\ve}, \hat{Y}_{k\delta}^{\varepsilon})$ is plugged into the superscript of $Y^{x,y}_t$.)
The starting point $(X_{k\delta}^{\ve}, \hat{Y}_{k\delta}^{\varepsilon})$ and $\bar{w}$ are independent. 
Then, $(\bar{Y}_{t/\ve})$ satisfies the following SDE:
\begin{align} 
\bar{Y}_{t/\ve}
&=
\hat{Y}_{k\delta}^{\varepsilon} 
+ \int_0^{t/\ve} g (X_{k\delta}^{\ve},   \bar{Y}_{s}) ds
+ \int_0^{t/\ve} h (X_{k\delta}^{\ve},   \bar{Y}_{s}) d^{{\rm I}}\bar{w}_s
\label{sde.0608}\\
&=
\hat{Y}_{k\delta}^{\varepsilon} 
+ \frac{1}{\ve}\int_0^{t} g (X_{k\delta}^{\ve},   \bar{Y}_{s/\ve}) ds
+ \frac{1}{\sqrt{\ve}}\int_0^{t} h (X_{k\delta}^{\ve},   \bar{Y}_{s/\ve}) 
d^{{\rm I}} (\sqrt{\ve}\bar{w}_{s/\ve}).
\nn
\end{align}
Note that $(\sqrt{\ve}\bar{w}_{s/\ve})_{t \ge 0}$ is 
again a  standard $e$-dimensional BM independent of 
$(X_{k\delta}^{\ve}, \hat{Y}_{k\delta}^{\varepsilon})$.
Compare SDE  \eqref{sde.0608} with SDE \eqref{def.0603}.
Then, we see that, restricted to the interval $[0,\delta)$,
$(\hat{Y}^\ve_{k\delta+ t})$ and $(\bar{Y}_{t/\ve})$
satisfy the same SDE with the same starting point 
(and the same frozen variable).
Hence, by the uniqueness of law of the frozen SDE, it holds that
\[
\bigl( X_{k\delta}^{\ve},    (\hat{Y}^\ve_{k\delta+ t})_{t \in [0,\delta)}\bigr)
=
\bigl( X_{k\delta}^{\ve},    (\bar{Y}_{t/\ve})_{t \in [0,\delta)}\bigr)
\qquad
\mbox{in law.}
\]
Hence, we have
\begin{align} 
\E \Bigl[ \Bigl| 
            \int_{k\delta}^{(k+1)\delta} \{ f (X_{k\delta}^{\ve}, \hat{Y}_{s}^{\varepsilon})-\bar{f} (X_{k\delta}^{\ve})\} ds  
                 \Bigr|^2\Bigr]
                 &=
                 \E \Bigl[ \Bigl| 
            \int_{0}^{\delta} \{ f (X_{k\delta}^{\ve}, \bar{Y}_{s/\ve})-\bar{f} (X_{k\delta}^{\ve})\} ds  
                 \Bigr|^2\Bigr]
                 \nn\\
                 &=
                 \ve^2  \E \Bigl[ \Bigl| 
            \int_{0}^{\delta/\ve} \{ f (X_{k\delta}^{\ve}, \bar{Y}_{s})-\bar{f} (X_{k\delta}^{\ve})\} ds  
                 \Bigr|^2\Bigr].
                 \nn
                                  \end{align}
Under the conditional expectation
$\E [\,\,\cdot ~|\sigma \{X_{k\delta}^{\ve},  \hat{Y}_{k\delta}^{\varepsilon}\}]$, 
$\bar{w}$ is still a standard BM
and $(X_{k\delta}^{\ve},  \hat{Y}_{k\delta}^{\varepsilon})$
can be viewed as constants, i.e. non-random. 
(Here, $\sigma \{X_{k\delta}^{\ve},  \hat{Y}_{k\delta}^{\varepsilon}\}$
stands for the sub-$\sigma$-field generated by 
$(X_{k\delta}^{\ve},  \hat{Y}_{k\delta}^{\varepsilon})$.)
Hence, we can use \eqref{in.0609-1} to see that
the right hand side above is dominated by
\[
\ve^2 \E [ 
 C'  (1 + |X_{k\delta}^{\ve}|^{\eta_3/2} + |\hat{Y}_{k\delta}^{\varepsilon}|) (\delta/\ve)]
 \le  c_3 \ve \delta,
\]
where Proposition \ref{prop.0601} and Lemma \ref{lem.0606-1} were used.
Thus, we have seen that 
\begin{equation} \label{in.0609-2}
\max_{0\le k \le \lfloor T/\delta\rfloor -1} 
 \E \Bigl[ \Bigl| 
            \int_{k\delta}^{(k+1)\delta} \{ f (X_{k\delta}^{\ve}, \hat{Y}_{s}^{\varepsilon})-\bar{f} (X_{k\delta}^{\ve})\} ds  
                 \Bigr|^2\Bigr]
                  \le c_3 \ve\delta,
                 \end{equation}
                    which completes the proof of the lemma.
\end{proof}

%

Now we prove our main theorem.
\begin{proof}[Proof of Theorem \ref{thm.main}]
Applying
Proposition \ref{prop.0703}  to \eqref{def.230428} and \eqref{eq.0603-1}, 
we obtain that
\begin{align} 
 \|X^{\ve} - \bar{X} \|_{\beta} &\le 
C \exp ( C \vertiii{B}_{\alpha}^\nu )  \,  \| M\|_{2\beta}
  \nn
 \end{align}
  for certain positive constant $C$ and $\nu$ which are 
  independent of $\ve, \delta$. (Below, $C$ 
     and $\nu$ may vary from line to line.)
By Lemmas \ref{lem.0606-2} and \ref{lem.0606-3}, we have
\[
\E [ \| M\|_{2\beta}^2 ] \le C (\delta^{2\beta} + \delta^{2(1-2\beta)} 
+ \delta^{-4\beta} \ve).
\]
Therefore, if we set $\delta := \ve^{1/(4\beta)} \log \ve^{-1}$ for example,
then $\| M\|_{2\beta}$ converges to $0$ in $L^2$-sense as $\ve\searrow 0$.
It immediately follows that $\|X^{\ve} - \bar{X} \|_{\beta}^p$ 
converges to $0$ in probability as $\ve\searrow 0$
for every $p \in [1,\infty)$.

On the other hand, 
we see from Proposition \ref{prop.0702} that
\[
\sup_{0<\ve \le 1}
\E
[ \|X^{\ve} - \bar{X}\|_{\beta}^p]  \le 
 2^{p-1}\sup_{0<\ve \le 1}
\E
[ \|X^{\ve}\|_{\beta}^p] + 2^{p-1}\E [\|\bar{X}\|_{\beta}^p]
\le C  (\E[\vertiii{B}_{\alpha}^{\nu p} ]+1) <\infty
\]
for every $p \in [1, \infty)$.
This implies that $\{ \|X^{\ve} - \bar{X}\|_{\beta}^p \}_{0<\ve \le 1}$
are uniformly integrable for each fixed $p$.
Hence, we have 
$\E [\|X^{\ve} - \bar{X}\|_{\beta}^p] \to 0$ as $\ve\searrow 0$.
This completes the proof of the main theorem.
\end{proof}

%

\appendix\section{Gronwall's inequality}

In this appendix
we recall two versions of Gronwall's inequality.
Note that in the differential version, $b$ can be negative.

\begin{lemma} \label{lem.Gron}
Let $T>0$.
\\
\noindent
{\rm (1)}~If a Lebesgue-integrable Borel measurable function
$f\colon [0,T] \to [0,\infty)$ satisfies for $a, b \in [0,\infty)$ 
that
\[
f_t \le a +b \int_0^t f_s ds, \qquad t \in  [0,T],
\]
then we have
\[
f_t \le  a e^{b t},
\qquad t \in  [0,T].
\]
\\
\noindent
{\rm (2)}~Let $f\colon [0,T] \to \R$ be an absolutely continuous 
function. Suppose that  for $b, c \in \R$ we have
\[
f^{\prime}_t \le b f_t + c, \qquad \mbox{for almost all $t \in [0,T]$}.
\]
Then, we have
\[
f_t \le (f_0 +\frac{c}{b}) e^{bt} - \frac{c}{b}, \qquad t \in  [0,T].
\]
\end{lemma}

\begin{proof} 
(1) is well-known. See Kusuoka \cite[p. 214]{ku} for example. 
We now prove (2). First, consider the case $c=0$.
For a.a. $t$, we have
$(f_t e^{-bt})^{\prime} = (f^{\prime}_t - b f_t)e^{-bt} \le 0$.
By integrating, we have
$f_t e^{-bt} \le f_0$ for every $t$.
If $c\neq 0$,
just set $g_t = f_t + (c/b)$. Then, $g^{\prime}_t \le b g_t$
and hence $g_t e^{-bt} \le g_0$ for every $t$.
This proves (2).
\end{proof}

%

\section{A simple application of It\^o's formula}

In this appendix 
we calculate $|Y_t|^q$, $q \in [2,\infty)$, for  a multi-dimensional 
It\^o process $(Y_t)$ by using It\^o's formula.

Set $F (y) = |y|^{q}= \{(y^1)^2 +\cdots+ (y^n)^2\}^{q/2}$ 
for $y\in \R^n$ and $q \in [2,\infty)$.
It is well-known $F\in C^2 (\R^n, \R)$.
First, we calculate its first and second order partial derivatives.
We write $\partial_i = \partial/  \partial y^i$ for simplicity.

\begin{lemma}
Let $F$ be as above. Then, for all $1 \le i, j \le n$ and 
$y\in \R^n$, we have the following:
\begin{align*}
\partial_i F (y) &= q |y|^{q-2} y^i, 
\\
\partial_i^2 F (y) &= q |y|^{q-2}  +q(q-2)|y|^{q-4} (y^i)^2,
\\
\partial_i \partial_j F (y) &= q (q-2)|y|^{q-4} y^iy^j
\qquad
\mbox{for $i\neq j$}.
\end{align*}
(Note that the right hand sides are well-defined for every 
$q\ge 2$ and $y \in \R^n$.
When $q>2$, $\partial_i \partial_j F (0)=0$ for all $(i, j)$. 
When $q=2$, $(q-2)$ times ``any quantity" is understood to be $0$.)
\end{lemma}

\begin{proof}
The proof is quite straightforward and  therefore is omitted.
\end{proof}

Let $(\Omega, \cF, {\mathbb P}, \{\cF_t \}_{0\le t \le T})$ be
a filtered probability space satisfying the usual condition 
and $(w_t)_{0\le t \le T}$ be an $e$-dimensional
standard $\{\cF_t \}$-BM defined on it.
Let $(G_t)_{0\le t \le T}$ and $(H_t)_{0\le t \le T}$ be 
progressively measurable processes which take values in 
$\R^n$ and $L(\R^e, \R^n)$, respectively. 
($L(\R^e, \R^n)$ is the set of $n\times e$ real matrices.)
We assume that 
\[
{\mathbb P} \Bigl(\int_0^T |G_s |ds <\infty \Bigr) =1
  = {\mathbb P} \Bigl(\int_0^T |H_s |^2 ds <\infty \Bigr) 
  \]
and define
\[
Y_t = Y_0 +\int_0^t G_s ds +\int_0^t H_s dw_s, \qquad t \in [0,T],
\]
where $Y_0$ is an $\cF_0$-measurable  $\R^n$-valued random variable.
Equivalently, in the component form, it reads
\[
Y_t^i = Y_0^i +\int_0^t G_s^i ds +\sum_{k=1}^e\int_0^t H_s^{ik} dw_s^k, 
  \qquad t \in [0,T], \,\, 1\le i \le n.
\]
One can easily see that
$
d \langle Y^{i},Y^{j} \rangle_t 
=
(H_t H_t^\top)^{ij} dt$ for all $i, j$. 
Here, $H^\top$ stands for the transpose of $H$
and $(HH^\top)^{ij}$ is the $(i,j)$-component of 
the matrix $HH^\top$.

\begin{lemma} \label{lem.apen.ito}
Let $(Y_t)_{0\le t \le T}$ be as above and $q\ge 2$.
Then, we almost surely have 
\begin{align*}
|Y_t|^q &= 
|Y_0|^q + 
\mbox{{\rm (a local martingale)}} 
\nn\\
&\quad +
q  \int_0^t   \{ 
|Y_s|^{q-2}   \langle Y_s, G_s\rangle 
+
\frac12 |Y_s|^{q-2} |H_s|^2
+
\frac{q-2}{2} |Y_s|^{q-4} (Y_s^\top H_sH_s^\top Y_s )
\}ds 
\nn\\
&\le |Y_0|^q + \mbox{{\rm (a local martingale)}} 
\nn\\
&\quad +
\frac{q}{2} \int_0^t   
|Y_s|^{q-2}  \{ 2\langle Y_s, G_s\rangle 
+
(q-1) |H_s|^2
\}ds, 
\qquad   t \in [0,T].
\end{align*}
Here, $|H_s|$ stands for the Hilbert-Schmidt norm of the matrix $H_s$. 
\end{lemma}

\begin{proof}
Since $0\le Y_s^\top H_sH_s^\top Y_s \le |H_s|^2 |Y_s|^2$,
the inequality is trivial.
Now we prove the equality by It\^o's formula.
\begin{align}
dF(Y_t)  
&= 
\sum_i \partial_i F (Y_t) dY_t^i 
   + \frac12 \sum_{i,j} \partial_i\partial_j F (Y_t) d \langle Y^{i},Y^{j} \rangle_t 
   \nn\\
     &= 
        \sum_i  q |Y_t|^{q-2} Y_t^i \sum_k H_t^{ik} dw_t^k 
          +
            \sum_i  q |Y_t|^{q-2} Y_t^i G_t^i dt 
             \nn\\
               &\quad + \frac12 \sum_i  q |Y_t|^{q-2}  \sum_k (H_t^{ik})^2 dt 
                 +
                 \frac12 \sum_{i,j}q(q-2) |Y_t|^{q-4} Y_t^i Y_t^j \sum_k H_t^{ik}H_t^{jk}dt.
                 \nn
\end{align}
The first term on the right hand side is clearly a local martingale.        
In the second term, $\sum_i Y_t^i G_t^i = \langle Y_s, G_s\rangle$.            
In the third term, $\sum_{i,k} (H_t^{ik})^2 = |H_t|^2$.
In the fourth term, $ \sum_{i,k,j} Y_t^i  H_t^{ik}H_t^{jk} Y_t^j =
Y_t^\top H_tH_t^\top Y_t$.             
          Noting these, we complete the proof.   
\end{proof}

%

\section{Frozen SDE}\label{appen.FR}

In this appendix, we recall basic facts on the frozen SDE 
associated with the slow-fast system  \eqref{def.SFeq}.
We basically follow \cite[Subsection 3.3]{lrsx}.
Let $(\Omega, \cF, {\mathbb P}, \{\cF_t \}_{t \ge 0})$ be
a filtered probability space satisfying the usual condition 
and $(w_t)_{t \ge 0}$ be an $e$-dimensional
standard $\{\cF_t \}$-BM defined on it.
(The time interval is $[0,\infty)$.
The setting in this appendix is not exactly the same as 
that in the main part of this paper.)
Since we are interested only in the law of 
the frozen SDE, the choice of $(\Omega, \cF, {\mathbb P}, \{\cF_t \}_{t \ge 0})$ that supports $(w_t)$ does not matter.
In this appendix, we simply write $dw_t$ for the standard It\^o integral
instead of $d^{{\rm I}}w_t$.

For any given $x\in \R^m$, we consider the following 
$\R^n$-valued SDE on $[0,\infty)$:
\begin{equation} \label{def.frozenSDE}
dY^{x,y}_t = g(x, Y^{x,y}_t)dt  +  h(x, Y^{x,y}_t) dw_t,
\qquad
Y^{x,y}_0 =y \in \R^n.
\end{equation}
We assume throughout this appendix
 that $h\in C(\R^m \times \R^n, L(\R^e, \R^n))$ is globally 
 Lipschitz continuous and $g \in C(\R^m \times \R^n, \R^n)$
 is locally Lipschitz continuous.
 
Under ${\bf (H6)}_2$,  SDE \eqref{def.frozenSDE}
has a unique strong (global) solution $(Y^{x,y}_t)_{t\ge 0}$
for every $x$ and $y$.
(The proof of this fact is essentially the same as the corresponding part 
 of Proposition \ref{prop.0601}.)
For every $x\in \R^m$, the family of processes 
$\{Y^{x,y}_{\cdot}\}_{y \in \R^n}$ indexed by $y$
becomes a diffusion process on $\R^n$.
The corresponding semigroup is denoted by $(P^x_t)_{t\ge 0}$, 
that is, 
\[
P^x_t \varphi (y) := {\mathbb E} [\varphi (Y^{x,y}_t) ],
\qquad t\ge 0,  \, y \in \R^n
\]
for every bounded Borel-measurable function on $\varphi$ on $\R^n$.
The next lemma (with $q=2$) and the standard Krylov-Bogoliubov argument
yield the existence of an invariant probability measure for 
$(P^x_t)_{t\ge 0}$ for every $x$.

\begin{lemma} \label{lem.lrsx3.6}
Assume ${\bf (H6)}_q$ for some $q \ge 2$. 
Then, there exists a positive  constants  $C_q$ 
(which is independent of $t, x, y$) such that
\[
{\mathbb E} [|Y^{x,y}_t|^q] \le e^{- tq\gamma_1 /4} |y|^q 
 +C_q  (1+|x|^{q\eta_3/2}),
\qquad t\ge 0, \, x\in \R^m, \, y \in \R^n.
\]
Here, $\gamma_1>0, \eta_3\ge 0$ are the constant in ${\bf (H6)}_q$.\end{lemma}

\begin{proof} 
This is a special case of \cite[Lemma 3.6]{lrsx}.
Alternatively, we can prove this lemma by slightly modifying
the proof of Proposition \ref{prop.0601}.
\end{proof}

\begin{lemma} \label{lem.lrsx3.7}
Assume ${\bf (H6)}_2$ and ${\bf (H7)}$. 
Then, it holds that
\[
{\mathbb E} [|Y^{x,y_1}_t - Y^{x,y_2}_t|^2] 
\le 
  e^{-\gamma_2 t} |y_1 -y_2|^2,
\qquad t\ge 0, \, x\in \R^m, \, y_1, y_2 \in \R^n.
\]
Here, $\gamma_2>0$ is the constant in ${\bf (H7)}$.
\end{lemma}

\begin{proof} 
This is a special case of \cite[Lemma 3.7]{lrsx}.
Alternatively, we can prove this lemma by slightly modifying
the proof of Proposition \ref{prop.0601}.
\end{proof}

%

\begin{lemma} \label{lem.lrsx3.8}
Assume ${\bf (H6)}_q$ for some $q \ge 2$ and ${\bf (H7)}$. 
Then, for every $x\in \R^m$, the semigroup $(P^x_t)_{t\ge 0}$
has a unique invariant probability measure $\mu^x$. 
Furthermore, the following two estimates hold: 
\\
{\rm (i)}~There exists a positive  constant $C_q$ 
(which is independent of $x$) such that
\[
\int_{\R^n} |z|^q \mu^x (dz) \le C_q  (1+|x|^{q\eta_3/2}),
\qquad  x\in \R^m.
\]
Here, $\eta_3\ge 0$ is the constant in ${\bf (H6)}_q$.
\\
{\rm (ii)}~There exists a positive  constant $C'$ such that 
for every $t\ge 0$, $x\in \R^m$, $y \in \R^n$ and Lipschitz function
$\varphi \colon \R^n\to \R$,
\[
\Bigl|
P^x_t \varphi (y)
-\int_{\R^n}  \varphi (z) \mu^x (dz)
\Bigr|
\le C' {\rm Lip}(\varphi) e^{-\gamma_2 t/2} (1+ |x|^{\eta_3/2} + |y|).
\]
Here,  ${\rm Lip}(\varphi)$ is the Lipschitz constant of $\varphi$
and $\gamma_2>0$ is the constant in ${\bf (H7)}$.
The constant $C'$ is independent of $t, x, y, \varphi$.
\end{lemma}

\begin{proof} 
This is a special case of \cite[Proposition 3.8]{lrsx}.
\end{proof}

\begin{lemma} \label{lem.lrsx3.10}
Assume ${\bf (H5)}_r$  
and ${\bf (H6)}_{2(r \vee 1)}$ for some $r\ge 0$.
Then, there exist positive  constants $C$ and $\xi$
independent of $t, x_1, x_2, y$ such that
\[
{\mathbb E} [|Y^{x_1,y}_t - Y^{x_2,y}_t|^2] 
\le
C |x_1- x_2|^2 (1+  |x_1|^{2\xi}+|x_2|^{2\xi} +|y|^{2(r\vee 1)}) 
\]
for all $t \ge 0$, $x_1, x_2 \in \R^m$ and  $y \in \R^n$.
\end{lemma}

\begin{proof} 
This is a special case of \cite[Lemma 3.10]{lrsx}.
\end{proof}

%
Under the assumptions of Lemma \ref{lem.lrsx3.8},  we set
\[
\bar{f} (x) =\int_{\R^n}   f(x,y) \mu^x (dy), \qquad  x\in\R^m
\]
for $f \colon \R^m \times \R^n\to \R^m$ whenever the integral  
on the right hand side is well-defined.

\begin{proposition} \label{prop.av.drft}
We assume ${\bf (H2)}$, ${\bf (H5)}_r$
and ${\bf (H6)}_{2(r \vee 1)}$ for some $r\ge 0$.
Then, the map $\bar{f}\colon \R^m \to  \R^m$ defined as above 
is bounded and locally Lipschitz continuous.
More precisely, we have $\|\bar{f} \|_\infty \le \|f\|_\infty$ 
and 
\[
| \bar{f} (x_1) -  \bar{f} (x_2)| \le C | x_1- x_2| 
 (1+ | x_1|^{\xi}+|x_2|^{\xi}), 
\qquad  x_1, x_2 \in \R^m
 \]
 for some constant $C >0$ independent of $x_1, x_2$.
 Here, $\xi >0$ is the constant which 
 appears in Lemma \ref{lem.lrsx3.10}.
 \end{proposition}

\begin{proof} 
Clearly, $\|\bar{f} \|_\infty \le \|f\|_\infty$ holds.
We now show the local Lipschitz property.
For every Lipschitz function $\varphi \colon \R^n\to \R$, 
we see from Lemmas \ref{lem.lrsx3.8} and \ref{lem.lrsx3.10} that
\begin{align*}
\Bigl|
\int_{\R^n}  \varphi (y) \mu^{x_1} (dy) - \int_{\R^n}  \varphi (y) \mu^{x_2} (dy)\Bigr|
&\le
\lim_{t\to\infty} | P^{x_1}_t  \varphi (0)- P^{x_2}_t  \varphi (0)|
\nn\\
&\le
\lim_{t\to\infty} | {\mathbb E} [\varphi (Y^{x_1, 0}_t) ] -
   {\mathbb E} [\varphi (Y^{x_2, 0}_t) ]|
   \nn\\
&\le
{\rm Lip} (\varphi) \limsup_{t\to\infty}
   {\mathbb E} [|Y^{x_1, 0}_t -Y^{x_2, 0}_t| ]
   \nn\\
    &\le 
       {\rm Lip} (\varphi)C |x_1- x_2| (1+  |x_1|^{\xi}+|x_2|^{\xi}).
                 \end{align*}
                 Then, we have
                 \begin{align*}
                            | \bar{f} (x_1) -  \bar{f} (x_2)| 
                            &\le    
                            \int_{\R^n}   |f(x_1,y)- f(x_2,y)| \mu^{x_1} (dy)
                            \nn\\
                             &\qquad + 
                            \Bigl|
\int_{\R^n}  f (x_2,y) \mu^{x_1} (dy) - \int_{\R^n}  f(x_2,y) \mu^{x_2} (dy)\Bigr|
 \nn\\ 
   &\le {\rm Lip} (f) |x_1-x_2| +{\rm Lip} (f) 
      C |x_1- x_2| (1+  |x_1|^{\xi}+|x_2|^{\xi}).                                                           \end{align*}
     This completes the proof of the proposition.
                 \end{proof}

\begin{proposition} \label{prop.estJ}
Assume ${\bf (H2)}$,  ${\bf (H7)}$ and ${\bf (H6)}_{q}$ for some $q\ge 2$.  Then, we have
 \[
 {\mathbb E} \Bigl[ \Bigl| 
            \int_{0}^{t} \{ f (x, Y^{x,y}_{s})-\bar{f} (x)\} ds  
                 \Bigr|^2\Bigr]
                   \le C  (1 + |x|^{\eta_3/2} + |y|) t,
                   \quad
                    (x, y)\in \R^{m+n}, \, t \ge 0   
                    \]
                 for some constant $C>0$ independent of $x, y, t$.
                  Here, $\eta_3\ge 0$ is the constant in ${\bf (H6)}_q$.
\end{proposition}

  \begin{proof}
   The expectation with respect to the law of $Y^{x,y}$ 
  is denoted by $\hat{\mathbb E}^{x,y}$.
  For  $s \le u$, we see from the Markov property that
 \begin{align}  
 \lefteqn{
|{\mathbb E}[\langle f (x, Y^{x,y}_{u})-\bar{f} (x), 
  f (x, Y^{x,y}_{s})-\bar{f} (x) \rangle ]|
  }
  \nn\\
&=
  |\hat{\mathbb E}^{x,y} 
  [\langle f (x, Y_{u})-\bar{f} (x),  f (x, Y_{s})-\bar{f} (x) \rangle ]|
     \nn\\
     &=  
       | \hat{\mathbb E}^{x,y} \bigl[
             \langle 
                \hat{\mathbb E}^{x, Y_s} [f (x, Y_{u-s})-\bar{f} (x)],  
    f (x, Y_{s})-\bar{f} (x) \rangle \bigr] |
     \nn\\
     &\le    
       2\|f\|_{\infty} \hat{\mathbb E}^{x,y} \bigl[
             \bigl|
                \hat{\mathbb E}^{x, Y_s} [f (x, Y_{u-s})-\bar{f} (x) ] \bigr|\bigr]
                   \nn\\
     &\le    
       2\|f\|_{\infty}        \hat{\mathbb E}^{x,y} \bigl[
           C' {\rm Lip} (f) e^{-\gamma_2 (u-s)/2} (1+ |x|^{\eta_3/2} + |Y_s|)           \bigr]              
                \nn\\
                 &\le    
       2C' \|f\|_{\infty}{\rm Lip} (f)  e^{-\gamma_2 (u-s)/2}
             {\mathbb E}[1+ |x|^{\eta_3/2} + |Y^{x,y}_s|  ]
              \nn\\
              &\le
                  C_1 
                   e^{-\gamma_2 (u-s)/2}
             (1+ |x|^{\eta_3/2} +  |y|)
              \nn
                            \end{align}
                            for some constant $C_1 >0$ which is independent of $x, y, t$.
                       Here, we used Lemma \ref{lem.lrsx3.8} {\rm (ii)} 
                       to the third to the last inequality 
                       and Lemma \ref{lem.lrsx3.6}  to the last inequality.

A simple application of Fubini's theorem yields that   
 \begin{align*}
 \lefteqn{
  {\mathbb E} \Bigl[ \Bigl| 
           \int_{0}^{t} \{ f (x, Y^{x,y}_{s})-\bar{f} (x)\} ds  
                 \Bigr|^2\Bigr]
                 }
                  \\
                 &=
                  2\int_0^t  ds \int_s^t  du
                   {\mathbb E}[
                     \langle f (x, Y^{x,y}_{u})-\bar{f} (x), 
  f (x, Y^{x,y}_{s})-\bar{f} (x) \rangle ]
     \\
                 &\le
                   C_2   (1+ |x|^{\eta_3/2} +  |y|)
                     \int_0^t  ds \int_s^t  du \, e^{-\gamma_2 (u-s)/2}    
                     \\
                      &\le C_3       (1+ |x|^{\eta_3/2} +  |y|)t                                                                    
\end{align*}
 for some constant $C_3 >0$ which is independent of $x, y, t$.
    \end{proof}

%

%

\bigskip
\begin{flushleft}
  \begin{tabular}{ll}
    Yuzuru \textsc{Inahama}
    \\
    Faculty of Mathematics,
    \\
    Kyushu University,
    \\
    744 Motooka, Nishi-ku, Fukuoka, 819-0395, JAPAN.
    \\
    Email: {\tt inahama@math.kyushu-u.ac.jp}
  \end{tabular}
\end{flushleft}

\newpage

\begin{center}
{\Large {\bf Erratum}}
\end{center}

\bigskip

Lemma 5.3 is wrong. 
In the proof, the $1$-variation norms are estimated instead 
of the $1$-H\"older norms.
\footnote{
The published version of this paper (Tohoku mathematical Journal, 2025)
contains the same error.}

Fortunately, we can fix this error as in what follows.
This modification affects no other parts of this paper.

\medskip

Lemma 5.3 should be replaced by the following:
\begin{lemma} \label{lem.5.3_modified}
Assume the same condition as in Lemma \ref{lem.0603-1}
and let $1/2 <\gamma <1$.
Then, there exists a constant $C=C_\gamma >0$ independent of $\ve, \delta$
such that
\begin{align*}
\E \Bigl[
\Bigl\| \int_{0}^{\cdot} \{
             f (X_{s}^{\varepsilon}, Y_{s}^{\varepsilon} )
                -f (X_{s(\delta)}^{\varepsilon}, Y_{s}^{\varepsilon}) 
           \} d s \Bigr\|_{1}^2
           +
           \Bigl\| 
                \int_{0}^{\cdot} \{\bar{f} (X_{s(\delta)}^{\ve} )-\bar{f} (X_{s}^{\ve})) 
                 \} ds
                                \Bigr\|_{1}^2 
                               \Bigr]
                                \le C\delta^{2\beta},
\\ 
 \E \Bigl[
\Bigl\| 
 \int_{0}^{\cdot} \{f(X_{s(\delta)}^{\varepsilon}, 
                      Y_{s}^{\varepsilon})
              -f(X_{s(\delta)}^{\varepsilon}, \hat{Y}_{s}^{\varepsilon})
                \} d s
 \Bigr\|_{\gamma}^2
 \Bigr]   \le        
                          C\delta^{2\beta (1-\gamma)}.
\end{align*}
Here, $\|\cdot \|_{\gamma}$ stands  for the $\gamma$-H\"older (i.e. Lipschitz) norm.
\end{lemma}

\begin{proof} 
In this proof, an unimportant positive constant is denoted by $C$,
which may vary from line to line.
Generally, $\| \int_{0}^{\cdot} A_s ds\|_1 \le  \sup_{0\le s\le T} |A_s|$ holds.

Firstly, consider the case
$A_s := f (X_{s}^{\varepsilon}, Y_{s}^{\varepsilon} )
                -f (X_{s(\delta)}^{\varepsilon}, Y_{s}^{\varepsilon})$.
Then, by the Lipschitz continuity of $f$, it holds that
$\sup_{0\le s\le T} |A_s| 
 \le C | X_{s}^{\varepsilon}- X_{s(\delta)}^{\varepsilon}|
  \le C \delta^{\beta} \| X^{\varepsilon}\|_\beta.$

Secondly, consider the case
$A_s :=\bar{f} (X_{s(\delta)}^{\ve} )-\bar{f} (X_{s}^{\ve}))$.
This term is more difficult 
since $\bar{f}$ may not be globally Lipschitz continuous. 
Since  $q \ge 2\eta_2$, we see from Proposition \ref{prop.av.drft}
that 
\[
|A_s|=
| \bar{f} (X_{s(\delta)}^{\ve} ) -  \bar{f} (X_{s}^{\ve})|
 \le C' | X_{s(\delta)}^{\ve}- X_{s}^{\ve}| 
 (1+ |X_{s(\delta)}^{\ve} |^{\xi}+|X_{s}^{\ve}|^{\xi})
 \]
 for certain positive constants $C$ and $\xi$ which do not depend on 
 $\ve, \delta, s$. 
Hence, 
\[
\sup_{0\le s\le T} |A_s| 
\le 
 C \delta^{\beta} \| X^{\varepsilon}\|_\beta (1+ \| X^{\varepsilon}\|_\beta^{\xi}).
\]
By \eqref{in.0602-1} in Proposition \ref{prop.0601}, we can easily prove the first assertion.

The most difficult is the case $A_s :=f(X_{s(\delta)}^{\varepsilon}, 
                      Y_{s}^{\varepsilon})
              -f(X_{s(\delta)}^{\varepsilon}, \hat{Y}_{s}^{\varepsilon})$.
Noting that $|A_s|$ is bounded, we have 
\[
\frac{|A_t - A_s|}{(t-s)^\gamma} 
\le 
(t-s)^{-\gamma} \Bigl( \int_s^t |A_u| du \Bigr)^{1-\gamma}  \,
\Bigl( \int_s^t  2 \|f\|_\infty du \Bigr)^{\gamma}
\le
(2 \|f\|_\infty)^\gamma \Bigl( \int_0^T |A_u| du \Bigr)^{1-\gamma}
\]
by H\"older's inequality. 
Note that the right hand side is independent of $(s,t)$.
Using H\"older's inequality twice, we have
\begin{align*}
 \E \Bigl[
\Bigl\| 
 \int_{0}^{\cdot} A_s d s
 \Bigr\|_{\gamma}^2
 \Bigr]
 &\le 
 C \Bigl( \int_0^T \E [ |A_u|^2 ] du \Bigr)^{1-\gamma}
 \\
 &\le 
 C \Bigl( \sup_{0 \le u \le T} \E [ |A_u|^2 ] \Bigr)^{1-\gamma}
 \\
 &\le 
 C \Bigl( \sup_{0 \le u \le T} \E [ |Y_{u}^{\varepsilon}-\hat{Y}_{u}^{\varepsilon} |^2 ] \Bigr)^{1-\gamma}
 \le 
 C\delta^{2\beta (1-\gamma)},
 \end{align*}
 where we used Lemma \ref{lem.0603-1} and
 the Lipschitz continuity of $f$.
 This proves the second assertion.
\end{proof}

\bigskip

Due to the above modification of Lemma 5.3, 
the estimate of $\E [ \| M\|_{2\beta}^2 ]$ in Page 36, Line 5 needs 
a slight adjustment: 
\[
\E [ \| M\|_{2\beta}^2 ] \le 
C (\delta^{2\beta} + \delta^{2\beta (1-2\beta)}  +\delta^{2(1-2\beta)} 
+ \delta^{-4\beta} \ve).
\]
Note that $\delta^{2\beta (1-2\beta)}$ was newly inserted.
Since this estimate still implies $\E [ \| M\|_{2\beta}^2 ] \to 0$ as $\ve\searrow 0$ (if we choose $\delta := \ve^{1/(4\beta)} \log \ve^{-1}$ for example),
no other parts of this paper needs correcting.


\begin{thebibliography}{00}


\bibitem{bsyy}
J. Bao, Q. Song, G. Yin and C. Yuan, 
Ergodicity and strong limit results for two-time-scale functional stochastic differential equations,
Stoch. Anal. Appl. 35 (2017), no. 6, 1030--1046. 

\bibitem{bgs1}
S. Bourguin, S. Gailus and K. Spiliopoulos, 
Typical dynamics and fluctuation analysis of slow-fast systems driven by fractional Brownian motion,
Stoch. Dyn. 21 (2021), no. 7, Paper No. 2150030, 30 pp. 

\bibitem{bgs2}
S. Bourguin, S. Gailus and K. Spiliopoulos, 
Discrete-time inference for slow-fast systems driven by fractional Brownian motion,
Multiscale Model. Simul. 19 (2021), no. 3, 1333--1366. 

\bibitem{bgs3}
S. Bourguin, T. Dang and K. Spiliopoulos, 
Moderate deviation principle for multiscale systems driven by fractional Brownian motion,
Preprint (2022). arXiv:2206.06794.


\bibitem{cf}
S. Cerrai and M. Freidlin, 
Averaging principle for a class of stochastic reaction-diffusion equations,
Probab. Theory Related Fields 144 (2009), no. 1-2, 137--177. 

\bibitem{duc}
L. H. Duc,  
Controlled differential equations as rough integrals,
Preprint (2020). arXiv: 2007.06295.

\bibitem{fw}
M. Freidlin and  A. D. Wentzell, 
Random perturbations of dynamical systems,
Translated from the 1979 Russian original by Joseph Sz\"ucs,
 Third edition, Grundlehren der mathematischen Wissenschaften 260, Springer, Heidelberg, 2012. 


\bibitem{fh}
P. Friz and M. Hairer,
A course on rough paths, Springer, Cham, 2014. 

\bibitem{fv}
P. Friz and N. Victoir, 
Multidimensional stochastic processes as rough paths, 
Cambridge University Press, Cambridge, 2010.


\bibitem{giv}
D. Givon, 
Strong convergence rate for two-time-scale jump-diffusion stochastic differential systems,
Multiscale Model. Simul. 6 (2007), no. 2, 577--594. 


\bibitem{gkk}
D. Givon,  I. G. Kevrekidis and R. Kupferman, 
Strong convergence of projective integration schemes for singularly perturbed stochastic differential systems,
Commun. Math. Sci. 4 (2006), no. 4, 707--729. 

\bibitem{gol}
J. Golec,
Stochastic averaging principle for systems with pathwise uniqueness, 
Stochastic Anal. Appl. 13 (1995), no. 3, 307--322. 

\bibitem{gl}
J. Golec and  G. Ladde, 
Averaging principle and systems of singularly perturbed stochastic differential equations,
J. Math. Phys. 31 (1990), no. 5, 1116--1123. 


\bibitem{hl}
M. Hairer and X. Li, 
Averaging dynamics driven by fractional Brownian motion,
Ann. Probab. 48 (2020), no. 4, 1826--1860. 


\bibitem{hxpw}
M. Han, Y. Xu, B. Pei and J.-L. Wu, 
Two-time-scale stochastic differential delay equations driven by multiplicative fractional Brownian noise: averaging principle,
J. Math. Anal. Appl. 510 (2022), no. 2, Paper No. 126004, 31 pp. 

\bibitem{hy}
 J. Hu and C. Yuan, 
 Strong convergence of neutral stochastic functional differential equations with two time-scales,
 Discrete Contin. Dyn. Syst. Ser. B 24 (2019), no. 11, 5831--5848. 
 
 
\bibitem{khas}
R. Z. Khas'minski\u{\i}, 
On the principle of averaging the It\^o's stochastic differential equations,
Kybernetika 4 (1968), 260--279. 


\bibitem{ku}
S. Kusuoka, 
Stochastic analysis, Springer, Singapore, 2020. 



\bibitem{lixm}
X. Li, 
Perturbation of conservation laws and averaging on manifolds, Computation and combinatorics in dynamics, stochastics and control, 499--550, 
Abel Symp. 13, Springer, Cham, 2018. 

\bibitem{ls}
X. Li and J. Sieber, 
Slow-fast systems with fractional environment and dynamics,
Preprint (2022). arXiv:2012.01910.

\bibitem{liu}
 D. Liu,  
 Strong convergence of principle of averaging for multiscale stochastic dynamical systems,
 Commun. Math. Sci. 8 (2010), no. 4, 999--1020.
 
\bibitem{liu2}
D. Liu,  
Strong convergence rate of principle of averaging for jump-diffusion processes,
Front. Math. China 7 (2012), no. 2, 305--320. 

 
\bibitem{lrsx}
W. Liu, M. R\"ockner, X. Sun and Y. Xie, 
Averaging principle for slow-fast stochastic differential equations with time dependent locally Lipschitz coefficients, 
J. Differential Equations 268 (2020), no. 6, 2910--2948. 


\bibitem{pix1}
B. Pei, Y. Inahama and Y. Xu, 
Averaging Principles for Mixed Fast-Slow Systems Driven by Fractional Brownian Motion,
To appear in Kyoto J. Math. (2023).  arXiv: 2001.06945.

\bibitem{pix2}
B. Pei, Y. Inahama and Y. Xu,  
Averaging principle for fast-slow system driven by mixed fractional Brownian rough path,
J. Differential Equations 301 (2021), 202--235. 

\bibitem{rsx1}
M. R\"ockner, X. Sun and Y. Xie, 
Strong and weak convergence in the averaging principle for SDEs with H\"older coefficients,
Preprint (2019).  arXiv: 1907.09256.


\bibitem{rsx2}
M. R\"ockner, X. Sun and Y. Xie, 
Strong convergence order for slow-fast McKean-Vlasov stochastic differential equations,
Ann. Inst. Henri Poincar\'e Probab. Stat. 57 (2021), no. 1, 547--576.  
 
\bibitem{sxx}
X. Sun, L. Xie and Y. Xie, 
Strong and weak convergence rates for slow-fast stochastic differential equations driven by $\alpha$-stable process,
Preprint (2021).  arXiv: 2004.02595.
  
\bibitem{ve}
A. Yu. Veretennikov, 
On an averaging principle for systems 
of stochastic differential equations,
Math. USSR-Sb. 69 (1991), no. 1, 271--284.

\bibitem{wxy}
R. Wang, Y. Xu and H. Yue,  
Stochastic averaging for the non-autonomous mixed stochastic differential equations with locally Lipschitz coefficients,
Statist. Probab. Lett. 182 (2022), Paper No. 109294, 11 pp. 

\bibitem{wy1}
F. Wu and G. Yin, 
An averaging principle for two-time-scale stochastic functional differential equations,
J. Differential Equations 269 (2020), no. 1, 1037--1077. 

\bibitem{wy2}
F. Wu and G. Yin,
Fast-slow-coupled stochastic functional differential equations,
J. Differential Equations 323 (2022), 1--37. 

\bibitem{xlm}
J. Xu, J. Liu and Y. Miao, 
Strong averaging principle for two-time-scale SDEs with non-Lipschitz coefficients,
J. Math. Anal. Appl. 468 (2018), no. 1, 116--140. 



\bibitem{xllm}
J. Xu, J. Liu, J. Liu and  Y. Miao, 
Strong averaging principle for two-time-scale stochastic McKean-Vlasov equations,
Appl. Math. Optim. 84 (2021), suppl. 1, S837--S867. 


\bibitem{xm}
J. Xu and Y. Miao, 
$L^p\,(p>2)$-strong convergence of an averaging principle for two-time-scales jump-diffusion stochastic differential equations,
Nonlinear Anal. Hybrid Syst. 18 (2015), 33--47. 


\bibitem{zfwl}
B. Zhang, H. Fu, L. Wan and J. Liu, 
Weak order in averaging principle for stochastic differential equations with jumps,
Adv. Difference Equ. 2018, Paper No. 197, 20 pp. 
\end{thebibliography}
\end{document}